\documentclass[12pt]{amsart}
\usepackage{amsthm}
\usepackage{epsfig}
\usepackage{amssymb}
\usepackage{latexsym}
\usepackage{amsmath}
\newtheorem{theorem}{Theorem} 
\newtheorem{definition}[theorem]{Definition}

\newtheorem{corollary}[theorem]{Corollary}
\newtheorem{lemma}[theorem]{Lemma}

\newtheorem{proposition}[theorem]{Proposition}
\newtheorem{remark}[theorem]{Remark}

\def\tg{\tilde{\gamma}}
\def \vp{\varphi}
\def\Se{S^{(8)}}
\def \tZZ{\tilde{Z}}

\def\R{{\Bbb R}}

\def\C{{\Bbb C}}

\def\Z{{\Bbb Z}}
\def\P{{\Bbb P}}

\def \ta{\tau}
\def \ta1{\tau_1}

\def \dl{\delta}
\def \g{\gamma}
\def \G{\Gamma}

\def \vp{\varphi}

\def \la{\langle}
\def \ra{\rangle}

\def \ra{\rangle}

\def \CP{\Bbb C \Bbb P}

\def \tC{\tilde{C}}

\def \tP{\tilde{\pi}}

\def \prodl{\prod\limits}
\def\hF{\hat{F}}
\def\uZ{\underline{Z}}
\def\bZ{\bar{Z}}
\def\bZut{\bar{Z}^2}
\def\uZut{\underline{Z}^2}
\def\uZumt{\underline{Z}^{-2}}

\def\eg{{\it{e.g.}}}

\def \zovera {
    \mathop{\lower 10pt \hbox{${\buildrel{\displaystyle\bar{z}} \over {\scriptstyle{(a)}}} $}}
    {\lower 4pt \hbox{${\scriptstyle{ij}}$}} 
} 

\def\vK{{van Kampen}}
\def\sub{{\subseteq}}

\newcommand\set[1]{{\{{#1}\}}}
\def\suchthat{{\,:\,}}

\def\wrt{{with respect to }}

\def\st{{such that }}
\def\pitil{{\tilde{\pi}_1}}

\newcommand\begintable[1][] {{}}

\newif\ifXY 
\XYtrue     
\newif\ifbigmatrices
\bigmatricestrue 

 \ifXY
 \usepackage{xy}
 \fi
 \ifXY
 \xyoption{all}
 \fi

\begin{document}

 \title {On the degeneration, regeneration and braid monodromy of $T \times T$}

 \author{ Meirav Amram \and Mina Teicher}



 \renewcommand{\subjclassname}{%
       \textup{2000} Mathematics Subject Classification}


\date{\today}

\maketitle

\begin{abstract}
This paper is the first  in a series of three papers concerning the surface 
 $T \times T$.  In this paper we study the degeneration of $T \times T$ and the regeneration of its degenerated object.  We study the braid monodromy  and its regeneration.
\end{abstract}

\section{Introduction}\label{intro}
Algebraic surfaces are classified by discrete and continuous
invariants.
Fixing the discrete invariants gives a family of algebraic surfaces
parameterized by algebraic variety called the {\bf moduli space}.
All surfaces in the same moduli space have the same homotopy type and
therefore the same fundamental group.
So, fundamental groups are discrete invariants of the surfaces 
and form a central tool in their classification.

We have no algorithm to compute the fundamental group of a given
algebraic surface $X$.
But we can cover $X$ by a surface $X_{Gal}$ with a computable
fundamental group.

The given surface $X$ is projective.
We embed it in $\C\P^N$.
The projection of $X$ from a ``general'' point in $\C\P^N-X$ will
map $X$ onto a surface in $\C\P^{N-1}$ which we project again into
$\C\P^{N-2}$, etc., until we finally project $X$ onto
$\C\P^2$.
This gives a ``generic'' projection $f\colon X\to\C\P^2$.
The ramification locus of $f$ is a curve $S$ in $\C\P^2$,
called the branch curve.
The singular points of $S$ are ordinary nodes or cusps.

Put $n=\deg(f)$.
Consider the fibred product
$$
X\times_f\cdots\times_fX
=\{(x_1,\ldots,x_n)\mid x_i\in X,\
f(x_i)=f(x_j)\ \forall i,j\}
$$
and the diagonal
$$
\Delta =\{(x_1,\ldots,x_n)\in X\times_f\cdots\times_f X
\mid x_i=x_j
\hbox{ for some } i\ne j\}.
$$
The surface $X_{Gal}$ is the Galois cover of $X$ with respect to $f$.
That is, the Zariski closure of the complement of $\Delta$:
$$
X_{Gal}=\overline{X\times_f\cdots\times_fX-\Delta}.
$$

In order to compute $\pi_1(X_{Gal})$ we must first compute the
fundamental groups
$\pi_1(\C^2-S,*)$ and $\pi_1(\C\P^2-S,*)$,
using tools which Moishezon-Teicher developed and applied in [MoTe4],
[MoTe5], and [MoTe8].

\begin{enumerate}
\item Degeneration of surfaces to unions of planes.
\item 
The braid group and the braid monodromy algorithm.
\item 
 Regeneration of degenerated objects and regeneration of braids.
\item 
The van Kampen Method for cuspidal curve.
\end{enumerate}

Having degenerated $X$ into a union $X_0$ of planes, we project $X_0$
onto $\C\P^2$ simultaneously with the projection of $X$.
The branch curve now is a line arrangement, that is, a union of
lines.
We call the intersection points of these lines singular points.
We use the braid monodromy to compute braids for the singular
points.
Now we regenerate the branch curve $S$ from the line arrangement
$S_0$.
This regenerates the braids.
Each regenerated braid gives, by van Kampen, a relation in canonical
generators $\Gamma_j$, $\Gamma_{j'}$ of $\pi_1(\C^2-S,*)$.

Now we project $\C\P^2-S$ into $\C\P^1$ from a ``general''
point in $\C\P^2$.
The fibre over a ``good'' point in $\C\P^1$ is a line which cuts
$S$ in $n$ points.
The group $\pi_1(\C^2-S,*)$ acts on these points.
This leads to a permutation representation
$\psi\colon\pi_1(\C^2-S,*)\to S_n$.
The images of $\Gamma_j$ and $\Gamma_{j'}$
under $\psi$ are transpositions.
Let $N_j$ be the normal subgroup of $\pi_1(\C^2-S,*)$ generated by
$\Gamma_j^2$ and $\Gamma_{j'}^2$.
This leads to a short exact sequence:
\begin{equation}
1\longrightarrow \frac{ \ker\psi}{N_j}\longrightarrow
\frac{ \pi_1(\C^2-S,*)}{N_j}
\longrightarrow S_n \longrightarrow 1.
\end{equation}
Moishezin-Teicher [Am2, Theorem 1.76] then say:
$$
\pi_1(X_{Gal}^{Aff})\cong \frac{ ker\psi}{N_j},
$$
where $X_{Gal}^{Aff}$ denotes the affine part of $X_{Gal}$.

Finally we use the Reidemeister-Schreier Method to compute the kernel
in (1).
The group $\pi_1(X_{Gal})$ has one more relation, namely, the product
of all generators equals $1$.

\section{The surface $T \times T$}\label{surface}

Let $T_1$ be a complex torus, \eg\ $T_1 = \set{(x,y)\suchthat y^2=x^3-6x} \sub \C^2$.
In order to embed $T_1 \times T_1$ in a projective space, we first transfer to homogeneous coordinates by substituting $\frac{x}{z},\frac{y}{z}$ into the equation. Then we get the projective torus $T = \set{[x,y,z] \suchthat y^2z=x^3-6xz^2} \sub \CP^2$.

Now $T$ is a curve of dimension $1$ and degree $3$ in $\CP^2$. 

The Segre map $T\times T \rightarrow \CP^8$ is defined by
$$[a_1, a_2, a_3]\ast [b_1, b_2, b_3] \longmapsto [a_1 b_1, a_2b_1, a_3b_1, a_1b_2, a_2b_2, a_3b_2, a_1b_3, a_2b_3,a_3b_3].$$
This embeds $T \times T$ as a surface of dimension $2$ and degree $18$ in $\C \P^8$(see Section \ref{degen}).

We remark that the general form of a projective torus is
$$\set{[x,y,z] \suchthat y^2z = Ax^3 + Bxz^2 + Cz^3},$$
see \cite{Sha1} for details.

\section{Degeneration of $T \times T$}\label{degen} 

In this section we describe the process of degeneration of $X$ to
 a surface $X_0$ which is a union of planes (each homeomorphic to
 $\CP^2$) for which we are able to describe the degenerated branch
 curve $S_0$ (a union of lines) and compute the induced
 (degenerated) braid monodromy. 

We give a defenition of a degeneration.

\begin{definition} \label{df1} {\bf  Projective degeneration}
Let $k  : Y \rightarrow \C \P^n, k'  : X \rightarrow \C \P ^n$ be
 projective embeddings.  We say that $k'$ is   
{\bf  a projective degeneration} of $k$ if there exist an algebraic variety $V$
and $\pi  : V \rightarrow \C,$ \st for every irreducible component 
$V_i$ of $V, \ \pi(V_i) = \C$ and $\pi^{-1}(0)
\simeq X, \ \pi^{-1}(1) \simeq Y, \ \pi^{-1}(1)$ is a generic
fiber. Moreover there exists a regular morphism $F  : \ V \rightarrow
\C \P^n \times \C,$ \st $F(\pi^{-1}(t))\subset  \C\P^n \times t,
\ F_t = F \mid_{\pi^{-1}(t)} \ : \pi^{-1}(t) \rightarrow \C \P^n
\times t$ is a projective embedding of $\pi^{-1}(t). 
\ F_0 = k'$
and $F_1 = k$ under the identification of $\pi^{-1}(0)$ and
$\pi^{-1}(1)$ with $X$ and $Y$, respectively.
\end{definition}

\begin{definition}\label{df2}
{\bf Total degeneration} Let $V$ be an algebraic variety
of dimension $n$.  A projective degeneration  $V \rightarrow W$
will be called  {\bf a total degeneration}  of $V$ if $W$ is a
union of linear spaces of dimension $n$.  
\end{definition}

Let $k_i$ and $l_i$ ($i=1,2,3$) denote the three lines into which two distinct copies of $T$ degenerate, then $T\times T$ is degenerated
 to a union of the nine quadrics $\set{k_i \times l_j}$, each
 homeomorphic to $\CP^1 \times \CP^1$. In Figure 1 below, each square represents one of the quadrics $k_i \times l_j$. Since $k_i$ intersects with $k_j$ for every $i,j$ (and likewise for $l_i,l_j$), we identify the extreme edges (right and left, up and
 down) in the figure. To get a total degeneration, we further divide each quadric to a union of 2 planes, each homeomorphic to $\C \P^2$.

\begin{figure}[h]\label{fig1}
\epsfxsize=3cm 
\epsfysize=3cm 
\begin{minipage}{\textwidth}
\begin{center}
\epsfbox {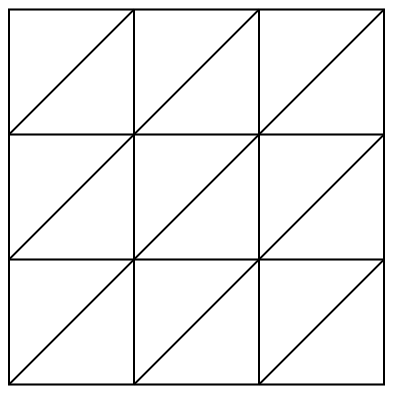}
\end{center}
\end{minipage}
\caption{}\end{figure}

 In Figure 1 the triangles represent copies of $\CP^2$, where $X = T \times T$
 is degenerated into the union of the $18$ planes. The lines represent the intersection of planes. Every plane intersects exactly with three others.
 Note again that we identify extreme edges, so there are $27$ distinct lines.
 The union of the intersection lines is the
 branch curve $S_0$ of the degenerated object of $X$.

 For the regeneration process we have to choose systematically the order of the curves on the fiber over each point. We fix a numeration of the
 vertices, and use it to numerate the intersection lines.

The order of vertices is chosen to be lexicographic, see Figure 2.
\begin{figure}[h]\label{fig2}
\epsfxsize=4cm 
\epsfysize=4cm 
\begin{minipage}{\textwidth}
\begin{center}
\epsfbox {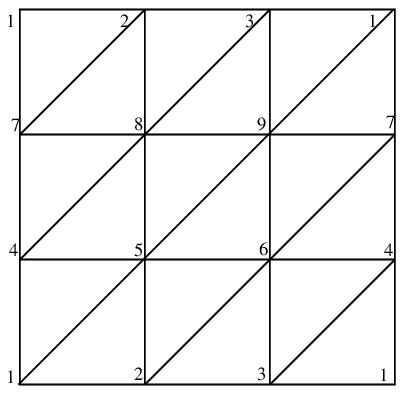}
\end{center}
\end{minipage}
\caption{}\end{figure}

 Now let $L_1$ and $L_2$ be two edges, where $L_i$ has vertices $\alpha_i < \beta_i$. We set
 $L_1 < L_2$ iff $\beta _1 < \beta _2$, or $\beta _1 = \beta _2$ and $\alpha _1 < \alpha _2$ (see Figure 3).

\begin{figure}[h]\label{c05}
\epsfxsize=4cm 
\epsfysize=4cm 
\begin{minipage}{\textwidth}
\begin{center}
\epsfbox {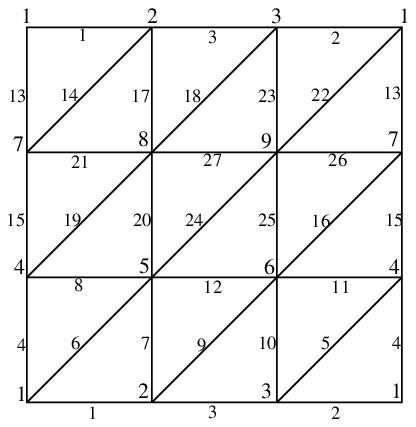}
\end{center}
\end{minipage}
\caption{}\end{figure}

In a similar way we numerate also the quadrics and planes, see Figures 
4 and 5.

\begin{figure}[h]\label{c06}
\epsfxsize=4cm 
\epsfysize=4cm 
\begin{minipage}{\textwidth}
\begin{center}
\epsfbox {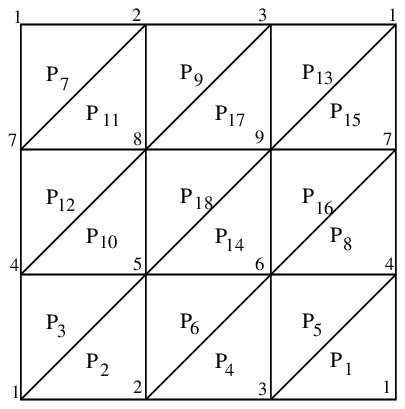}
\end{center}
\end{minipage}
\caption{}\end{figure}

\begin{figure}[h]\label{c07}
\epsfxsize=4cm 
\epsfysize=4cm 
\begin{minipage}{\textwidth}
\begin{center}
\epsfbox {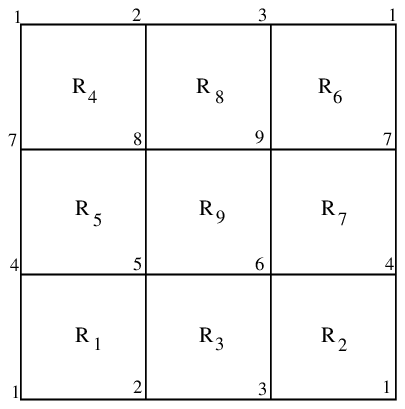}
\end{center}
\end{minipage}
\caption{}\end{figure}

We have nine points in the degenerated object  (see Figure 6).  Each one of them is a 6-point. We numerate each point in a local numeration from 1 to 6.  This local numeration is compatible with the global one.

\begin{figure}[h]\label{6points}
\epsfxsize=26cm 
\epsfysize=11cm 
\begin{minipage}{\textwidth}
\begin{center}
\epsfbox {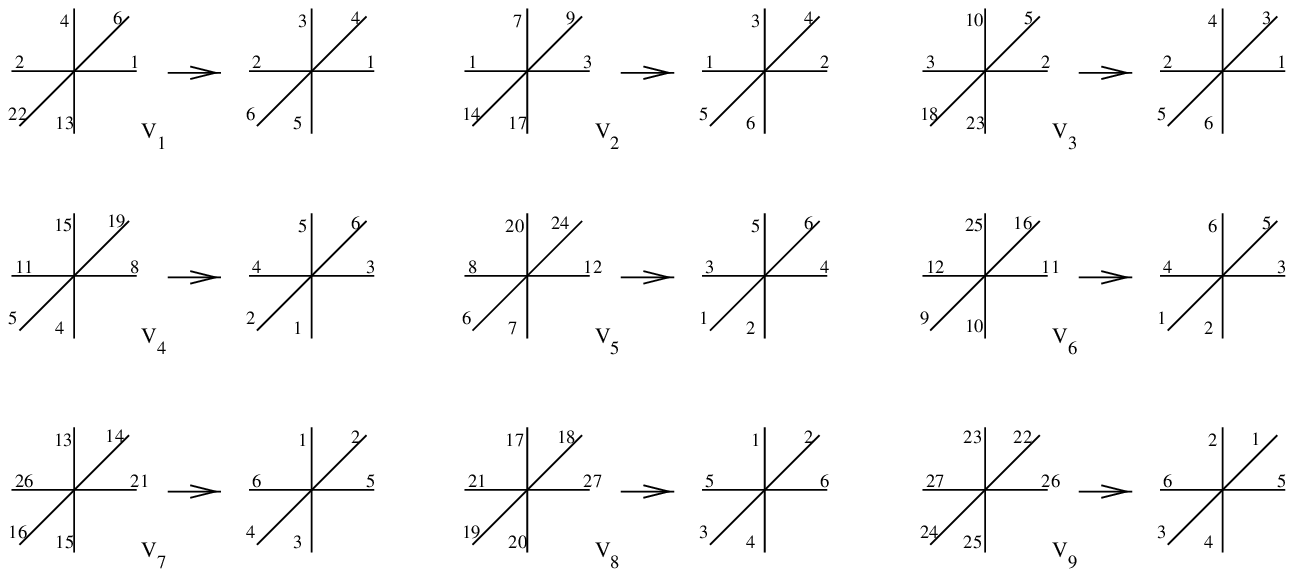}
\end{center}
\end{minipage}
\vspace{-4cm}
\caption{}\end{figure}

We can summerize the resulting degenration as follows.
We have a sequence of projective degenerations:

$T \times T = Z^{(0)} \leadsto Z^{(1)} \leadsto \cdots Z^{(i)}
\leadsto Z^{(i+1)} \leadsto \cdots \leadsto Z^{(8)}$. 

$Z^{(8)}$ is
the degenerated object, see Figure 7. 
Let  $\{\hat{L}_i \}^{27}_{i=1}$ be the 27 lines in $Z^{(8)}$
and  $\{j\}^9_{j=1}$  their intersection points (see Figure 3).

Take generic projections $\pi^{(i)} : Z^{(i)} \rightarrow \C\P^2$
for $0 \leq i \leq 8$.  Let $S^{(i)}$ be the branch curve
of the generic projection for each
$i$ and let $S^{(i+1)}$ be a degeneration of $S^{(i)}$ for $0 \leq  i \leq 7$.

The degenerated branch curve $S^{(8)}$ has a degree of 27.

\begin{figure}[h]\label{fig4}
\epsfxsize=11cm 
\epsfysize=9cm 
\begin{minipage}{\textwidth}
\begin{center}
\epsfbox {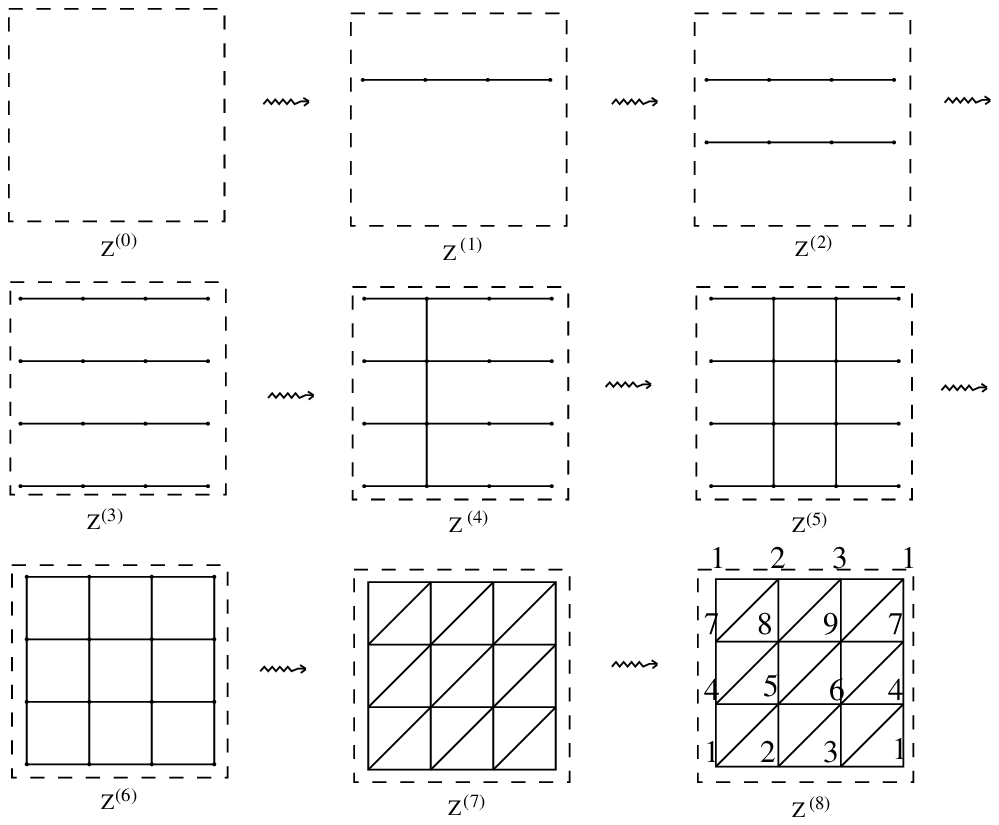}
\end{center}
\end{minipage}
\caption{}
\end{figure}

\section{The braid monodromy}\label{bm}

 Consider the following situation (Figure 8).
 $S$ is an algebraic curve in $\C^2\ , \ p = \mbox{deg } S$.
 $\pi: \C^2 \rightarrow \C$ is a generic projection on the first 
coordinate, 
 $K(x) = \{y \mid (x,y) \in S\}$ is the projection to the $y$-axis. Let
 $N = \{x \mid \# K(x) < p \}$ and 
 $M' = \{ x \in S \mid \pi_{\mid x} \mbox{ is not \'{e}tale at } x \}$ \st  $\pi (M') = N$.
 Assume $\# (\pi^{-1}(x) \cap M') =1, \forall x \in N$. Let $E$
 (resp. $D$) be a closed disk  on $x$-axis (resp. $y$-axis), \st
 $M' \subset E \times D, N \subset \mbox{ Int } (E)$. We choose $u \in \partial E, x \ \ll u \;\;\; \forall x \in N$, and define $\C_u = \pi^{-1}(u)$, \st $K = \{q_1 , \cdots , q_p\}$.
Let $\set{A_j}^q_{i=1}$ be the set of singular points and $\set{x_j}^q_
{i=1}$ their projection on the $x$-axis. 

We construct a g-base for the fundamental group $\pi_1(E - N, u)$.
Take a set of paths  $\set{\g_j}^q_{i=1}$ which connect $u$ with the points $\set{x_j}^q_{i=1}$. Now encircle each $x_j$ with a small orriented counterclockweise circle $c_j$. Call the path segment from $u$ to the boundary 
of this circle as $\g'_j$. We define an element (a loop) in 
the g-base as  
$\ell(\g_j)= \g'_j c_j {\g'}^{-1}_j$. Shortly we denote it as $\delta_j$.

 Let $B_p[D,K]$ be the braid group, and let $H_1 , \cdots , H_{p-1}$ be its frame (for complete definitions, see
 \cite[III.1/2]{MoTe4}).  

\begin{figure}[h]\label{fig5}
\epsfxsize=7cm 
\epsfysize=6cm 
\begin{minipage}{\textwidth}
\begin{center}
\epsfbox {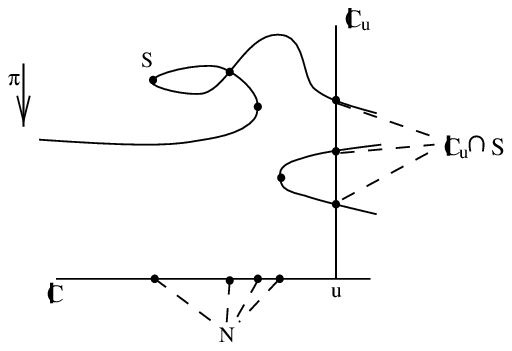}
\end{center}
\end{minipage}
\caption{}\end{figure}

We are now able to introduce {\em  the braid monodromy}.

\label{df9}
\begin{definition} {\bf Braid monodromy of an affine  curve}
 $\mbox{\boldmath $S$}$ {\bf \wrt } $\mbox{\boldmath $E \times D, \pi, u$}$

The braid monodromy of $S$ is a map $\varphi: \pi_1(E - N, u) \rightarrow 
B_p[D,K]$ defined as follows: 
every loop in $E - N$ starting at $u$ has liftings to a system of $p$ paths in $(E - N) \times D$ starting at $q_1, \cdots ,q_p$.
 Projecting them to $D$ we get $p$ paths in $D$ defining a motion $\{q_1(t), \cdots , q_p(t)\}$ of $p$ points in $D$ starting and ending at $K, 0 \leq t \leq 1$. This motion defines a braid in $B_p [D , K]$.
 \end{definition}

 \label{th13}
 \begin{theorem}{\bf The Artin Theorem [A2]}

 Let $S$ be a curve and let $\delta _1, \cdots, \delta _q$ be a g-base of 
$\pi _1 (E -N, u)$.
 Assume that the singularities of $S$ are cusps, nodes and  tangent points of a parabola/hyperbola with a line and branch points. Let 
$\varphi: \pi_1(E - N,u) \rightarrow B_p$ be the braid monodromy.
 Then for all $i$, there exist a halftwist $Z_j \in B_p$ and $r_j \in
 \Z$, \st
 $\varphi (\dl _j)  = Z_j ^{r_j}$
  and $r_j$ depends on the type of the singularity:
 $
 r_j =  1,2,3,4$ for  branch point, node,  cusp,  tangent point respectively.
 \end{theorem}

In order to know what $Z_j$ is, we recall the braid monodromy algorithm
(the algorithm appears in detail in [Am2] and [Am5]).

We take a singular point $A_j$ and project it to the $x$-axis to get $x_j$. We choose a point $x'_j$ close to  $x_j$, \st $\pi^{-1}(x'_j)$ is a typical fiber. We choose it in this way:\\
If $A_j$ is  a node, a cusp, a tangent point or a branch point(locally defined by $y^2-x=0$), then   $x'_j$ is on the right side of $x_j$.
If  $A_j$ is  a branch point defined locally by  $y^2+x=0$, then  $x'_j$ is on the left side of $x_j$.\\
We encircle $A_j$ with a small circle, and we consider the 
maximal and minimal irreducible components in the circle which 
meet at $A_j$ and intersect the fiber $\pi^{-1}(x'_j)$ . 
We connect them by a skeleton $\xi_{x'_j}$ and start to move 
from one typical fiber to another by a Lefschetz diffeomorphism $\Psi$ 
(see [Am2, Subsection 1.9.5]), applying  diffeomorphisms 
which correspond to the singularities
we pass near by([Am2, Definition 1.29])on $\xi_{x'_j}$  
till we arrive the typical fiber $\C_u$. 
We get a skeleton L.V.C.$(\g_j)=(\xi_{x'_j}) \Psi$ corresponding to the halftwist $Z_j$. We denote  $Z_j=\Delta < (\xi_{x'_j}) \Psi>$. 

\begin{remark}
We compute the braid monodromy using models, since computations are easier
in this way([Am2, Subsection 1.9.4]). There exists a continuous family of
diffeomorphisms $\set{\beta_x}$ which are transition functions from $(D,K)$
to its models([Am2, Lemma 1.37]).
\end{remark}

\begin{corollary} \label{co38}
$\varphi(\delta_j) = \Delta < (\xi_{x'_j}) \Psi>^{r_j}$.
\end{corollary}

We compute a presentation for $\pitil$ by applying the \vK \ Theorem on  the braid monodromy factorization.  Now we define {\em a braid monodromy factorization}.

\begin{proposition}{\bf [MoTe4, Proposition VI.2.1]}\label{pr14}
 Let $S$  be an algebraic curve of degree $p$ in $\C \P^2$.  Let
 $\pi, u, D, E, \C_u$ be as above.  Let $\varphi$ be the braid
 monodromy of $S$ \wrt   $\pi, u.$  Let $\delta_1, \cdots ,
 \delta_q$ be a g-base of $\pi_1(E - N,u)$.  Then
 $
 \prodl^q_{j=1} \varphi(\delta_j) = \Delta^2_{p} = \Delta^2_{p} 
[u \times D,K \cap S].$\\
 Such a presentation of $\Delta^2_{p}$ as a product of (positive)
 elements is called a {\em braid monodromy
 factorization of $\Delta^2_{p}$ associated to $S$}.
\end{proposition}

 \begin{theorem}{\bf [MoTe4, Prop. V.2.2]}\label{th19}
 Let $\{H_i\}^{p-1}_{i=1}$ be a frame of $B_p$.  Then for $p \geq 2$, the center
 of $B_p$ is generated by $\Delta^2_p = (H_1
 \cdot \ \cdots \ \cdot H_{p-1})^p$.  Moreover, deg$\Delta^2_p = p \cdot (p-1)$.
 \end{theorem}

\section{$S^{(8)}$ and $\varphi^{(8)}$}\label{ssec:232}
\indent
Recall that we project the degenerated object $Z^8$ to $\C\P^2$. We get a branch curve which is a line arrangement, 
$S^{(8)} = \bigcup\limits^{27}_{i=1} L_i$. $\{j\}^9_{j=1}$ are the vertices in $Z^{(8)}$.  Let $V_j = \pi^{(8)}(j)$ for $j = 1, \cdots , 9,$ so $V_j$ are the singular points of $S^{(8)}$.

Let $C$ be the union of all lines connecting pairs of the $V_j$'s.
 $S^{(8)}$ is a subcurve of $C$.   $C$ can be chosen as discussed in  [MoTe4, Theorem  IX.2.1].
This theorem gives a description of a braid monodromy for $C \ :
\Delta^2_C = \prodl^9_{r=1} Cr \Delta^2_r$ with an appropriate
description of L.V.C.. We use this formula to obtain a description
of $\vp^{(8)}$ by deleting
all factors that involve lines which do not appear in $S^{(8)}$.
Thus, we get $\Delta^2_{S^{(8)}} = \prodl^9_{j=1} \tilde{C}_j \tilde{\Delta}^2_j$.
 We discuss each
factor separately.

\underline{ $\tilde{\Delta}^2_j$:}  Each one of the nine points in $\Se$ is
a 6-point. The numeration of every six lines intersecting at a point 
is shown in Figure 6. We
numerate each six lines by a local numeration from 1 to 6. The
local braid monodromies $\vp^{(8)}_j$ are
 $\tilde{\Delta}^2_j = \Delta^2 < 1, \cdots , 6>$ for $1 \le j \le 9$
 (recall [Am2, Definition 1.27]  of $\Delta < 1 , \cdots ,
 6>$).
 Each $\tilde{\Delta}^2_j$ is regenerated during the regeneration of  $Z^{(8)}$ to $Z^{(6)}$.

\underline{$\tC_j$:} We have  27 lines in the degenerated object.  
Each line $L_i, 1 \le i \le 27$, can be
represented as a pair, by its two end vertices.

Let $u$ be a point, $u \notin S$, \st $\#(\pi^{(8)})^{-1} (u) =
27$.

Let $q_i = L_i \cap \C_u$ be a real point.  We take two lines $(i,k) =
 L_p, (j,\ell) = L_t$ and apply [MoTe4, Theorem IX.2.1] by defining
 $D_t = \prodl_{\stackrel{p<t}{L_p \cap L_t = \emptyset}} \tilde{Z}^2_{pt}$
  ($\tZZ^2_{pt}$ are formulated in [MoTe4, p. 526]). $\tZZ^2_{pt}$ are 
related to  parasitic intersections, since there are lines which do not intersect in $\C\P^8$ but may intersect in $\C\P^2$.

\def\Zover4{\mathop{\lower 10pt \hbox{$ {\buildrel{\displaystyle\bar{Z}^2_{15}}\over{\hspace{-.2cm}\scriptstyle{(4)}}} $}}  \nolimits}
\def\Zoversix{\mathop{\lower 10pt \hbox{$ {\buildrel
{\displaystyle\bar{Z}^2_{i \; 7}}\over{\hspace{-.2cm}\scriptstyle{(6)}}} $}}  \nolimits}
\def\Zoverseven{\mathop{\lower 10pt \hbox{$ {\buildrel{\displaystyle\bar{Z}^2_{i \; 8}}\over{\hspace{-.2cm}\scriptstyle{(6)(7)}}} $}}  \nolimits}
\def\Zovernine{\mathop{\lower 10pt \hbox{$ {\buildrel{\displaystyle
\bar{Z}^2_{i \; 10}}\over{\hspace{-.2cm}\scriptstyle{(9)}}} $}}  \nolimits}
\def\Zovernineten{\mathop{\lower 10pt \hbox{$ {\buildrel{\displaystyle\bar{Z}^2_{i \; 11}}\over{\hspace{-.2cm}\scriptstyle{(9)(10)}}} $}}  \nolimits}
\def\Zovernineteneleven{\mathop{\lower 10pt \hbox{$ {\buildrel{\displaystyle\bar{Z}^2_{i \; 12}}\over{\hspace{-.2cm}\scriptstyle{(9)(10)(11)}}} $}}  \nolimits}
\def\Zoverthirteen{\mathop{\lower 10pt \hbox{$ {\buildrel{\displaystyle\bar{Z}^2_{i \; 14}}\over{\hspace{-.2cm}\scriptstyle{(13)}}} $}}  \nolimits}

\def\Zoverthirteenf{\mathop{\lower 10pt \hbox{$ {\buildrel{\displaystyle\bar{Z}^2_{i \; 15}}\over{\hspace{-.2cm}\scriptstyle{(13)(14)}}} $}}  \nolimits}

\def\Zoverthirteenff{\mathop{\lower 10pt \hbox{$ {\buildrel{\displaystyle\bar{Z}^2_{i \; 16}}\over{\hspace{-.2cm}\scriptstyle{(13)(14)(15)}}} $}}  \nolimits}

\def\Zoverseventeen{\mathop{\lower 10pt \hbox{$ {\buildrel{\displaystyle\bar{Z}^2_{i \; 18}}\over{\hspace{-.2cm}\scriptstyle{(17)}}} $}}  \nolimits}

\def\Zoverseventeene{\mathop{\lower 10pt \hbox{$ {\buildrel{\displaystyle\bar{Z}^2_{i \; 19}}\over{\hspace{-.2cm}\scriptstyle{(17)(18)}}} $}}  \nolimits}

\def\Zoverseventeenen{\mathop{\lower 10pt \hbox{$ {\buildrel{\displaystyle\bar{Z}^2_{i \; 20}}\over{\hspace{-.2cm}\scriptstyle{(17)(18)(19)}}} $}}  \nolimits}

\def\Zoverseventeentwenty{\mathop{\lower 10pt \hbox{$ {\buildrel{\displaystyle\bar{Z}^2_{i \; 21}}\over{\hspace{-.2cm}\scriptstyle{(17)-(20)}}} $}}  \nolimits}

\def\Zoverseventeentwt{\mathop{\lower 10pt \hbox{$ {\buildrel{\displaystyle\bar{Z}^2_{i \; 23}}\over{\hspace{-.2cm}\scriptstyle{(22)}}} $}}  \nolimits}

\def\Zovertwentytwothree{\mathop{\lower 10pt \hbox{$ {\buildrel{\displaystyle\bar{Z}^2_{i \; 24}}\over{\hspace{-.2cm}\scriptstyle{(22)(23)}}} $}}  \nolimits}

\def\Zovertwothreefour{\mathop{\lower 10pt \hbox{$ {\buildrel{\displaystyle\bar{Z}^2_{i \; 25}}\over{\hspace{-.2cm}\scriptstyle{(22)(23)(24)}}} $}}  \nolimits}

\def\Zovertwofive{\mathop{\lower 10pt \hbox{$ {\buildrel{\displaystyle\bar{Z}^2_{i \; 26}}\over{\hspace{-.2cm}\scriptstyle{(22)-(25)}}} $}}  \nolimits}

\def\Zovertwosix{\mathop{\lower 10pt \hbox{$ {\buildrel{\displaystyle\bar{Z}^2_{i \; 27}}\over{\hspace{-.2cm}\scriptstyle{(22)-(26)}}} $}}  \nolimits}

Thus: $D_1 = D_2 = D_3 = Id \ , \ D_4 = Z^2_{34} \ , \ D_5 = \Zover4 \ , \ D_6 = \bZ^2_{36} Z^2_{56}$   \ ,
$D_7  =   \prodl_{i=2,4,5} \Zoversix \ ,\\
  D_8 = \prodl^3_{i=1} \;\;
\Zoverseven \ , \ D_9=\prodl_{i=2,4,5,6,8} \bZ^2_{i \; 9} \ , \ D_{10} =
\prodl_{i = 1,4,6-8} \;\; \Zovernine \ , \
D_{11}  =  \prodl_{i=1-3,6,7} \;\; \Zovernineten \ , \\
D_{12} = \prodl^5_{i=1}
\;\; \Zovernineteneleven \ , \
D_{13} = \prodl^{12}_{\stackrel{\scriptstyle i=3}{\scriptstyle i \neq 4,6}}\;\;  \bar{Z}^2_{i \; 13} \ , \ D_{14} = \prodl^{12}_{\stackrel{\scriptstyle i=2}{\scriptstyle
i \neq 3,7,9}} \;\; \Zoverthirteen \ , \ $
$D_{15}  =   \prodl^{12}_{\stackrel{\scriptstyle i=1}{\scriptstyle i \neq 4,5,8,11}} \;\;  \Zoverthirteenf \ , \\ D_{16} = \prodl^8_{i=1} \;\; \Zoverthirteenff \ $ , \
$D_{17} = \prodl^{16}_{\stackrel{\scriptstyle i=2}{\scriptstyle i \neq 3,7,9,14}} \bZ^2_{i \; 17} \ , \
$
$
D_{18} = \prodl^{16}_{\stackrel{\scriptstyle i=1}{\scriptstyle i \neq 2,3,5,10}} \;\; \Zoverseventeen \ , \ $
$D_{19}=\prodl^{16}_{\stackrel{\scriptstyle i=1}{\scriptstyle i \neq 4,5,8,11,15}} \;\;  \Zoverseventeene  \ $ , \\
$D_{20} = \prodl^{16}_{\stackrel{\scriptstyle i=1}{\scriptstyle i \neq 6-8,12}} \;\; \Zoverseventeenen \ , \ D_{21} = \prodl^{12}_{i=1} \;\; \Zoverseventeentwenty \ , \ D_{22}=\prodl^{21}_{\stackrel{\scriptstyle i=3}{\scriptstyle i \neq 4,6,13}} \;\; \bZ^2_{i \; 22}$ \ , \
$D_{23} = \prodl^{21}_{\stackrel{\scriptstyle i=1}{\scriptstyle
i \neq 2,3,5,10,18}} \;\; \Zoverseventeentwt  \ ,$ \\$
D_{24} = \prodl^{21}_{\stackrel{\scriptstyle i=1}{\scriptstyle
i \neq 6-8,12,20}} \;\; \Zovertwentytwothree \ ,
\ D_{25}=\prodl^{21}_{\stackrel{\scriptstyle i=1}{\scriptstyle i \neq 9-12,16}} \;\;
\Zovertwothreefour $\ , \
$D_{26} = \hspace{-.1cm}\prodl^{20}_{\stackrel{\scriptstyle i=1}{\scriptstyle
i \neq 13-16}} \ \  \Zovertwofive$ \ , \
$D_{27} = \prodl^{16}_{i=1} \;\; \Zovertwosix$.

Define $\tilde{C}_j = \prodl_{V_j \in L_t}\hspace*{-.2cm}D_t$, where  $V_j$ is the \underline{small} vertex among two vertices in $L_t$.\\
$\tC_1 = D_1 \cdot D_2  \cdot D_4 \cdot D_6  \cdot D_{13} \cdot D_{22} \ ; \
 \tC_2 = D_3 \cdot D_7  \cdot D_9 \cdot D_{14} \cdot D_{17}\ ; \ \tC_3 = D_5 \cdot D_{10} \cdot  D_{18} \cdot D_{23}\ ; \
\tC_4 = D_8 \cdot D_{11} \cdot  D_{15} \cdot D_{19}\ ; \
\tC_5 =  D_{12} \cdot  D_{20} \cdot D_{24}\ ; \
\tC_6 =  D_{16} \cdot  D_{25}\ ; \ \tC_7 =  D_{21} \cdot  D_{26} \ ; \
\tC_8 =  D_{27}\ ; \
\tC_9 = Id$.

We start with a suitable description
 of the degenerated braid monodromy factorization $\vp^{(8)}$.
 We compare degrees of the factorized  expressions given in $\vp^{(8)} = \prodl^9_{j=1} \tC_j \tilde{\Delta}^2_j$.\\ Deg $S^{(8)} = 27$, so: $\vp^{(8)} = \Delta^2_{27}$ and deg$\Delta^2_{27} =
27 \cdot 26 = 702$. $\sum\limits^9_{j=1} \mbox{deg} \tC_j = 432$
and $\sum\limits^9_{j=1} \mbox{deg} \tilde{\Delta}^2_j = 9 \cdot (6 \cdot 5)
= 270$. Thus $ \sum\limits^9_{j=1} (\mbox{deg} \tC_j \tilde{\Delta}^2_j) =
702$. Therefore deg$\vp^{(8)} = \mbox{deg} \prodl^9_{j=1} \tC_j
\tilde{\Delta}^2_j$.

\section{The regeneration of $S^{(8)}$}

Recall [Am2, Definition 1.47] for 
a regeneration of a braid group, a frame and an embedding.

The degenerated branch curve $S_0$ has a degree of 27.  We regenerate $S_0$ 
in a way that each point on the typical fiber is replaced by two close points, 
since each line in $S_0$ becomes a conic or two parallel lines.  

During the regeneration, each branch point is replaced by two
branch points, a node is replaced by two or four nodes and a
tangent point is replaced by three cusps(see [MoTe5]).
We obtain a regenerated branch curve $S$, which has a degree of 54.

$B_{27}$ is the braid group
corresponding to the degenerated curve $S_0$, therefore $B_{54}$
is the braid group corresponding to the regenerated curve $S$.

In Section \ref{degen} we constructed a series of algebraic
varieties $Z^{(i)}$ for $i = 0, \cdots , 8$.  $Z^{(0)} = T \times
T$ and each $Z^{(i+1)}$ is a degeneration of $Z^{(i)}$ for $i = 0,
\cdots , 7$.

For every $Z^{(i)}$ we take a generic projection $\pi^{(i)}$ to
$\C\P^2$.  We denote by $S^{(i)}$ the corresponding branch curve.
One can choose a family of generic projections to $\C\P^2$, \st
$S^{(i+1)}$ is a projective degeneration of $S^{(i)}$.

Let us choose a line at infinity transversal to $S^{(i)},$ $0 \le i \le 8,$ and a system of coordinates $x;y$ in the affine part $\C^2,$ which is in a 
general position to each $S^{(i)},$ $0 \le i 
\le 8$. Let $m^{(i)} = \deg S^{(i)} = \deg \pi^{(i)}$. Let 
$N^{(i)} = \{x\in x\mbox{-axis} \bigm| \#{\pi^{(i)}}^{-1}(x) \le m^{(i)} \}.$ $u$ is a real point in the $x$-axis, $N^{(i)} \ll u \ \ \forall i$. Let $\C_u = \{(u,y) \mid y \in \C \}$ and $K^{(i)} = S^{(i)} \cap \C_u.$ Denote $B^{(i)} = B_{m_i} \left[\C_u, K^{(i)}\right].$ $\vp^{(i)} : \pi_1 (x\mbox{-axis} - N^{(i)}, u) \rightarrow  B_{m_i}[\C_u,K^{(i)}]$ is the braid monodromy w.r.t. $S^{(i)}, \pi^{(i)}, u.$

We want to compute the regenerated braid monodromy factorization $\vp^{(0)}.$

By the construction of $Z^{(8)},$\ $K^{(8)}$ is a set of real points.
We start with $\vp^{(8)}$ and use the below three regeneration rules  
and the regeneration of a frame to obtain a formula for $\vp^{(0)}$,
proceeding step by
step from $i = 8$ to $i = 0.$

For every point $c$ in $K^{(8)}$ there exists an $i$, \st $c$ is replaced in $K^{(k)}$
for $k\le i$ by two points close to each other.
In $K^{(0)}$ we obtain twice as many points as in $K^{(8)}.$ This implies that $B^{(i)}$ is a natural regeneration of $B^{(i + 1)}.$ We choose a real frame in $B^{(8)}$ and inductively  a frame of $B^{(i)}$ which is a regeneration of the frame of $B^{(i + 1)}.$
We describe every  $\vp^{(i)}$ using these frames.

Let us recall from [Am2] the way to compute
$\varphi^{(0)}$.  For every 6-point $V_j$ in $S^{(8)}, 1 \leq j
\leq 9$, we shall take a small neighbourhood of $V_j$ and analyze
the local braid monodromy $\varphi_j^{(i)}$ of $S^{(i)}_j$ which
results from the singularities of this neighbourhood from $i = 8$
to $i = 0$.  $S^{(i)} = \bigcup\limits^9_{j=1} S^{(i)}_j$, thus
$\varphi^{(i)}$ compounds the braid monodromies
$\varphi_j^{(i)}$.  Each one of these $\varphi_j^{(i)}$ is determined by the regeneration of the
embedding $B_6 [\C_u,K^{(8)}] \hookrightarrow B_{27}
[\C_u,K^{(8)}]$ to $B_{12} [\C_u,K^{(0)}] \hookrightarrow B_{54}
[\C_u,K^{(0)}]$. $\varphi^{(8)}_j$ is induced by braids in $B_6
[\C_u,K^{(8)}]$, which are embedded in $B_{27} [\C_u,K^{(8)}]$. We
take a local real frame of $B_6 [\C_u,K^{(8)}]$ and embed it into
$B_{27} [\C_u,K^{(8)}]$ by a sequence of consecutive simple paths
connecting the points in $K^{(8)}$ below the real line.  When
regenerating each 6-point, we get $\varphi^{(0)}_j$.
$\varphi^{(0)}_j$ is induced by the braids in $B_{12}
[\C_u,K^{(0)}]$, which are embedded in $B_{54} [\C_u,K^{(0)}]$ in
the same way.

We discuss the regeneration in a neighbourhood of every point. We
numerate every  6-point locally from 1 to 6(see Figure 6).  The
point $V_1$ is an intersection point of the lines $L_1, L_2, L_4,
L_6, L_{13}, L_{22}$.  The local numeration is compatible with the
global one:\\
$L_1 \rightarrow 1,  L_2 \rightarrow 2,  L_4 \rightarrow 3, L_6 \rightarrow 4, L_{13} \rightarrow 5, L_{22} \rightarrow 6$.  We also numerate lines which intersect at $\{V_j\}^9_{j=2}$ locally from 1 to 6.

In order to regenerate the degenerated object of $T \times T$ and
in order to compute $\vp^{(0)}$, we follow the sequence of degenerations in 
Figure 7.  
Proceeding in the opposite direction of the arrows, we obtain the
sequence of projective regenerations in the neighbourhoods of the 6-points. 
In the second regeneration ($Z^{(6)}$ in Figure 7) 
each diagonal becomes a conic, which is tangent to
the two lines in which it intersected. It holds by the following lemma.

\begin{lemma}{\bf [MoTe8, Lemma 1]}.\label{lm:24}   Let $V$ be a projective algebraic surface, and $D'$ be a curve in $V$.  Let $f: V \rightarrow \C \P^2$ be a generic projection.  Let $S \subseteq \C \P^2, S' \subset V$ be the branch curve of $f$ and the corresponding ramification curve.  Assume $S'$ intersects $D'$ in $\alpha'$.  Let $D = f(D')$ and $\alpha = f(\alpha ')$.  Assume that there exist neighbourhoods of
$\alpha $ and $\alpha '$, \st their $f_{\mid_{S'}}$ and
$f_{\mid_{D'}}$ are isomorphisms.  Then $D$ is tangent to $S$ at
$\alpha $.
\end{lemma}

We obtain the curves $S^{(6)}_j, \ 1 \le j \le 9$. Each
$S^{(6)}_j$ compounds two conics and four lines, thus $\deg
S^{(6)}_j$ has a degree of 8.

Now we concentrate on the neighbourhoods of the remaining
4-points. 
From the third regeneration on, the vertical lines are doubled and the
horizontal lines become hyperbolas. 
This is w.r.t. their local
neighbourhoods and to the local numeration of lines(see Figures 9, 10, 11).

\begin{figure}[h]\label{fig6}
\epsfxsize=8cm 
\epsfysize=7cm 
\begin{minipage}{\textwidth}
\begin{center}
\epsfbox {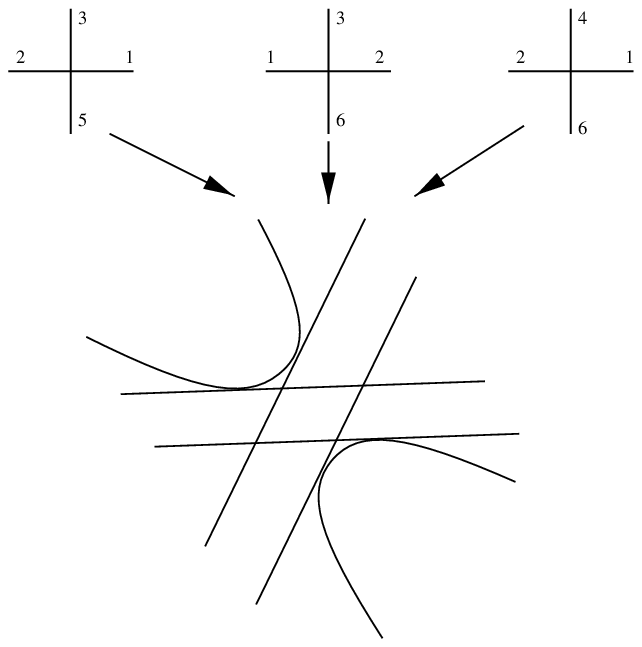}
\end{center}
\end{minipage}
\caption{}\end{figure}

\begin{figure}[h]\label{fig7}
\epsfxsize=8cm 
\epsfysize=7cm 
\begin{minipage}{\textwidth}
\begin{center}
\epsfbox {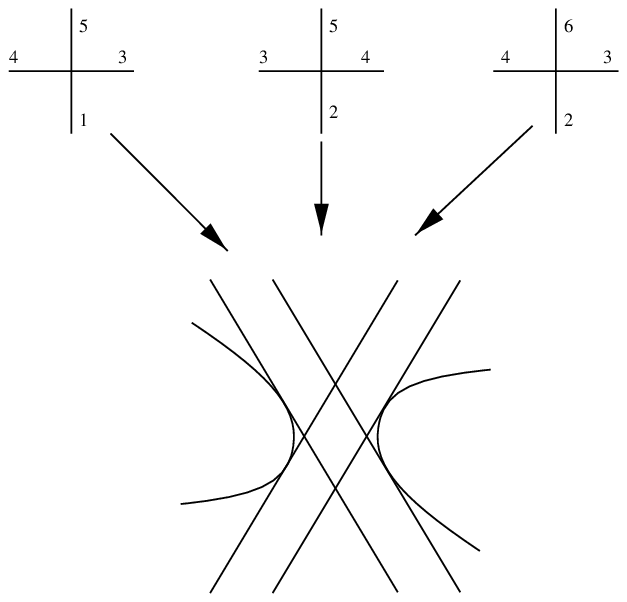}
\end{center}
\end{minipage}
\caption{}\end{figure}

\begin{figure}[h]
\epsfxsize=8cm 
\epsfysize=7cm 
\begin{minipage}{\textwidth}
\begin{center}
\epsfbox {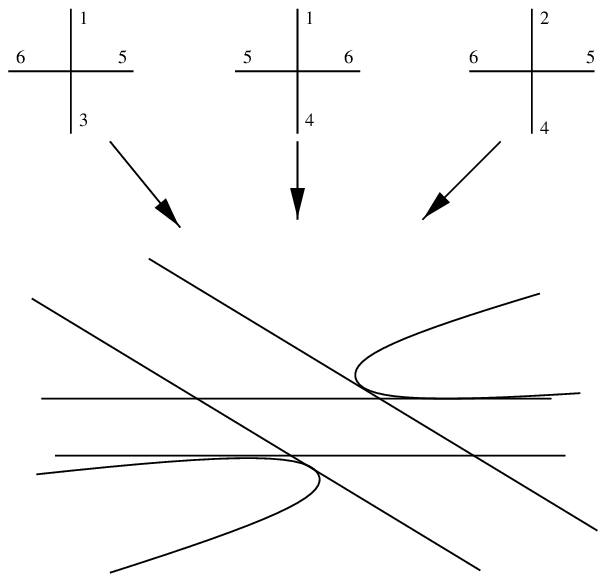}
\end{center}
\end{minipage}
\caption{}\end{figure}

\subsection{The regeneration rules}\label{ssec:1104}

The regeneration regenerates also the braids.

Recall that $u \in \C$, $u$ real, \st $x \ll u \ \ \forall x \in N$. 
Recall also that $K$ is the set of the intersection points of the curve $S$ 
with $\C_u,
\# K = 54$.  Choose $M \in \C_u, M \notin S$ below the real line far enough, 
\st $B_{54}$ does not move $M$.

\begin{lemma}{\bf [Am2, Lemma 2.5]}\label{lm:25}  \ Let 
$\{a_i, a_{i'}, \cdots , a_j, a_{j'}\}$ 
be a finite regenerated set in $K$.\\
Then: (a) $Z^2_{ii',j} = Z^2_{i'j} Z^2_{ij}$;(b)
$Z^{2}_{i',jj'} = Z^{2}_{i'j'} Z^{2}_{i'j}$ ;  (c)
$Z^{-2}_{i',jj'} = Z^{-2}_{i'j} Z^{-2} _{i'j'}$  ; (d)
$\bar{Z}^{-2}_{i',jj'} = \bar{Z}^{-2}_{i'j'} \bar{Z}^{-2} _{i'j}$
 ;  (e) $Z^{-2}_{ii',j} = Z^{-2}_{ij} Z^{-2} _{i'j}$  ;  
(f)
$Z^2_{ii',jj'} = Z^2_{i',jj'} Z^2_{i,jj'}$  ; (g)
$Z^{-2}_{ii',jj'} = Z^{-2}_{i,jj'} Z^{-2}_{i',jj'}$.
\end{lemma}

In what follows we shall describe what happens to a factor in a
factorized expression of a braid monodromy by a regeneration of
the braid group and of the frame.  Following Definition 1.49 of a halftwist 
in [Am2], we have the following regeneration rules.

\begin{theorem}\label{th:151}
{\bf First regeneration rule [Mote5, Lemma 3.1]}\\
A factor of the braid monodromy of the form $Z_{ij}$ is replaced 
in the regeneration  by $Z_{ij'} \cdot Z_{i'j}$. It is invariant
 under $(Z_{ii'} Z_{jj'})^r \ \ \forall r \in \Z$.
\end{theorem}

\begin{theorem}\label{th:152} 
{\bf Second regeneration rule [Mote5, Lemma 3.2]}\\
A factor of the form $Z^2_{ij}$ is  replaced by a factorized  expression  
$Z^2_{ii',j}=Z^2_{i'j} \cdot Z^2_{ij}$ , by $Z^2_{i,jj'}=Z^2_{ij'} \cdot 
Z^2_{ij}$ or by  $Z^2_{ii',jj'}$.\\
$Z^2_{ii',j}$ is invariant under
$Z^p_{ii'}$, and $Z^2_{ii',jj'}$
 is invariant under $Z^p_{ii'} Z^q_{jj'}$.
\end{theorem}

\begin{theorem}\label{th:153}
{\bf Third regeneration rule [Mote5, Lemma 3.3]}\\
A factor of the form $Z^4_{ij}$ in the braid monodromy
factorized expression is replaced by $Z^{3}_{i,jj'}=
(Z^3_{ij})_{\rho_j} \cdot (Z^3_{ij}) \cdot
(Z^3_{ij})_{\rho^{-1}_j}$ (where $\rho_j=Z_{jj'}$).
$Z^{3}_{i,jj'}$ is invariant under $Z^r_{jj'}$.
\end{theorem}

\section{Technical details related to the braid monodromy algorithm}\label{ssec:241}

\def\pitil{\tilde{\pi}_1}
In the end of the regeneration we obtain the regenerated branch
curve $S^{(0)}=S, \ S$ has a degree of 54.  We denote the
regenerated braid monodromy factorization $\varphi^{(0)}$ as
$\Delta^2_{54}$.

$\Delta^2_{54} = \prodl^9_{i=1} C_i H_{V_i}$, where
$C_i$ are the regenerated $\tilde{C}_i$ and $H_{V_i}$ are the
local braid monodromies $\varphi^{(0)}_i , \ 1 \leq i \leq 9$.

\subsection{Conjugation property}
We apply the following Conjugation Property on some braids in order 
to simplify some of the paths, corresponding to  braids.
\def\vK{van Kampen}

\begin{proposition}{\bf [Mo2]} \label{pp:21}
If $Z_1$ and $Z_2$ are two consecutive braids in a list of braids,
then $Z^{-1}_1 Z_2 Z_1$ induces a relation on $\pi_1(\C^2-S,M)$.
\end{proposition}

\subsection{Complex Conjugation}\label{ssec:243}

All explicit curves which we use are defined over $\R$.  
Therefore, in all our arguments we can change roles of
upper and lower halfplanes in the  $x$-axis and in the typical fiber by simultaneous complex conjugation of the $x$-axis and the $y$-axis. Such a
change requires also to replace by opposite the order of factors
in all our expressions for the braid monodromy. By [MoTe8, Lemma 19], the resulting 
object is also a braid monodromy factorization. 

\subsection{Some invariance properties}\label{invariance}

\begin{corollary}\label{co:158}
From the regeneration rules we obtain the following:\\
Invariance Rule I:\ \ $Z^{1}_{ij}$ is invariant under $(Z_{ii'} Z_{jj'})^q$; \\
Invariance Rule II:\ \ $Z^{2}_{i,jj'}$ is invariant under
$Z_{jj'}^q$ \ ; \
$Z^{2}_{ii',jj'}$ is invariant under $Z_{ii'}^p Z^q_{jj'}$ ;\\
Invariance Rule III:\ \ $Z^{3}_{i,jj'}$ is invariant under
$Z_{jj'}^q$ .
\end{corollary}

\begin{lemma}\label{lm:159}{\bf Chakiri.} Let $g = g_1 \cdots g_k$ be a
factorized expression in a group $H$.  Then $g_1 \cdots g_k$ is
invariant under $g^m$ for any $m \in \Z$.
\end{lemma}

\begin{remark}\label{rm:160}{\bf [MoTe8, Invariance Remarks]}
\begin{itemize}

\item[{(i)}] To prove invariance of $g_1 \cdots g_k$ under $h$ it is enough to prove that $g_1 \cdots g_t$ is invariant under $h$ and $g_{t+1} \cdots g_k$ is invariant under $h$.  Thus we can divide a factorization into subfactorizations and prove invariance on each part separately.

\item[{(ii)}] A factorized expression of one element $g_1$ is invariant under $h$ iff $g_1$ commutes with $h$.

\item[{(iii)}] If the product of elements that commutes with $h$ is invariant under $h$, the relevant factorizations are actually equal.

\item[{(iv)}] If $\sigma_1 \cap \sigma_2 = \emptyset$ then $H(\sigma_1)$ commutes with $H(\sigma_2)$.

\item[{(v)}]If $g$ is invariant under $h_1$ and $h_2$ then $g$ is invariant under $h_1 h_2$.
\end{itemize}
\end{remark}

Consider the above statements. Let  $z_{ij}$ be a path connecting $q_i$ or $q_{i'}$ with $q_j$ or
$q_{j'}$ and $Z_{ij} = H(z_{ij})$ its corresponding halftwist.  We can
conjugate $Z_{ij}$ by $\rho_j, \rho^{-1}_j, \rho_i, \rho^{-1}_i$.
These conjugations are the actions of $\rho_j$ and $\rho^{-1}_j$
on the ``tail'' of $z_{ij}$ within a small circle around $q_j$ and
$q_{j'}$ and the action of $\rho_i$ and $\rho^{-1}_i$ on the
``head'' of $z_{ij}$ in a small circle around $q_i$ and $q_{i'}$.
The ``body'' of $z_{ij}$ does not change under $\rho^{\pm 1}_j$
and $\rho^{\pm 1}_i$ and in particular not under $\rho^m_j$ and
$\rho^n_i$ for $m, n \in \Z$.

\def\zovera{\mathop{\lower 10pt \hbox{$ {\buildrel{\displaystyle\bar{z}} \over
{\scriptstyle{(a)}}} $}}{\lower 4pt
\hbox{${\scriptstyle{ij}}$}}\nolimits}
\def\Zovera{\mathop{\lower 10pt \hbox{$ {\buildrel{\displaystyle\bar{Z}} \over
        {\scriptstyle{(a)}}} $}}{\lower 4pt \hbox{${\scriptstyle{ij}}$}}\nolimits}
\def\zoverka{\mathop{\lower 10pt \hbox{$ {\buildrel{\displaystyle\bar{z} } \over
        {\scriptstyle{(k)}}} $}}{\lower 4pt \hbox{${\scriptstyle{m \ k+1}}$}}\nolimits}
\def\zoverkb{\mathop{\lower 10pt \hbox{$ {\buildrel{\displaystyle\bar{z} }
\over
        {\scriptstyle{(k)}}} $}}{\lower 4pt \hbox{${\scriptstyle{k+1 \ m}}$}}\nolimits}
\def\zoverkc{\mathop{\lower 10pt \hbox{$ {\buildrel{\displaystyle\bar{z} } \over
        {\scriptstyle{(k)}}} $}}{\lower 4pt \hbox{${\scriptstyle{n \ k+2}}$}}\nolimits}
\def\zoverkd{\mathop{\lower 10pt \hbox{$ {\buildrel{\displaystyle{z} } \over
        {\scriptstyle{(k)}}} $}}{\lower 4pt \hbox{${\scriptstyle{k-1 \ k+2}}$}}\nolimits}
\def\zoverke{\mathop{\lower 10pt \hbox{$ {\buildrel{\displaystyle z } \over
        {\scriptstyle{(k)}}} $}}{\lower 4pt \hbox{${\scriptstyle{k-1 \ k+2}}$}}\nolimits}
\def\zoverkf{\mathop{\lower 10pt \hbox{$ {\buildrel{\displaystyle z} \over
        {\scriptstyle{(k)(k+2)}}} $}}{\lower 4pt \hbox{${\scriptstyle{k-1 \ k+4}}$}}\nolimits}
\def\zoverkg{\mathop{\lower 10pt \hbox{$ {\buildrel{\displaystyle \bar{z}} \over
        {\scriptstyle{(k+2)}}} $}}{\lower 4pt \hbox{${\scriptstyle{k-2 \ k+1}}$}}\nolimits}
\def\zoverkh{\mathop{\lower 10pt \hbox{$ {\buildrel{\displaystyle \bar{z}} \over
        {\scriptstyle{(k)}}} $}}{\lower 4pt \hbox{${\scriptstyle{k-2 \ n-3+{i \over 2}}}$}}\nolimits}
\def\zoverki{\mathop{\lower 10pt \hbox{$ {\buildrel{\displaystyle z} \over
        {\scriptstyle{(k)}}} $}}{\lower 4pt \hbox{${\scriptstyle{k-1 \ k+2}}$}}\nolimits}
\def\zoverkj{\mathop{\lower 10pt \hbox{$ {\buildrel{\displaystyle z} \over
        {\scriptstyle{(k)}}} $}}{\lower 4pt \hbox{${\scriptstyle{k-1 \ k+2}}$}}\nolimits}

In the following three sections we compute $H_{V_1}$, $H_{V_4}$ and $H_{V_7}$.
The others have similar computations, which appear in [Am2, Chapter 3].
The resulting $H_{V_i}$ for $i=2,3,5,6,8,9$ with the tables of braids 
appear in [AmTe2].

\section{The computation of $H_{V_1}$}\label{sec:31} 

\indent
We follow Figure 6 when computing the local braid monodromy for $V_1$, 
denoted by
$\vp^{(0)}_1$, but first we need to compute its factors which come from singularities in
the neighbourhood of $V_1$.    We start with $\vp^{(6)}_1$.  Then we compute $S^{(5)}_1$
and $\vp^{(5)}_1$ in the neighbourhood of $V_1$.  Notice that $\vp^{(5)}_1 =
\vp^{(4)}_1 = \vp^{(3)}_1$ since the fourth and fifth regenerations do not change the neighbourhood of $V_1$.
 Finally we compute $\vp^{(2)}_1.$  \ $\vp^{(2)}_1 =  \vp^{(1)}_1 = \vp^{(0)}_1$ since the seventh and eighth regenerations act similarly on the neighbourhood of $V_1$.  Hence we get a factorized expression for $\vp^{(0)}_1$.

{\bf Theorem 8.1} 
{\em In the neighbourhood of $V_1$, the local braid monodromy of $S^{(6)}_1$ around $V_1$ is given by:\\
$\vp^{(6)}_1 = (\uZ^2_{4'i})_{i = 2,5,6,6'}  \cdot Z^4_{34} \cdot
(\uZ^2_{4i})^{Z^2_{34}}_{i = 2,5,6,6'} \cdot \uZ^4_{14} \cdot
(Z_{44'})^{Z^2_{34}\uZ^2_{14}} \cdot Z^4_{56} \cdot (\uZ^2_{i6})^{Z^2_{i2}Z^2_{34}}_{i = 1,3} \cdot \uZ^4_{26} \cdot
(\uZ^2_{i6'})^{\uZ^2_{i6}\uZ^2_{i5}Z^2_{i2}Z^2_{34}}_{i = 1,3} \cdot
(Z_{66'})^{Z^2_{56}\uZ^2_{26}} \cdot (\Delta^2<1,2,3,5>)^{Z^2_{56}Z^2_{34}}$
 , where
$ (Z_{44'})^ {Z^2_{34}\uZ^2_{14}}$  is determined by\\
  $\vcenter{\hbox{\epsfbox{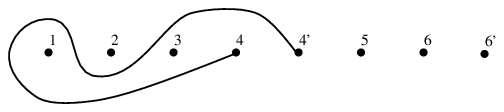}}}$,  $(Z_{66'})^{Z^2_{56} \uZ^2_{26}}$
is determined by $\vcenter{\hbox{\epsfbox{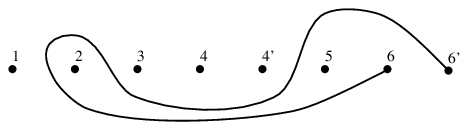}}}$   and the} {\rm   L.V.C.}
{\em corresponding to
$\Delta^2<1,2,3,5>$ is given by 
$\vcenter{\hbox{\epsfbox{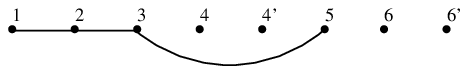}}}$}.

\underline{Proof:}
 Let $\hat{L}_i, 1 \le i \le 6$,
 be the lines in $Z^{(8)}$ and in $Z^{(7)}$ intersecting in a point 1.
  Let $L_i = \pi^{(8)} (\hat{L}_i)$.  Then $L_i, 1 \le i \le 6$
, intersect in $V_1$.  $S^{(8)}$ (and $S^{(7)}$) around $V_1$ is
 $\bigcup\limits^6_{i=1} L_i$.  After the second regeneration, $S^{(6)}_1$ in a
 neighbourhood of $V_1$ is as follows:  the lines $L_4$ and $L_6$ are replaced by
  two conics $Q_4$ and $Q_6$, where $Q_4$ is tangent to the lines $L_1$ and $L_3$,
  $Q_6$ is tangent to $L_2$ and $L_5$ .
 This follows from the fact that $P_2 \cup P_3$ (resp. $P_{13} \cup P_{15})$
  is regenerated to $R_1$ (resp. $R_6)$, and the tangency follows from Lemma 
\ref{lm:24}.  Thus, in a neighbourhood of $V_1$, $S^{(6)}_1$ is of the form: $S^{(6)}_1 = L_1 \cup L_2 \cup L_3 \cup Q_4 \cup L_5 \cup Q_6$.

$S^{(6)}_1$ is shown in Figure 12.
\begin{figure}[!h]
\epsfxsize=11cm 
\epsfysize=9cm 
\begin{minipage}{\textwidth}
\begin{center}
\epsfbox {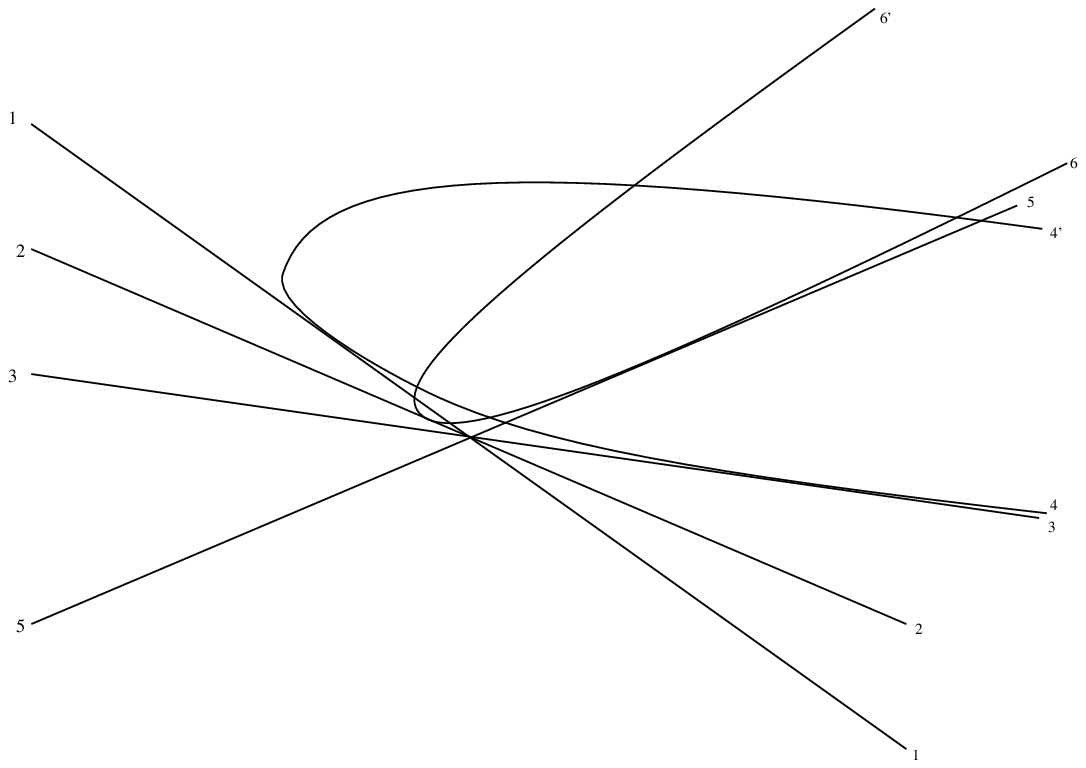}
\end{center}
\end{minipage}
\caption{}
\end{figure}

In order to prove the theorem we devide the proof to parts and we present the following propositions and remarks.

\noindent
{\bf Proposition 8.1.1.} {\em Consider in $\C^2$ the following figure:  $C_1 = L_1 \cup L_3 \cup Q_4 \cup L_5$, where $L_i$ are lines, $i = 1,3,5, Q_4$ is a conic tangent to $L_1$ and $L_3$, it intersects $L_5$ and the 3 lines meet at one point (see Figure 13).

Then the braid monodromy of $C_1$ is given by:\\
$\vp_{C_1} = Z^2_{4'5} \cdot Z^4_{34} \cdot (\uZ^2_{45})^{Z^2_{34}} \cdot
(\Delta^2<1,3,5>)^{Z^2_{34}} \cdot \uZ^4_{14} \cdot (Z_{44'})^{Z^2_{34}\uZ^2_{14}}$, where
$(Z_{44'})^{Z^2_{34}\uZ^2_{14}}$ is determined by
$\vcenter{\hbox{\epsfbox{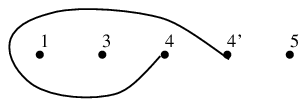}}}$ and the L.V.C. corresponding to $\Delta^2<1,3,5>$ is given by
{$\vcenter{\hbox{\epsfbox{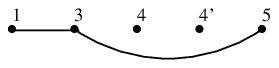}}}$}}.
\begin{figure}[h]\label{fig11'}
\epsfxsize=11cm 
\epsfysize=9cm 
\begin{minipage}{\textwidth}
\begin{center}
\epsfbox {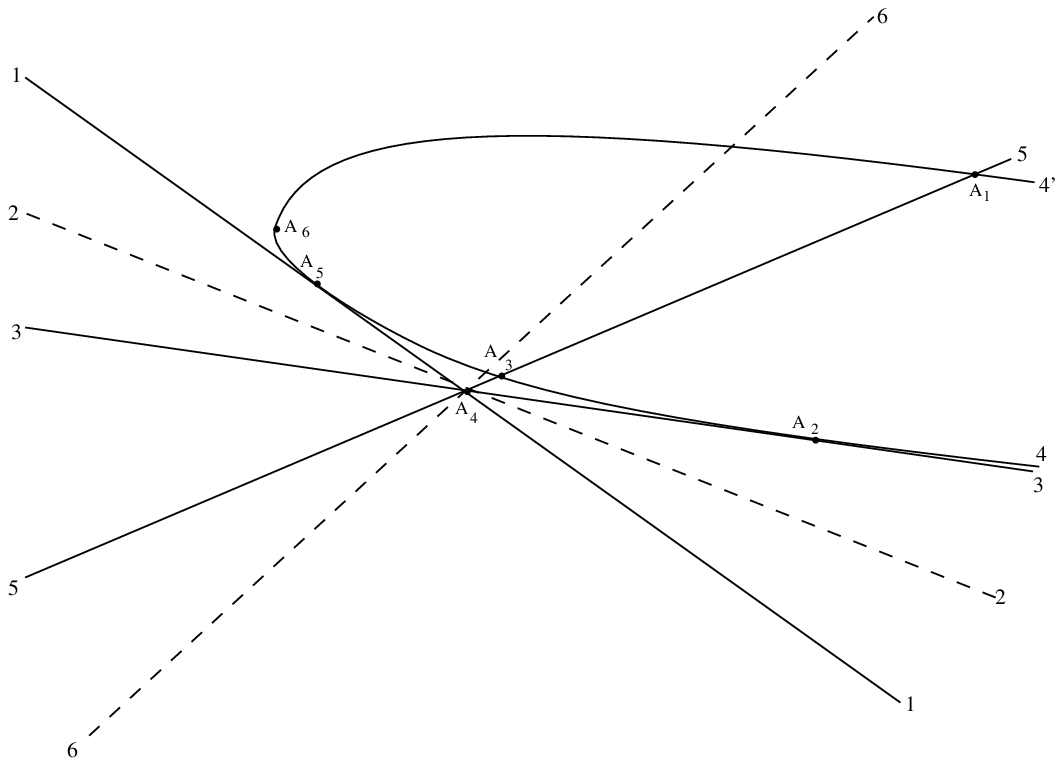}
\end{center}
\end{minipage}
\caption{}
\end{figure}

\underline{Proof:}
Let $\pi_1 : E \times D \rightarrow E$ be the projection to $E$.

Let $\{A_j\}^{6}_{j=1}$ be  singular points of $\pi_1$ as follows:\\
$A_1, A_3 $ are intersection points of $Q_4 $ with the line $L_5$.\\
$A_2,A_5$ are tangent points of $Q_4$ with the lines $L_3, L_1$ respectively. \\
$A_4$ is an intersection  point of the  lines $L_1, L_3, L_5$.\\
$A_6$ is a point of the type $a_1$ in $Q_4$.

Let $N = \{ x(A_j) = x_j \mid 1 \le j \le 6 \}$, \st $N \subset E - \partial E, N \subset E_\R$.  Let $M$ be a real point on the $x$-axis, \st $x_j <<M \;\; \forall x_j \in N, 1 \le j \le 6$.

We claim that there exists a $g$-base $\{\ell(\g_j)\}^{6}_{j=1}$ of $\pi_1(E - N,u)$, \st each path $\g_j$ is below the real line and the values of $\vp_M$ w.r.t. this base and $E \times D$ are the ones given in the proposition.

Recall that $K = K(M)$.  Let $i = L_i \cap K$ for $i = 1,3,5$.
 Let $\{4,4'\} = Q_4 \cap K$.
 $K = \{1,3,4,4',5\}$, \st 1,3,4,4',5 are  real points.
Moreover: $1 < 3 < 4< 4' < 5$.
Let $\beta_M$ be the diffeomorphism defined by: $\beta_M(1) = 1,\ \beta_M(3) =2,\ \beta_M(4) = 3, \ \beta_M(4') = 4, \ \beta_M (5) = 5$. Moreover: deg$C_1 = 5$ and $\# K_\R (x) = 5 \ \forall x$.

We are looking for $\vp_M(\ell (\g_j))$ for $j = 1, \cdots , 6$.  
We choose a $g$-base $\{\ell(\g_j)\}^{6}_{j=1}$ of $\pi_1(E - N,u)$, \st each path $\g_j$ is below the real line.

Put all data in the following table:
\vspace{-1cm}
\begin{center}
\begin{tabular}{cccc} \\
$j$ & $\lambda_{x_j}$ & $\epsilon_{x_j}
$ & $\delta_{x_j}$ \\ \hline
1 & $<4,5>$ & 2 & $\Delta<4,5>$\\
2 & $<2,3>$ & 4 & $\Delta^2<2,3>$\\
3 & $<3,4>$ & 2 & $\Delta<3,4>$\\
4 & $<1,3>$ & 2 & $\Delta<1,3>$\\
5 & $<3,4>$ & 4 & $\Delta^2<3,4>$\\
6 & $<4,5>$ & 1 & $\Delta^{\frac{1}{2}}_{I_2\R}<4>$\\
\end{tabular}
\end{center}

For computations, we use the formulas in [Am2, Theorems 1.41,1.44].

\noindent
\underline{Remark}:  $\beta_{x'_j} (K(x'_j)) = \{1,2,3,4,5\}$ for $1 \leq j \leq 6$.\\

\noindent
$(\xi_{x'_1}) \Psi_{\g'_1} = <4,5> \beta^{-1}_M = z_{4'5}$\\
$<4,5>$ \hspace{.5cm}$\vcenter{\hbox{\epsfbox{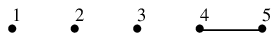}}}$ \hspace{.4cm}
$\beta^{-1}_M$\hspace{.5cm}  $\vcenter{\hbox{\epsfbox{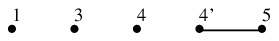}}}$ \hspace{.5cm}
$\vp_M(\ell(\g_1)) = Z^2_{4'5}$ \\

\medskip

\noindent
$(\xi_{x'_2}) \Psi_{\g'_2} = <2,3> \Delta<4,5>\beta^{-1}_M = z_{34}$\\
$<2,3> \Delta <4,5>$   \ \ $\vcenter{\hbox{\epsfbox{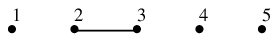}}}$ \ \ $\beta^{-1}_M$ \ \  $\vcenter{\hbox{\epsfbox{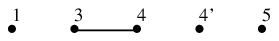}}}$ \ \ $\vp_M(\ell(\g_2)) = Z^4_{34}$ \\

\medskip
\noindent
$(\xi_{x'_3}) \Psi_{\g'_3}
 = <3,4> \Delta^2<2,3> \Delta<4,5> \beta^{-1}_M =
\underline{z}_{45}^{Z^2_{34}}$

\noindent
$<3,4> \vcenter{\hbox{\epsfbox{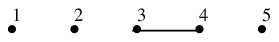}}}$ \ \
$\Delta^2<2,3> \vcenter{\hbox{\epsfbox{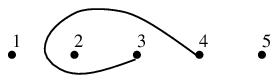}}}$
\ \ $\Delta<4,5> \vcenter{\hbox{\epsfbox{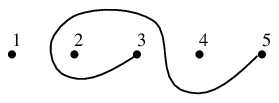}}}$ \\
$\beta^{-1}_M$ $\vcenter{\hbox{\epsfbox{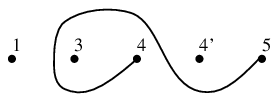}}}$ \ \
$\vp_M(\ell(\g_3)) = (\uZ^2_{45})^{Z^2_{34}}$\\

\medskip

\noindent
$(\xi_{x'_4}) \Psi_{\g'_4}  = <1,3>  \Delta<3,4>\Delta^2<2,3> \Delta<4,5>\beta^{-1}_M =
(\Delta<1,3,5>)^{Z^2_{34}}$\\
$<1,3> \vcenter{\hbox{\epsfbox{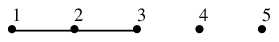}}}$ \ \
$\Delta<3,4> \  \
\vcenter{\hbox{\epsfbox{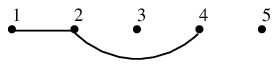}}}$ \ \
$ \Delta^2<2,3> \ \
 \vcenter{\hbox{\epsfbox{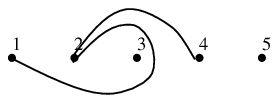}}}$ \\
$ \Delta<4,5>
 \vcenter{\hbox{\epsfbox{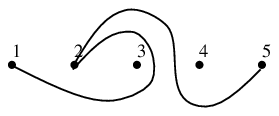}}}$
\ \
$\beta^{-1}_M$ $\vcenter{\hbox{\epsfbox{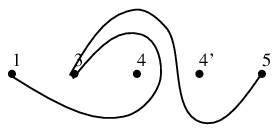}}}$ \ \
$\vp_M(\ell(\g_4)) =  (\Delta^2<1,3,5>)^{Z^2_{34}}$\\

\medskip

\noindent
$(\xi_{x'_5}) \Psi_{\g'_5}
 = <3,4> \Delta<1,3> \Delta<3,4>\Delta^2<2,3> \Delta<4,5> \beta^{-1}_M =
\underline{z}_{14}$\\
$<3,4>$  $\vcenter{\hbox{\epsfbox{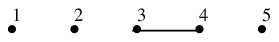}}}$ \ \ $\Delta<1,3>$
 $\vcenter{\hbox{\epsfbox{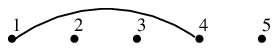}}}$ \ \
$\Delta<3,4>$
$\vcenter{\hbox{\epsfbox{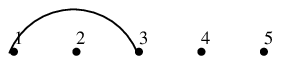}}}$ \\
$\Delta^2<2,3>\Delta<4,5>$
 $\vcenter{\hbox{\epsfbox{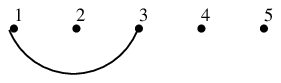}}}$ \ \
$\beta^{-1}_M$
$\vcenter{\hbox{\epsfbox{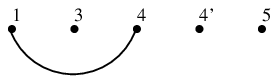}}}$ \ \
$\vp_M(\ell(\g_5)) = \uZ^4_{14}$\\

\medskip
\noindent
$(\xi_{x'_6}) \Psi_{\g'_6}
 = <4,5> \Delta^2<3,4>\Delta<1,3>  \Delta<3,4>\Delta^2<2,3> \Delta<4,5>\beta^{-1}_M = z^{Z^2_{34}\uZ^2_{14}}_{44'}$\\
\noindent
$<4,5>$ $\vcenter{\hbox{\epsfbox{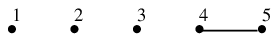}}}$ \ \
$\Delta^2<3,4>$ $\vcenter{\hbox{\epsfbox{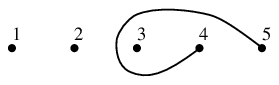}}}$ \ \
$\Delta<1,3>$ $\vcenter{\hbox{\epsfbox{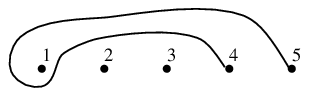}}}$\\
$\Delta<3,4>$
$\vcenter{\hbox{\epsfbox{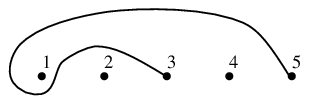}}}$ \ \
$\Delta^2<2,3>$
$\vcenter{\hbox{\epsfbox{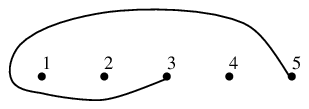}}}$ \ \
$\Delta<4,5>$
$\vcenter{\hbox{\epsfbox{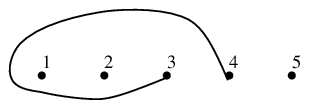}}}$\\
$\beta^{-1}_M$
$\vcenter{\hbox{\epsfbox{1v025.ps}}}$ \ \
$\vp_M(\ell(\g_6)) = (Z_{44'})^{Z^2_{34}\uZ^2_{14}}$

{\bf Remark 8.1.2.} {\em Since $C_1$ is in $\C^2 \subset \C \P^2$, we shall use a slight modification of Proposition 8.1.1 (using an embedding to the projective plane).  In the modification there are 2 more lines $L_2, L_6$.  $L_2$ intersects the conic $Q_4$ in 2 complex points.

Following the proof of Proposition 8.1.1, one can easily see that the formulation of the result is very similar to the formulation of Proposition 8.1.1.  The changes in the result are as follows:\\
(a) $(\Delta^2<1,3,5>)^{Z^2_{34}}$ is replaced by
$(\Delta^2<1,2,3,5,6>)^{Z^2_{34}}$ with the corresponding L.V.C.
which appears as follows:
$\vcenter{\hbox{\epsfbox{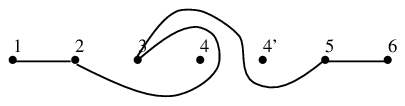}}}$ .\\
(b) $Z^2_{4'5}$ is replaced by $(\uZ^2_{4'i})_{i = 2,5,6}$ as follows:
 $\vcenter{\hbox{\epsfbox{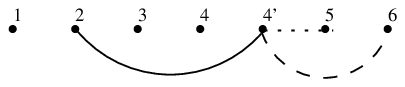}}}$ .\\
(c) $(\uZ^2_{45})^{Z^2_{34}}$ is replaced by $(\uZ^2_{4i})_{i = 2,5,6}^{Z^2_{34}}$ as follows:
 $\vcenter{\hbox{\epsfbox{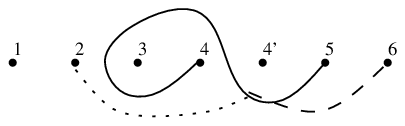}}}$}.

\noindent
{\bf Proposition 8.1.3.} {\em Consider in $\C^2$ the following figure: $C_2 = L_1 \cup L_2 \cup L_5 \cup Q_6$, where $L_i$ are lines, $i = 1,2,5, Q_6$ is a conic tangent to $L_2$ and $L_5$, it intersects $L_1$ and the 3 lines meet at one point (see Figure 14).

Then the braid monodromy of $C_2$ is given by:\\
$\vp_{C_2} = Z^4_{56} \cdot (\Delta^2<1,2,5>)^{Z^2_{56}}
\cdot (\uZ^2_{16})^{Z^2_{12}}\cdot \uZ^4_{26}\cdot
(\uZ^2_{16'})^{\uZ^2_{16}\uZ^2_{15}Z^2_{12}}
\cdot (Z_{66'})^{Z^2_{56}\uZ^2_{26}}$, where
$(Z_{66'})^{Z^2_{56}\uZ^2_{26}}$  is determined by
$\vcenter{\hbox{\epsfbox{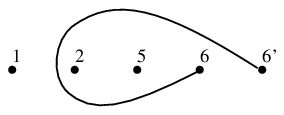}}}$   and the {\rm   L.V.C.}  
corresponding to
$\Delta^2<1,2,5>$ is determined by
$\vcenter{\hbox{\epsfbox{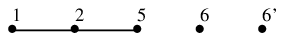}}}$}.

\begin{figure}[h]
\epsfxsize=11cm 
\epsfysize=9cm 
\begin{minipage}{\textwidth}
\begin{center}
\epsfbox {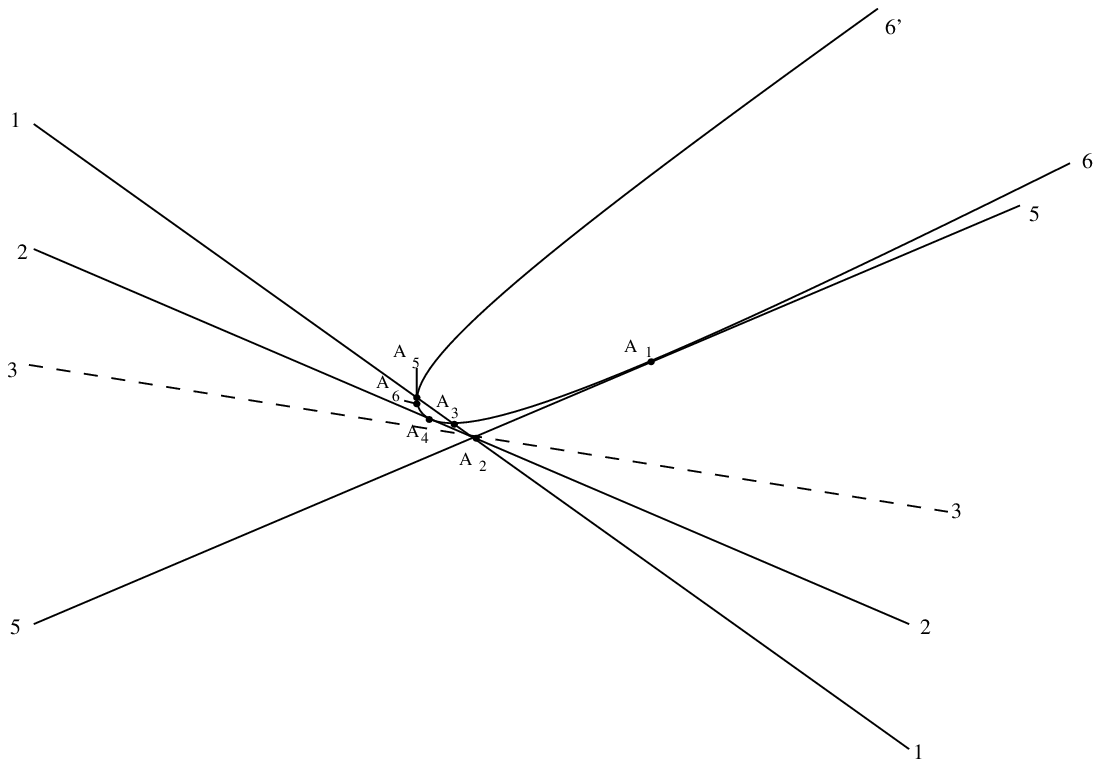}
\end{center}
\end{minipage}
\caption{}
\end{figure}
\underline{Proof}: Let $\pi_1: E \times D \rightarrow E$ be the projection to $E$.

Let $\{A_j\}^{6}_{j=1}$ be singular points of $\pi_1$ as follows:\\
$A_1, A_4$ are tangent points of $Q_6$ with the lines $L_5, L_2$ respectively.\\
$A_3, A_5$ are intersection points of $Q_6$ with the line $L_1$.\\
$A_2$ is an intersection point of the lines $L_1, L_2, L_5$.\\
$A_6$ is a point of the type $a_1$ in $Q_6$.

Let $N = \{x(A_j) = x_j \mid 1 \leq j \leq 6 \}$, \st $N \subset E - \partial E, N \subset E_\R$.  Let $M$ be a real point on the $x$-axis, \st $x_j << M \; \; \forall x_j \in N, 1 \leq j \leq 6$.

We claim that there exists a $g$-base $\{\ell (\g_j)\}^{6}_{j=1}$ of $\pi_1(E - N,u)$, \st each path $\g_j$ is below the real line and the values of $\vp_M$ w.r.t. this base and  $E \times D$ are the ones given in the proposition.

Recall that $K = K(M)$.  Let $i = L_i \cap K$ for $i = 1,2,5$.  Let
$\{6,6'\} = Q_6 \cap K$. $K = \{1,2,5,6,6'\}$, \st $1,2,5,6,6'$ are real points.  Moreover: $1 < 2 < 5 < 6 < 6'$.  Let $\beta_M$ be the diffeomorphism defined by $\beta_M(1) = 1, \ \beta_M(2) = 2, \ \beta_M(5) = 3, \ \beta_M(6) = 4, \ \beta_M(6') = 5$.  Moreover: deg$C_2 = 5$ and
 $\# K_\R (x) = 5 \ \forall x$.

We are looking for $\vp_M(\ell (\g_j))$ for $j = 1, \cdots , 6$.  We choose a $g$-base $\{\ell(\g_j)\}^{6}_{j=1}$ of $\pi_1(E - N,u)$, \st each path $\g_j$ is below the real line.

Put all data in the following table:
\vspace{-1cm}
\begin{center}
\begin{tabular}{cccc} \\
$j$ & $\lambda_{x_j}$ & $\epsilon_{x_j}
$ & $\delta_{x_j}$ \\ \hline
1 & $<3,4>$ & 4 & $\Delta^2<3,4>$\\
2 & $<1,3>$ & 2 & $\Delta<1,3>$\\
3 & $<3,4>$ & 2 & $\Delta<3,4>$\\
4 & $<2,3>$ & 4 & $\Delta^2<2,3>$\\
5 & $<4,5>$ & 2 & $\Delta<4,5>$\\
6 & $<3,4>$ & 1 & $\Delta^{\frac{1}{2}}_{I_2\R}<3>$\\
\end{tabular}
\end{center}

For computations, we use the formulas in [Am2, Theorems 1.41, 1.44].

\underline{Remark}:  $\beta_{x'_j} (K(x'_j)) = \{1,2,3,4,5\}$ for $1 \leq j \leq 6$.\\

\noindent
$(\xi_{x'_1}) \Psi_{\g'_1} = <3,4> \beta^{-1}_M = z_{56}$\\
$<3,4>$ \hspace{.5cm}$\vcenter{\hbox{\epsfbox{1v005.ps}}}$ \hspace{.4cm}
$\beta^{-1}_M$\hspace{.5cm}  $\vcenter{\hbox{\epsfbox{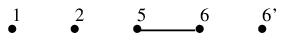}}}$ \hspace{.5cm}
$\vp_M(\ell(\g_1)) = Z^4_{56}$ \\

\medskip

\noindent
$(\xi_{x'_2}) \Psi_{\g'_2} = <1,3> \Delta^2<3,4>\beta^{-1}_M = (\Delta<1,2,5>)^{Z^2_{56}}$\\
$<1,3>   \vcenter{\hbox{\epsfbox{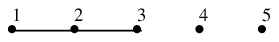}}}$ \ \
$\Delta^2<3,4>  \ \vcenter{\hbox{\epsfbox{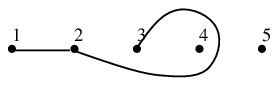}}}$ \ \
$\beta^{-1}_M$ \ \ $\vcenter{\hbox{\epsfbox{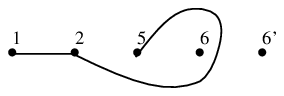}}}$\\
$\vp_M(\ell(\g_2)) = (\Delta^2<1,2,5>)^{Z^2_{56}}$\\

\medskip

\noindent
$(\xi_{x'_3}) \Psi_{\g'_3}
 = <3,4> \Delta<1,3>   \Delta^2 <3,4>\beta^{-1}_M =
\underline{z}^{Z^2_{12}}_{16}$\\
$<3,4>$  $\vcenter{\hbox{\epsfbox{1v005.ps}}}$ \ \
$\Delta<1,3>$  $\vcenter{\hbox{\epsfbox{1v015.ps}}}$ \ \
$\Delta^2<3,4>$  $\vcenter{\hbox{\epsfbox{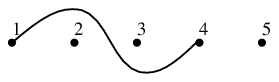}}}$ \\
\noindent
$\beta^{-1}_M$ $\vcenter{\hbox{\epsfbox{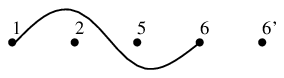}}}$ \ \
$\vp_M(\ell(\g_3)) = (\uZ_{16}^2)^{Z^2_{12}}$\\

\medskip

\noindent
$(\xi_{x'_4}) \Psi_{\g'_4}  = <2,3> \Delta<3,4>\Delta<1,3> \Delta^2<3,4>
\beta^{-1}_M = \underline{ z}_{26}$\\
$<2,3>$  $\vcenter{\hbox{\epsfbox{1v003.ps}}}$ \ \
$\Delta<3,4>$ \
 $\vcenter{\hbox{\epsfbox{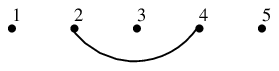}}}$ \ \
$\Delta<1,3>$ \ \
 $\vcenter{\hbox{\epsfbox{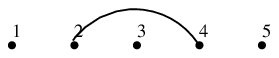}}}$ \\
$\Delta^2<3,4>$ \ \
 $\vcenter{\hbox{\epsfbox{1v035.ps}}}$ \ \
$\beta^{-1}_M$ $\vcenter{\hbox{\epsfbox{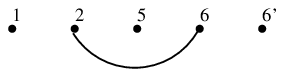}}}$ \ \
$\vp_M(\ell(\g_4)) = \uZ^4_{26}$\\

\medskip

\noindent
$(\xi_{x'_5}) \Psi_{\g'_5}
 = <4,5> \Delta^2<2,3>  \Delta<3,4>\Delta<1,3> \Delta^2<3,4> \beta^{-1}_M =
\underline{ z}^{\uZ^2_{16}\uZ^2_{15}Z^2_{12}}_{16'}$\\
$<4,5> \Delta^2<2,3>$   $\vcenter{\hbox{\epsfbox{1v001.ps}}}$ \ \
$\Delta<3,4>$
 $\vcenter{\hbox{\epsfbox{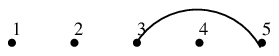}}}$ \\
$\Delta<1,3> \Delta^2<3,4>$
 $\vcenter{\hbox{\epsfbox{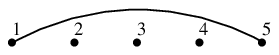}}}$\ \
$\beta^{-1}_M$
$\vcenter{\hbox{\epsfbox{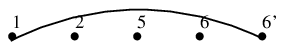}}}$ \ \
$\vp_M(\ell(\g_5)) =  (\uZ^2_{16'})^{\uZ^2_{16}\uZ^2_{15}Z^2_{12}}$\\

\medskip
\noindent
$(\xi_{x'_6}) \Psi_{\g'_6}
 = <3,4> \Delta<4,5>\Delta^2<2,3> \Delta<3,4>\Delta<1,3> \Delta^2<3,4> \beta^{-1}_M = z_{66'}^{Z^{2}_{56}\uZ^2_{26}}$

\medskip


\noindent
$<3,4>$ $\vcenter{\hbox{\epsfbox{1v005.ps}}}$ \ \
$\Delta<4,5>$ $\vcenter{\hbox{\epsfbox{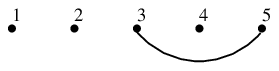}}}$ \ \
$\Delta^2<2,3>$ $\vcenter{\hbox{\epsfbox{1v007.ps}}}$\\
\noindent
$\Delta<3,4>$ $\vcenter{\hbox{\epsfbox{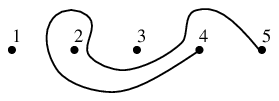}}}$ \ \
$\Delta<1,3>$ $\vcenter{\hbox{\epsfbox{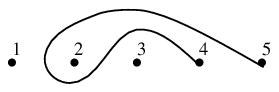}}}$ \ \
$\Delta^2<3,4>$ $\vcenter{\hbox{\epsfbox{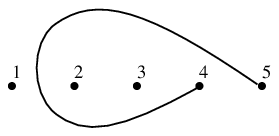}}}$\\
$\beta^{-1}_M$
$\vcenter{\hbox{\epsfbox{1v045.ps}}}$ \ \
$\vp_M(\ell(\g_6)) =  (Z_{66'})^{Z^{2}_{56}\uZ^2_{26}}$

\medskip
\noindent
{\bf Remark 8.1.4.} {\em Since $C_2$ is in $\C^2 \subseteq \C \P^2$, we shall use a slight modification of Proposition 8.1.3.  In the modification there is one more line $L_3$ which intersects $Q_6$ in 2 complex points.

Following the proof of Proposition 8.1.3, one can easily see that the formulation of the result is very similar to the formulation of Proposition 8.1.3.  The changes in the result are as follows:\\
(a) $(\Delta^2<1,2,5>)^{Z^2_{56}}$ is replaced by
$(\Delta^2<1,2,3,5>)^{Z^2_{56}}$ with the corresponding L.V.C.
which appears as follows:
$\vcenter{\hbox{\epsfbox{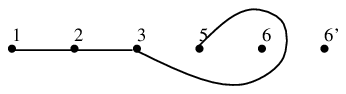}}}$. \\
(b) $(\uZ^2_{16})^{Z^2_{12}}$ is replaced by $(\uZ^2_{i6})_{i = 1,3}^{Z^2_{i2}}$ as follows: \  $\vcenter{\hbox{\epsfbox{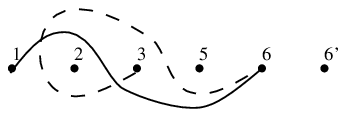}}}$ .\\
(c) $(\uZ^2_{16'})^{\uZ^2_{16}\uZ^2_{15}Z^2_{12}}$ is replaced by
$(\uZ^2_{i6'})_{i = 1,3}^{\uZ^2_{i6}\uZ^2_{i5}Z^2_{i2}}$ as follows: \
 $\vcenter{\hbox{\epsfbox{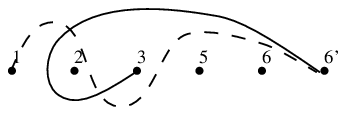}}}$ .}\\
\noindent
\underline{Proof of Theorem 8.1:}  Until now we computed the braid monodromy of all singularities (real and complex).  Each one of the intersection points of $L_6 \cap Q_4$ is replaced by 2 intersection points $(\subseteq Q_4 \cap Q_6)$ which are close to each other.  Take $M$ on the $x-axis$, \st $\pi^{-1} (M) \cap S^{(6)}_1 = \{1,2,3,4,4',5,6,6'\}$.

To build a part of the desired $g$-base that corresponds to singularities from Proposition 8.1.1 and Remark 8.1.2, we take all elements of the $g$-base that we constructed for Remark 8.1.2 (except for the one that corresponds to $x(\cap L_i))$ and we make on them the following changes:\\
(I) We replace each of the paths that corresponds to one of the points $Q_4 \cap L_6$ by 2 paths which follow former paths almost up to the end and then form a bush with 2 branches, see Figure 15.
\begin{figure}[h]
\begin{minipage}{\textwidth}
\begin{center}
\epsfbox {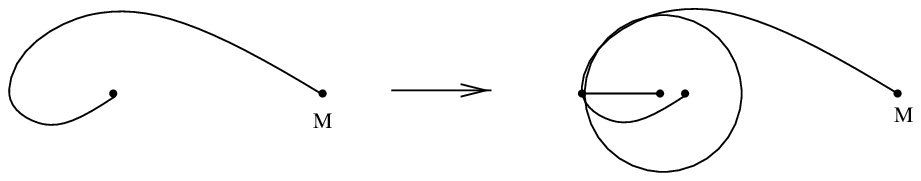}
\end{center}
\end{minipage}
\caption{}
\end{figure}

(II) For the loops that we multiplied we get the following braid monodromies: $\uZ^2_{4'6}$ is replaced by 2 elements $\uZ^2_{4'6}$ and $\uZ^2_{4'6'}$.  These
braids correspond to the paths $\underline{ z}_{4'6}$ and $\underline{ z}_{4'6'}$ respectively.  $(\uZ^2_{46})^{Z^2_{34}}$ is replaced by $(\uZ^2_{46})^{Z^2_{34}}$ and $(\uZ^2_{46'})^{Z^2_{34}}$.  These braids correspond to the paths $\underline{ z}_{46}^{Z^2_{34}}$ and $\underline{ z}^{Z^2_{34}}_{46'}$ respectively.\\
(III) Notation: $x(\cap L_i) = x_0$.\\
We choose a circular neighbourhood $E'$ of $x_0$ on the $x$-axis that includes all $x$-projections of singularities from Proposition 8.1.3 and Remark 8.1.4.  Let $M' = \max\{\partial E' \cap \R\}$.  Let $S' = \pi^{-1} (E') \cap (S^{(6)}_1 - Q_4)$.  To build the part of the desired $g$-base that corresponds to these
singularities, we start by building a $g$-base for $(E'-N,M)$ and computing
its braid monodromy.  We can apply Proposition 8.1.3 to Figure 8.3 restricted to $(E',M',\pi)$.

Thus we get a $g$-base in $(E'-N, M')$ whose braid monodromy is  given in Proposition 8.1.3, since a $g$-base is determined up to homotopical equivalence.  We can assume that the path $\g_{0}$ from $M$ to $x_0$, which is an element of the $g$-base of Remark 8.1.2 is entering $E'$ at the point $M'$.  We extend the $g$-base of $(E'-N,M')$ obtained above until $M$ by adding to it the part
$\g''_0$ of $\g_0$ outside of $E'$.  The braids of the prolonged loops under the braid monodromy of $S^{(6)}_1$ at $M$ are obtained from the braids of the loops in
$(E'-N,M')$ under the braid monodromy of $S'$ at $M'$, by applying on them a natural monomorphism $e: B[D,K(M')\cap S'] \rightarrow B[D,K(M')]$ followed by
$(\psi_{\g''_0})^\vee : B[D,K(M')] \rightarrow B[D,K(M)]$.  As always,
$(\psi_{\g''_0})^\vee \circ e$ is determined by its values on a linear frame $H = H(\sigma'_i)^5_{i=1}$ in  $B[D,K(M')\cap S']$ where $\xi_{M'} = \{\sigma'_i\}$ is a skeleton in $(D,K(M') \cap S')$.
$K(M') \cap S' = \{1*,2*,3*,5*,6*,6'*\}$,
$\xi_{M'} = \{ \sigma'_i\}^5_{i=1} = \vcenter{\hbox{\epsfbox{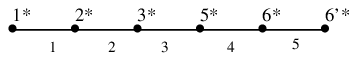}}}$.\\
Let us denote by $\{ \sigma_i\}^5_{i=1}$ the following skeleton in $(D,K(M),K(M)\cap S')$:\\
$\xi_{M} = \{ \sigma_i\}^5_{i=1} = \vcenter{\hbox{\epsfbox{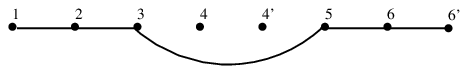}}}$ .

In Proposition 8.1.1 and Remark 8.1.2 we computed $(\xi_{x'_0})\psi_{\g'_0}$, where $\xi_{x'_0}$ is a skeleton in $(D, K(x'_0))$.  Let us write $\g'_{0} = \g'''_{0} \cup \g''_{0}$.  It is obvious from Remark 8.1.2 that $(\xi_{x'_0})
\psi_{\g'''_{0}} = \xi_{M'}$ and thus $(\xi_{M'}) \psi_{\g''_0}= (\xi_{x'_0})\psi_{\g'_0}$.  The last skeleton is computed in Proposition 8.1.1 and it equals
 $<(\sigma_i)^{Z^2_{34}}>$.

Thus $(\xi_{M'}) \psi_{\g''_0} = < (\sigma_i)^{Z^2_{34}}> =
\vcenter{\hbox{\epsfbox{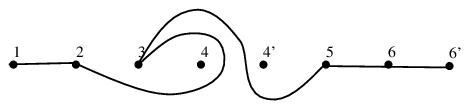}}} = \xi^{Z^2_{34}}_{M'}$.  Therefore, we have $(\psi_{\g'_0})^\vee \circ e(H(\sigma'_i)) = H((\sigma'_i)\psi_{\g''_0}) = \left\{ \begin{array}{ll} H(\sigma_i) & i = 1,4, 5 \\
(H(\sigma_i))^{Z^2_{34}} & i = 2,3. \end{array} \right.$

We conclude that there exist elements of the $g$-base of $(E-N,u)$ for $S^{(6)}_1$ corresponding to the singularities from Proposition 8.1.3 and Remark 8.1.4 whose braid monodromies are in conjugation with $Z^2_{34}$ of the corresponding
 braids for $S'$ and $(E', \pi,M')$.  Similarly, the related L.V.C.'s are obtained from the corresponding L.V.C.'s for $(S', E', \pi, M')$ by applying the braid $Z^2_{34}$.

Considering the above changes, we present here the list of all braids:

\begin{tabbing}
ggggffgfgfffffffgggvvddvvvvvv \= dggd corresponding to
$ (\uZ^2_{4'i})^2_{i = 2,5,6,6'}$ \kill
$\vcenter{\hbox{\epsfbox{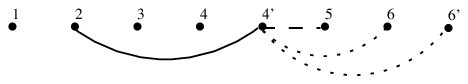}}}$ \> corresponding to
$ (\uZ^2_{4'i})_{i = 2,5,6,6'}$ \\
$\vcenter{\hbox{\epsfbox{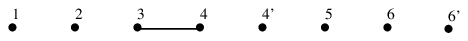}}}$ \> corresponding to
$Z^4_{34} $\\
$\vcenter{\hbox{\epsfbox{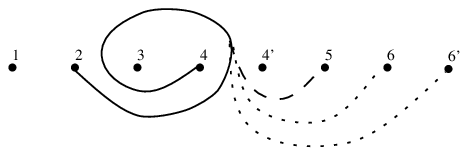}}}$ \> corresponding to
$(\uZ^2_{4i})^{Z^2_{34}}_{i = 2,5,6,6'}$\\ [.4cm]
$\vcenter{\hbox{\epsfbox{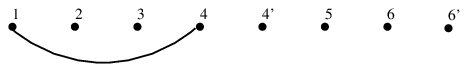}}}$ \> corresponding to
$\uZ^4_{14} $\\
$\vcenter{\hbox{\epsfbox{v101.ps}}}$ \> corresponding to
$(Z_{44'})^{Z^2_{34}\uZ^2_{14}} $\\
$\vcenter{\hbox{\epsfbox{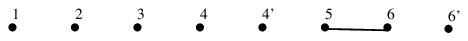}}}$ \> corresponding to
$Z^4_{56}$\\[.4cm]
$\vcenter{\hbox{\epsfbox{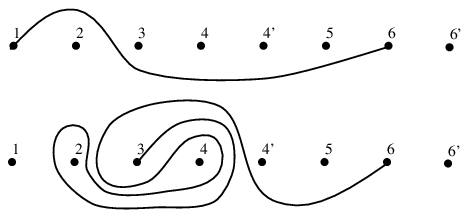}}}$ \>
 corresponding to $(\uZ^2_{i6})^{Z^2_{i2}Z^2_{34}}_{i = 1,3}$\\[.4cm]
$\vcenter{\hbox{\epsfbox{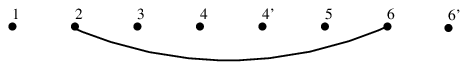}}}$ \> corresponding to $\uZ^4_{26}$\\[.4cm]
$\vcenter{\hbox{\epsfbox{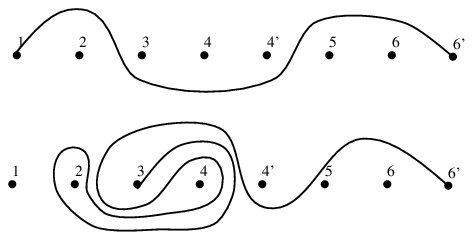}}}$ \> corresponding to $(\uZ^2_{i6'})^{\uZ^2_{i6}\uZ^2_{i5}Z^2_{i2}Z^2_{34}}_{i = 1,3}$\\[.4cm]
$\vcenter{\hbox{\epsfbox{v102.ps}}}$ \> corresponding to $(Z_{66'})^{Z^2_{56}\uZ^2_{26}}$\\
$\vcenter{\hbox{\epsfbox{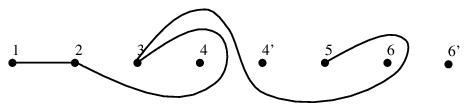}}}$ \> corresponding to $(\Delta^2<1,2,3,5>)^{Z^2_{56}Z^2_{34}}$ \hspace{3cm}$\Box$
\end{tabbing}

Now we compute $\vp^{(5)}_1$.  Recall that $\vp^{(5)}_1 = \vp^{(4)}_1 =\vp^{(3)}_1$, so in fact we compute $\vp^{(3)}_1$.

\noindent
{\bf Theorem 8.2.} {\em Let $\omega = L_1 \cup L_2 \cup L_3 \cup L_5$ in $S^{(6)}_1$. \\
(a) In a small neighbourhood $U_3$ of $\omega, T^{(3)} = S^{(3)}_1 \cap \ U_3$ resembles Figure 17, i.e., the singularities of $T^{(3)}$ are 4 nodes, 4 tangency points and 2 branch points.\\
(b) The local braid monodromy of $S^{(3)}_1$ in that neighbourhood is presented by $F_1(F_1)^{\rho^{-1}}$, where:\\
 $F_1 = Z^4_{23} \cdot Z^2_{3'5} \cdot (\uZ^2_{35})^{Z^2_{23}} \cdot
(\uZ^4_{25})^{Z^2_{23}}\cdot (Z_{12})^{\uZ^2_{25}Z^2_{23}} .\\
 \rho = Z_{33'} Z_{55'}$ and $(z_{12})^{\uZ^2_{25} Z^2_{23}} = $}
$\vcenter{\hbox{\epsfbox{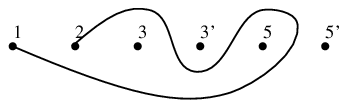}}}$ .

\noindent
\underline{Proof}:

(a) Let $\hat{\omega} = R_1 \cap R_2 \cap R_4 \cap R_6$ in $Z^{(6)}$.

For local analysis we can use holomorphic coordinates in a neighbourhood $\hat{U}$ of $\hat{\omega}$ in $\C \P^n$.  This allows us to consider $\hat{U}$ as a subvariety of a neighbourhood of the origin in $\C^4$ with coordinates $X, Y, Z, T$ defined by the following system of equations
$
\hat{U} : \left\{ \begin{array}{l}
XT = 0 \\
YZ = 0 \end{array} \right .$.

By abuse of notation, $\hat{L}_5$  is called the $x$-axis (see Figure 16), and then: $R_1 = \{X=0, Y=0\}, R_2 = \{X=0, Z=0\}, R_4 = \{T=0, Y=0\}, R_6 = \{T=0, Z=0 \}$.
\begin{figure}[h]
\epsfxsize=5cm 
\epsfysize=5cm 
\begin{minipage}{\textwidth}
\begin{center}
\epsfbox {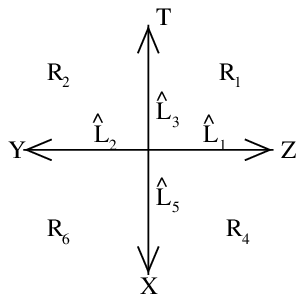}
\end{center}
\end{minipage}
\caption{}
\end{figure}

Consider the curve we obtain in [MoTe8, Lemma 6].  We can spin this curve 
to obtain Figure 17.
\begin{figure}[h]
\epsfxsize=11cm 
\epsfysize=9cm 
\begin{minipage}{\textwidth}
\begin{center}
\epsfbox {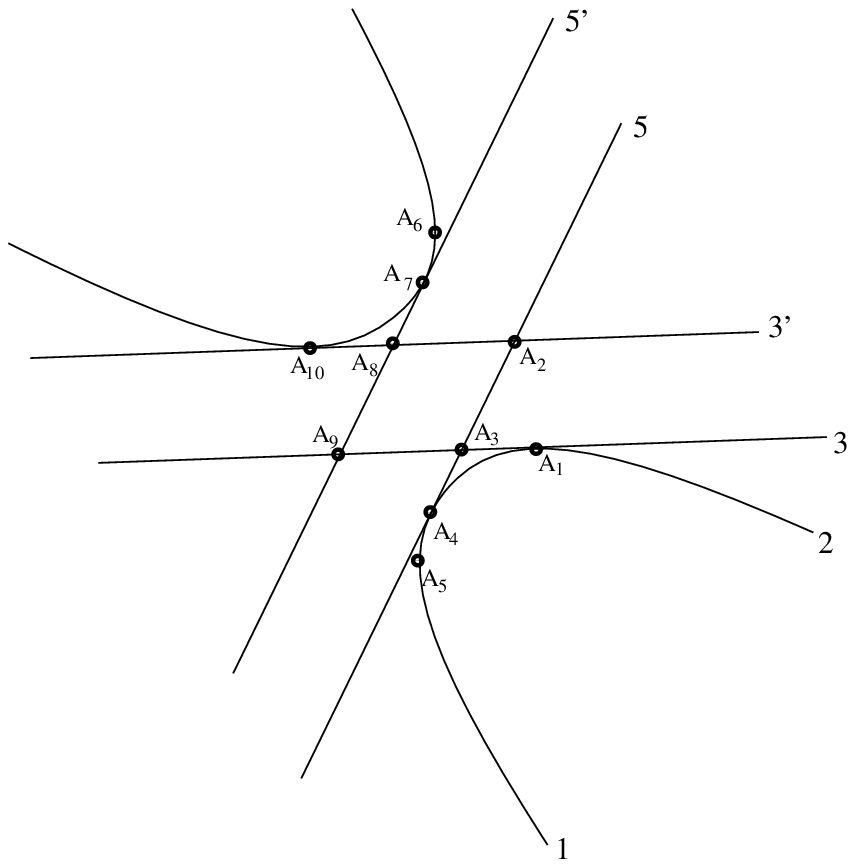}
\end{center}
\end{minipage}
\caption{}
\end{figure}

Let $M$ be a real point, $x << M \;\; \forall x \in N, K(M) = \{1,2,3,3',5,5'\}$.\\
(b) After the fifth regeneration, $S^{(3)}_1$ in a neighbourhood
of $\omega$ is as follows: the lines $L_1$ and $L_2$ are replaced
to a hyperbola $h_{12}$, the lines $L_3$ and $L_5$ are doubled and
each one of them is tangent to $h_{12}$.  This follows from (a).
Thus in a neighbourhood of $\omega, S^{(3)}_1$ is of the form:
$S^{(3)}_1 = h_{12} \cup L_3 \cup L_{3'} \cup L_5 \cup L_{5'}$,
see Figure 17.

The lines $L_3$ and $L_{3'}$ intersect in a point $1 << x(L_3 \cap L_{3'})$; the lines $L_5$ and $L_{5'}$ intersect in a point $x (L_5 \cap L_{5'})$, \st $x(L_3 \cap L_{3'}) < x (L_5 \cap L_{5'}) ; x(L_3 \cap L_{3'}), x(L_5 \cap L_{5'}) \notin E$ for $E$ a disk on the $x$-axis.

Let $\{A_j\}^{10}_{j=1}$ be singular points of $\pi_1$ as follows:\\
$A_1, A_4, A_7, A_{10}$ are tangent points of $h_{12}$ with the lines $L_3, L_5, L_{5'}, L_{3'}$ respectively.\\
$A_2, A_3$ are intersection points of the line $L_5$ with the lines $L_{3'}, L_3$ respectively.\\
$A_8, A_9$ are intersection points of the line $L_{5'}$ with the lines $L_{3'}, L_3 $ respectively.\\
$A_5$ is a point of the type $a_1$ in $h_{12}$.\\
$A_6$ is a point of the type $a_2$ in $h_{12}$.

Let $N = \{x(A_j) = x_j \mid 1 \leq j \leq 10 \}$, \st $N \subset E - \partial E, N \subset E_\R$.  Let $M$ be a real point on the $x$-axis, \st $x_j << M \; \; \forall x_j \in N, 1 \leq j \leq 10$.

We claim that there exists a $g$-base $\{\ell (\g_j)\}^{10}_{j=1}$ of $\pi_1(E - N,u)$, \st the first five paths $\g_j$ are below the real line and the last five paths $\g_j, 6 \leq j \leq 10$, have a part below the real line and a part above it, see Figure 18.

Moreover, we claim that $\vp_M$, the braid monodromy w.r.t. $E \times D$, is $\vp_M = F_1 (F_1)^{\rho^{-1}}$ where $F_1$ and $\rho$ are the ones given in the beginning of the theorem, (b).

\begin{figure}[h]
\epsfxsize=12cm 
\epsfysize=2.5cm 
\begin{minipage}{\textwidth}
\begin{center}
\epsfbox {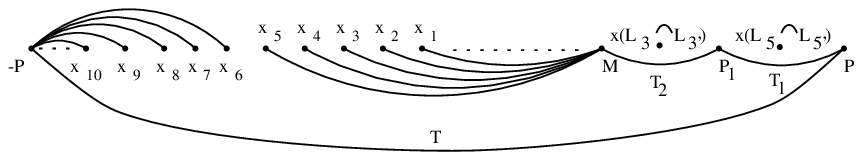}
\end{center}
\end{minipage}
\caption{}
\end{figure}

We will prove the two of them.

Let us choose a $g$-base of $\pi_1(E - N,u)$ as follows:
let us take the first five elements of the $g$-base the loops constructed from the standard bush from $M$ to $\{x_j\}^5_{j=1}$ below the real line, namely: $\ell(\g_1), \cdots , \ell(\g_5)$.
For the last five elements of the $g$-base we take $\ell(\g_6) , \cdots , \ell (\g_{10})$ where $\g_j$ are constructed as follows:  let $-P  ,  P_1, P$ be points on the real line which satisfy $M << x (L_3 \cap L_{3'}) < P_1 < x (L_5 \cap L_{5'}) < P  ,\ - P < < x_{10}$.

Let $T$ be a big semi-circle below the real line from $-P$ to $P$\ ; \ $T_1$ be a semi-circle below the real line from $P$ to $P_1$\ ; \ $T_2$ be a semi-circle below the real line from $P_1$ to $M$.

Let $\tg_j$ be paths above the real line from $x_j$ to $-P, \ 6 \le j \le 10$.  Let $\g_j$ be the paths of the form: $\g_j = \tg_j T T_1 T_2, \ 6 \le j \le 10$.

By Figure 17,
$h_{12} \cap K = \{1,2\}\ , \ L_3 \cap K = 3\ , \ L_{3'} \cap K = 3'\ ,\ L_5 \cap K = 5\ ,\ L_{5'} \cap K = 5'.$  So $K = \{1,2,3,3',5,5'\}$, \st $1< 2 < 3 < 3' < 5 < 5'$.  We choose a diffeomorphism $\beta_M$ which satisfies:
$\beta_M(1) = 1 \ , \ \beta_M(2) = 2\ , \  \beta_M(3) = 3 \ , \ \beta_M(3') = 4\ , \  \beta_M(5) = 5 \ , \ \beta_M(5') = 6.$

We have the following table for the first five points:
\begin{center}
\begin{tabular}{cccc}
$j$ & $\lambda_{x_j}$ & $\epsilon_{x_j}$ & $\delta_{x_j}$ \\
\hline
1 & $<2,3>$ & 4 & $\Delta^2<2,3>$\\
2 & $<4,5>$ & 2 & $\Delta<4,5>$\\
3 & $<3,4>$ & 2 & $\Delta<3,4>$\\
4 & $<2,3>$ & 4 & $\Delta^2<2,3>$\\
5 & $<1,2>$ & 1 & $\Delta^{\frac{1}{2}}_{I_2\R}<1>$
\end{tabular}
\end{center}

For computations we use the formulas from [Am2, Theorems 1.41, 1.44].

$(\xi_{x'_1}) \Psi_{\g'_1} = <2,3> \beta^{-1}_M = z_{23}$\\
$<2,3>$ $\vcenter{\hbox{\epsfbox{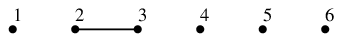}}}$  \ $\beta^{-1}_M$
$\vcenter{\hbox{\epsfbox{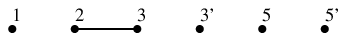}}}$ \
$\vp_M(\ell(\g_1)) = Z^4_{23}$ \\

\medskip

\noindent
$(\xi_{x'_2}) \Psi_{\g'_2} = <4,5> \Delta^2<2,3>\beta^{-1}_M = z_{3'5}$

\medskip
\noindent
$<4,5> \Delta^2<2,3>$
$\vcenter{\hbox{\epsfbox{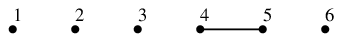}}}$  \ \
$\beta^{-1}_M$
$\vcenter{\hbox{\epsfbox{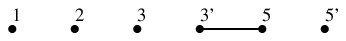}}}$  \ \
$\vp_M(\ell(\g_2)) = Z^2_{3'5}$ \\

\medskip

\noindent
$(\xi_{x'_3}) \Psi_{\g'_3}  = <3,4> \Delta<4,5>\Delta^2<2,3>\beta^{-1}_M =
\underline{z}_{35}^{Z^2_{23}}$\\
$<3,4>$ $\vcenter{\hbox{\epsfbox{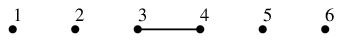}}}$  \ \
$\Delta<4,5>$
$\vcenter{\hbox{\epsfbox{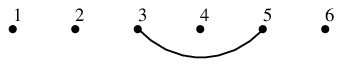}}}$  \ \
$\Delta^2<2,3>$
$\vcenter{\hbox{\epsfbox{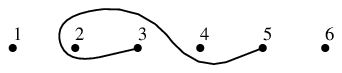}}}$  \ \
\ \  $\beta^{-1}_M$
$\vcenter{\hbox{\epsfbox{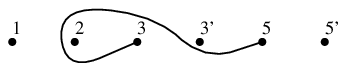}}}$  \hspace{1cm}
$\vp_M(\ell(\g_3)) = (\uZ_{35}^2)^{Z^2_{23}}$\\

\medskip

\noindent
$(\xi_{x'_4}) \Psi_{\g'_4} = <2,3> \Delta<3,4>\Delta<4,5> \Delta^2<2,3>\beta^{-1}_M = \underline{z}_{25}^{Z^2_{23}}$\\
$<2,3>$ $\vcenter{\hbox{\epsfbox{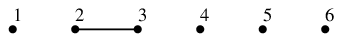}}}$  \ \
$\Delta<3,4>$ $\vcenter{\hbox{\epsfbox{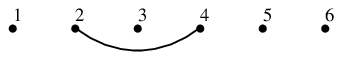}}}$  \ \
$\Delta<4,5>$ $\vcenter{\hbox{\epsfbox{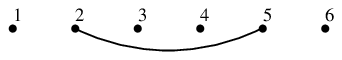}}}$  \\
$\Delta^2<2,3>$ $\vcenter{\hbox{\epsfbox{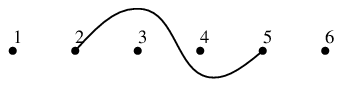}}}$  \ \
$\beta^{-1}_M$
 $\vcenter{\hbox{\epsfbox{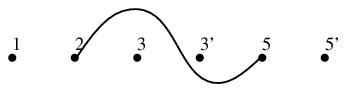}}}$
$\vp_M(\ell(\g_4)) = (\uZ^4_{25})^{Z^2_{23}}$\\

\medskip

\noindent
$(\xi_{x'_5}) \Psi_{\g'_5} = <1,2>\Delta^2<2,3>\Delta<3,4> \Delta<4,5>
\Delta^2<2,3> \beta^{-1}_M = z_{12}^{\uZ_{25}^2 Z^2_{23}}$\\
$<1,2>$  $\vcenter{\hbox{\epsfbox{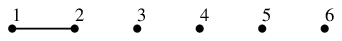}}}$  \ \   $\Delta^2<2,3>$
$\vcenter{\hbox{\epsfbox{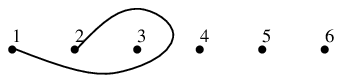}}}$  \ \
$\Delta<3,4>$
$\vcenter{\hbox{\epsfbox{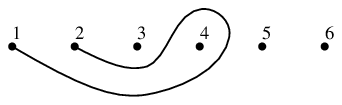}}}$\ \
$\Delta<4,5>$
$\vcenter{\hbox{\epsfbox{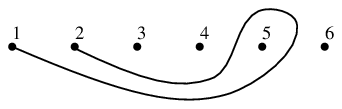}}}$ \ \
$\Delta^2<2,3>$
$\vcenter{\hbox{\epsfbox{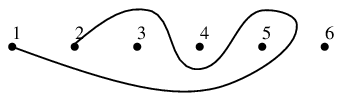}}}$ \ \
$\beta^{-1}_M$
$\vcenter{\hbox{\epsfbox{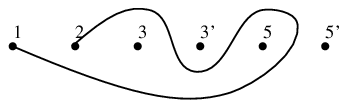}}}$ \\
$\vp_M(\ell(\g_5)) = (Z_{12})^{\uZ^2_{25}Z^2_{23}}$

The sequence of braids that we obtain is:\\
$Z^4_{23} \ , \ Z^2_{3'5} \ , \ (\uZ^2_{35})^{Z^2_{23}} \ , \ (\uZ^4_{25})^{Z^2_{23}} \ , \
(Z_{12})^{\uZ^2_{25}Z^2_{23}}$ which is the  sequence that is given by the factorized expression $F_1$.  So as a factorized expression, $F_1 = \prodl^5_{j=1} \vp_M(\ell(\g_j))$.

Now we want to compute $\beta^\vee_M (\vp_M(\ell(\g_j)))$  for $6 \le j \le 10$.
First, we have to compute $\beta^\vee_{-P} (\vp_{-P}(\ell(\tg_j)))$ for $6 \le j \le 10$.  Then we have to compute $\beta^\vee_{P} (\vp_{P}(\ell(\tg_j T)))$ for $6 \le j \le 10$.  Finally, we will compute $\beta^\vee_M (\vp_M(\ell(\g_j))) = \beta^\vee_M (\vp_M(\ell(\tg_j TT_1T_2)))$ for $6 \le j \le 10$.

In the computations of $\beta^\vee_{-P} (\vp_{-P}(\ell(\tg_j)))$, 
we apply $\delta_{x_j}$ in a reversed order (see [Am2, Remark 1.46]).  
The point $A_6$ is of the type $a_1$ w.r.t. the point $-P$.

We get the following table for the last five points:
\begin{center}
\begin{tabular}{cccc}
$j$ & $\lambda_{x_j}$ & $\epsilon_{x_j}$ & $\delta_{x_j}$ \\ \hline
10 & $<4,5>$ & 4 & $\Delta^2<4,5>$\\
9 & $<2,3>$ & 2 & $\Delta<2,3>$\\
8 & $<3,4>$ & 2 & $\Delta<3,4>$\\
7 & $<4,5>$ & 4 & $\Delta^2<4,5>$\\
6 & $<5,6>$ & 1 & $\Delta^{\frac{1}{2}}_{I_2\R}<5>$ \\
\end{tabular}
\end{center}

\medskip

Now we are computing $\beta^\vee_{-P} (\vp_{-P}(\ell(\tg_j )))$ for $j = 10, \cdots , 6$.\\

\noindent
L.V.C. $(\tg_{10}){\beta_{-P}} = z_{45}$

\medskip
\noindent
$<4,5>$ $\vcenter{\hbox{\epsfbox{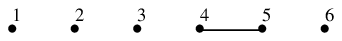}}}$ \hspace{1cm}
$\beta^\vee_{-P}(\vp_{-P}(\ell(\tg_{10})))= Z^4_{45}$\\

\medskip

\noindent
L.V.C. $(\tg_{9}){\beta_{-P}} = <2,3> \Delta^2<4,5> = z_{23}$

\medskip
\noindent
$<2,3>\Delta^2<4,5>$ $\vcenter{\hbox{\epsfbox{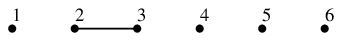}}}$ \hspace{.5cm}
$\beta^\vee_{-P}(\vp_{-P}(\ell(\tg_{9}))) = Z^2_{23}$\\

\medskip

\noindent
L.V.C. $(\tg_{8}){\beta_{-P}} = <3,4> \Delta<2,3>\Delta^2<4,5> = \bar{z}_{24}^{Z^2_{45}}$

\medskip
\noindent
$<3,4>$  $\vcenter{\hbox{\epsfbox{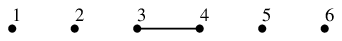}}}$ \hspace{.5cm}
$\Delta<2,3>$ $\vcenter{\hbox{\epsfbox{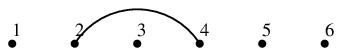}}}$ \ \
$\Delta^2<4,5>$  $\vcenter{\hbox{\epsfbox{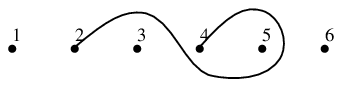}}}$ \\
\noindent
$\beta^\vee_{-P}( \vp_{-P}(\ell(\tg_{8}) ))= (\bZ^2_{24})^{Z^2_{45}}$\\

\medskip

\noindent
L.V.C. $(\tg_7)\beta_{-P} = <4,5>\Delta<3,4> \Delta<2,3>\Delta^2 <4,5> =
 \underline{z}_{25}^{Z^2_{23}}$

\medskip

\noindent
$<4,5>$ $\vcenter{\hbox{\epsfbox{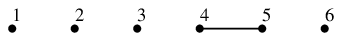}}}$ \hspace{.5cm}
$\Delta<3,4>$$\vcenter{\hbox{\epsfbox{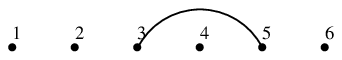}}}$ \hspace{.5cm}
$\Delta<2,3>$ $\vcenter{\hbox{\epsfbox{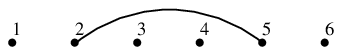}}}$\\
$\Delta^2<4,5>$
$\vcenter{\hbox{\epsfbox{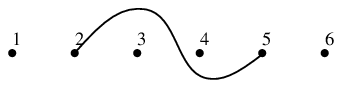}}}$\hspace{.5cm}
$\beta^\vee_{-P}(\vp_{-P}(\ell(\tg_{7}))) = (\uZ^4_{25})^{Z^2_{23}}$\\

\medskip

\noindent
L.V.C. $(\tg_6)\beta_{-P} = <5,6>\Delta^2<4,5>\Delta<3,4> \Delta<2,3>\Delta^2<4,5> = z_{56}^{\bZ^2_{25} Z^2_{45}}$

\medskip

\noindent
$<5,6>$   $\vcenter{\hbox{\epsfbox{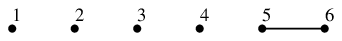}}}$ \hspace{.5cm}$\Delta^2<4,5>$
 $\vcenter{\hbox{\epsfbox{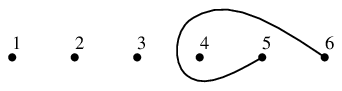}}}$ \hspace{.5cm}
$\Delta<3,4>$  $\vcenter{\hbox{\epsfbox{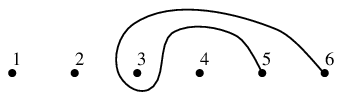}}}$ \\
$\Delta<2,3>$
$\vcenter{\hbox{\epsfbox{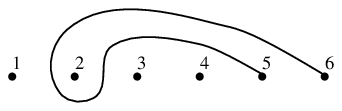}}}$
 \hspace{.5cm}$\Delta^2<4,5>$
$\vcenter{\hbox{\epsfbox{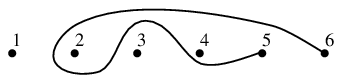}}}$  \\
$\beta^\vee_{-P}(\vp_{-P}(\ell(\tg_{6}))) = (Z_{56})^{\bZ^2_{25}Z^2_{45}}$\\

Now we compute L.V.C.$(\tg_jT)\beta_{P}$ for $j = 10, \cdots ,6$.\\
L.V.C.$(\tg_jT)\beta_P = {\rm L.V.C.} (\tg_j)\Psi_T \beta_P = {\rm L.V.C.} (\tg_j)\beta_{-P}(\beta^{-1}_{-P} \Psi_T \beta_P), \ \ j = 10 , \cdots , 6$.  Because $S^{(3)}_1$ does not have any singularities at $\infty, \beta_{-P} \Psi_T \beta_P$ is a 180$^0$ rotation centered at $3 \frac{1}{2}$.

So we apply $180^0$ rotation on L.V.C.$(\tg_j)\beta_{-P}$ to get L.V.C.$(\tg_jT)\beta_{P}$ for $j = 10, \cdots , 6$:


\noindent
L.V.C.$ (\tg_{10}T)\beta_P = \vcenter{\hbox{\epsfbox{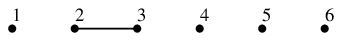}}}
  = z_{23}$ \   ;  \
L.V.C.$ (\tg_{9}T)\beta_P = \vcenter{\hbox{\epsfbox{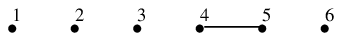}}}  = z_{45}$ \\
\noindent
L.V.C.$ (\tg_{8}T)\beta_P =  \vcenter{\hbox{\epsfbox{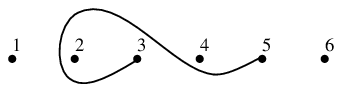}}}= \underline{z}_{35}^{Z^2_{23}}$   \  ;  \
L.V.C.$ (\tg_{7}T)\beta_P =   \vcenter{\hbox{\epsfbox{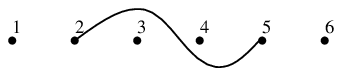}}} = \underline{z}_{25}^{Z^2_{23}}$ \\
\noindent
L.V.C.$ (\tg_{6}T)\beta_P = \vcenter{\hbox{\epsfbox{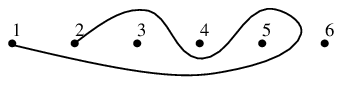}}} = z_{12}^{\uZ^2_{25}Z^2_{23}}$ .

Moreover  $\beta^\vee_P(\vp_P(\ell(\tg_jT))) = \Delta < {\rm L.V.C.} (\tg_jT)
\beta_P >^{\epsilon_{x_j}}$  for $j = 10, \cdots , 6$.\\
Thus:
$\beta^\vee_P(\vp_P(\ell(\tg_{10}T))) = Z^4_{23}$ \ ; \ $\beta^\vee_P(\vp_P(\ell(\tg_{9}T))) = Z^2_{45}$ \ ; \ $\beta^\vee_P(\vp_P(\ell(\tg_{8}T))) = (\uZ^2_{35})^{Z^2_{23}}$ \ ; \ $\beta^\vee_P(\vp_P(\ell(\tg_{7}T))) = (\uZ^4_{25})^{Z^2_{23}}$ \ ; \ $\beta^\vee_P(\vp_P(\ell(\tg_{6}T))) = Z_{12}^{\uZ^2_{25} Z^2_{23}}$.

By Figure 17, $L$-pair $(x(L_3 \cap L_{3'})) = (3,3')$ and $L$-pair $(x(L_5 \cap L_{5'})) = (5,5')$.  These two points are of the type c.  By [Am2, 
Theorem 1.40]:  $\beta^{-1}_P \Psi_{T_1T_2} \beta_M = (\Delta<3,3'> \Delta<5,5'>)^{-1}$.  Let us denote $\Delta<3,3'> \Delta<5,5'>$ by $\rho$.

Now, for $j = 10, \cdots , 6$ we have the following:\\
$\vp_M(\ell(\tg_jTT_1T_2)) = \Delta^{\epsilon_{x_j}}< $ L.V.C.$ (\tg_jTT_1T_2)> = \Delta^{\epsilon_{x_j}}< $ L.V.C.$ (\tg_jT)\Psi_{T_1T_2}> = \\
\Delta^{\epsilon_{x_j}}< $L.V.C.$ (\tg_jT)\beta_P (\beta^{-1}_P \Psi_{T_1T_2}\beta_M)\beta^{-1}_M> = \Delta^{\epsilon_{x_j}}< $ L.V.C.$ (\tg_jT)\beta_P \rho^{-1}\beta^{-1}_M> =\\
(\beta^{-1}_M)^\vee(\rho^{-1})^\vee \Delta^{\epsilon_{x_j}}< {\rm L.V.C.} (\tg_j T)
\beta_P> = (\beta^{-1}_M)^\vee(\rho^{-1})^\vee (\beta^\vee_P
(\vp_P(\ell(\tg_j T))))$.

So by comparing the beginning and the end of this long equality we have:
$\vp_M(\ell(\g_j)) = (\beta^{-1}_M)^\vee \beta^\vee_P ((\vp_P(\ell(\tg_j T)))^{\rho-1})$ for $j = 10, \cdots , 6$ and for $\rho=\Delta <3,3'> \Delta<5,5'>.$ Therefore:\\
$\vp_M(\ell(\g_{10})) = (Z^4_{23})^{\rho^{-1}}$ \\
$\vp_M(\ell(\g_{9})) = (Z^2_{3'5})^{\rho^{-1}}$ \\
$\vp_M(\ell(\g_{8})) = ((\uZ^2_{35})^ {Z^2_{23}})^{\rho^{-1}}$ \\
$\vp_M(\ell(\g_{7})) = ((\uZ^4_{25})^{Z^2_{23}})^{\rho^{-1}}$ \\
$\vp_M(\ell(\g_{6})) = ((Z_{12})^{\uZ^2_{25}Z^2_{23}})^{\rho^{-1}}$ .

From the first five points, one can easily see that as a factorized expression:\\
$F_1 = \prodl^6_{j=10} (\beta^{-1}_M)^\vee \beta^\vee_P(\vp_P(\ell(\tg_jT)))$.
From the last five points, one can easily see that as a factorized expression:
$(F_1)^{\rho^{-1}} = \prodl^6_{j=10} (\beta^{-1}_M)^\vee \beta^\vee_P((\vp_P(\ell(\tg_jT)))^{\rho^{-1}}).$
Therefore, the braid monodromy w.r.t. $E \times D$ is $\vp_M = F_1 (F_1)^{\rho^{-1}}$.
\hfill
$\Box$

\noindent
{\bf Proposition 8.3.} {\em The local braid monodromy for $S^{(3)}_1$ around $V_1$ is obtained from formula $\vp^{(6)}_1$ in Theorem 8.1 by the following replacements:\\
(i) Consider the following $( \ )^\ast$ as conjugations by the braids, induced from the following motions:

\noindent
$(\quad \quad)^{Z^{2}_{33',4}}\vcenter{\hbox{\epsfbox{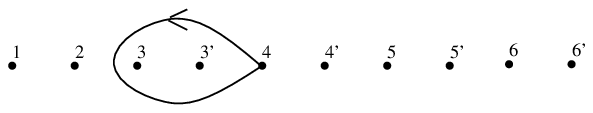}}}  $ ;

\noindent
$(\quad \quad)^{Z^2_{2,33'}}  \vcenter{\hbox{\epsfbox{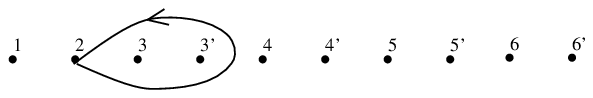}}}$ ;

\noindent
$(\quad \quad)^{\uZ^{2}_{33',6}}\vcenter{\hbox{\epsfbox{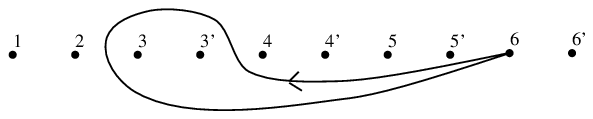}}}$ ;

\noindent
$(\quad \quad)^{\uZ^2_{33',55'}} \vcenter{\hbox{\epsfbox{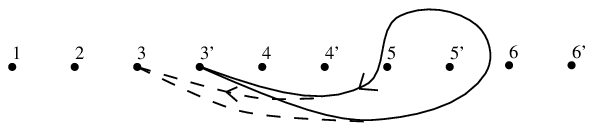}}}$ ;

\noindent
$(\quad \quad)^{Z^{2}_{55',6}}  \vcenter{\hbox{\epsfbox{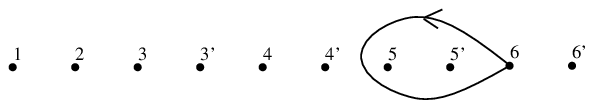}}} $

All the other conjugations do not change, since indices 3, 5 are not involved.

\noindent
(ii) $\Delta^2<1,2,3,5>$ is replaced by  $F_1(F_1)^{\rho^{-1}}$.\\
(iii) Each of the degree 4 factors in $\vp^{(6)}_1$ that involves index 3 or 5 is replaced by 3 cubes as in the third regeneration rule.\\
(iv) Each of the degree 2 factors in $\vp^{(6)}_1$ that involves index 3 or 5 is replaced by 2 degree 2 factors, where : 3 and 3' are replacing 3; 5 and 5' are replacing 5.

We call the formula that was obtained: (3)$_1$.}

\noindent
\underline{Proof}:  In the regeneration $S^{(6)}_1 \leadsto S^{(3)}_1, \hat{L}_3$ and $\hat{L}_5$ are replaced by 2 conics, each of which intersect the typical fiber twice in 2 points that are close to $\hat{L}_3 \cap \C_{V_1}$ and $\hat{L}_5 \cap \C_{V_1}$ respectively, namely 3,3',5,5'.\\
(i) By the second regeneration rule: $Z^2_{34}, Z^2_{23}, \uZ^2_{36}, \uZ^2_{35}, Z^2_{56}$
should be replaced by the ones appearing above.\\
(ii) By Theorem 8.2.\\
(iii) By the third regeneration rule.\\
(iv) By the second regeneration rule.
\hfill
$\Box$

Now we compute $\vp^{(2)}_1$.  Recall that $\vp^{(2)}_1 = \vp^{(1)}_1 = \vp^{(0)}_1$, so in fact we compute $\vp^{(0)}_1$.

\noindent
{\bf Theorem 8.4.}  {\em In the notation of Theorem 8.2:  let $T^{(0)}$ be the curve
 obtained from $T^{(3)}$ in the regeneration $Z^{(3)} \leadsto Z^{(0)}$. \\
(i) Then the local braid monodromy of $T^{(0)}$ is $\hat{F}_1 (\hat{F}_2)$, where:\\
$\hat{F}_1 = \uZ^3_{22',3} \cdot Z^2_{3'5}\cdot(\uZ^2_{35})^{Z^2_{22',3}}
 \cdot(\uZ^3_{22',5})^{Z^2_{22',3}} \; \; \cdot Z_{12'} \cdot Z_{1'2} ,$ \\
$\hat{F}_2 = (\uZ^3_{22',3})^{\rho^{-1}} \cdot(Z^2_{3'5})^{\rho^{-1}}
\cdot\left((\uZ^2_{35})^{Z^2_{22',3}}\right)^{\rho^{-1}}\cdot\left((\uZ^3_{22',5})^{Z^2_{22',3}}\right)^{\rho^{-1}} \; \; \cdot Z_{12'}^{\rho^{-1}} 
\cdot Z_{1'2}^{\rho^{-1}} $. \\
(ii) The singularities of the $x$-projection of \ $T^{(0)}$ are those arising from $T^{(3)}$ by the regeneration rules, namely: 4 nodes that exist in $T^{(3)}, 3 \times 4$ cusps arising from 4 tangency points in $T^{(3)}, 2 \times 2$ branch points for 2 branch points of $T^{(3)}$.\\
(iii) The braid monodromy of $T^{(0)}$ is $\hat{F}_1 (\hat{F}_2)$. \\
(iv) $\hat{F}_1 (\hat{F}_2) = \Delta^{-2}_8 Z^{-2}_{33'} Z^{-2}_{55'} Z^{-2}_{11'} Z^{-2}_{22'}$.}

\noindent
\underline{Proof}: We use the notation of Theorem 8.2 (formulation and proof).\\
(i) In the  regeneration process,  $Z^{(3)} \leadsto Z^{(0)}, W_1 \cup W_2$ becomes irreducible with equations $Z Y = \epsilon $ and $XT = \delta$ where the branch curve in $\C \P^2$ of its projection from $U_3$ to $\C \P^2$ is a regeneration of $T^{(3)}$, \st the conic is doubled in the generic point.  We apply the regeneration rules (I, III) on the local braid monodromy of $T^{(3)}$  on each factor that involves indices 1 or 2, to get
 $\hat{F}_1(\hat{F}_2)$.\\
(ii)  Let $V$ be an affine surface in $\C^4$, defined by: $YZ = \epsilon $ and $XT = \delta$.  Let $\bar{V}$ be its projective closure in $\C \P^4$, defined by: $YZ = \epsilon W^2$ and $XT = \delta W^2$.  It is easy to check that the only singularities of $\bar{V}$ are 4 non-degenerate quadric singularities at $\infty$, namely the intersection of $\hat{L}_i$ with $\{W=0\}, i = 1,2,3,5$.  Let $S_V$ be the branch  curve of the $x;y$ projection of $V$.  Let $\bar{S}_V$ be its projective closure.  Each of the singularities of $\bar{V}$ corresponds to a singularity of $\bar{S}_V$, namely: the intersection of $L_i$ with the line in infinity, $i = 1,2,3,5$.  Thus $\bar{S}_V$ has 4 nodes at $\infty$.  $T^{(0)}$ can be identified with a part of $S_V$.  In particular, each singularity of the $x$-projection of $T^{(0)}$ corresponds to a singularity of the $x$-projection of $S_V$.  We claim that there are no more singularities
of $\bar{S}_V$ besides those mentioned above: the 4 nodes at $\infty$ and those that come from $T^{(0)}$.  We prove it by comparing degrees of the factorized expressions given by the braid monodromy of $S_V$.

The degree of $S_V$ is 8, degree $\Delta^2_8 = 8 \cdot 7 = 56$. The degree of the product in $\Delta^2_8$ of the factors that come from the nodes at $\infty
$ is $4 \cdot 2 = 8$.  The singularities of $T^{(0)}$ that come from the regenerating singularity of $T^{(3)}$ contributes to the factorization of $\Delta^2_8$  factors whose added degrees are: $4 \cdot 2 + 12 \cdot 3 + 4 \cdot 1 = 48$.  Thus, we have factors together whose added degrees are 48 + 8 = 56 and there cannot be any additional singularities of the $x$-projection of $T^{(0)}$ which prove part (ii).  Moreover, $S_V$ has no additional singularities, and $\bar{S}_V$ has only the 4 nodes at $\infty$.\\
(iii)  Since $T^{(0)}$ has no additional singularities than those coming from regenerating $T^{(0)}$ , these contribute all the factors of $\hat{F}_1(\hat{F}_2)$.\\
(iv) To prove this point, we have to analyse the contribution of the nodes at $\infty$ to the braid monodromy of $\bar{S}_V$.

We move the line in infinity, \st  $\bar{S}_V$ is transversal to it.  Denote the new affine part of  $\bar{S}_V$ by $S'_V$.  Then the nodes of
 $\bar{S}_V$ are (transversal) intersection points of $L_i (i = 1,2,3,5)$ with a line at $\infty$. We want to compute the braid monodromy of $S'_V$.
In the regeneration process that created $S'_V$ these lines were doubled so that the L.V.C. corresponding to the nodes of $S'_V$
close to $\infty$ are $z_{11'}, z_{22'}, z_{33'}, z_{55'}$ and the contribution of those nodes to the braid monodromy of $S'_V$ is
$Z^2_{11'}, Z^2_{22'}, Z^2_{33'}, Z^2_{55'}$.  We get the following formula for the braid monodromy of
$S'_V: \Delta^2_8 = \hat{F}_1 \hat{F}_2 Z^2_{11'}  Z^2_{22'} Z^2_{33'} Z^2_{55'}$.
\hfill
$\Box$

The local braid monodromy of $S^{(0)}_1$ in $U_3$ is $\hat{F}_1 (\hat{F}_1)^{\rho^{-1}}$.  Applying regeneration rules on $F_1(F_1)^{\rho^{-1}}$, we obtain $\hat{F}_1 (\hat{F}_1)^{\rho^{-1}}$.

{\bf Corollary 8.5.} {\em The paths corresponding to $\hF_1(\hF_1)^{\rho^{-1}}$ (without their conjugation): \\
(a) \underline{The paths corresponding to the factors in $\hF_1$}:


\medskip
\noindent$\vcenter{\hbox{\epsfbox{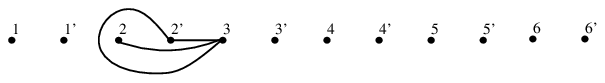}}}$ \hspace{.5cm}corresponding to $\uZ^3_{22',3}$

\medskip

\noindent$\vcenter{\hbox{\epsfbox{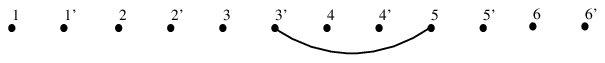}}}$ \hspace{.5cm} corresponding to $Z^2_{3'5}$

\medskip
\noindent$\vcenter{\hbox{\epsfbox{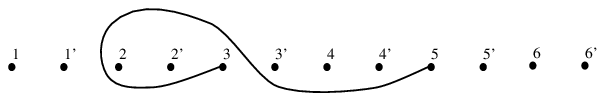}}}$ \hspace{.5cm} corresponding to $(\uZ^2_{35})^{Z^2_{22',3}}$

\medskip
\noindent$\vcenter{\hbox{\epsfbox{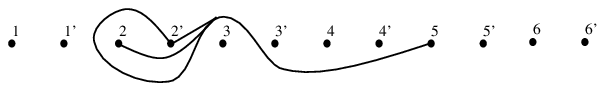}}}$ \hspace{.5cm} corresponding to $(\uZ^3_{22',5})^{Z^2_{22',3}}$

\medskip
\noindent$\vcenter{\hbox{\epsfbox{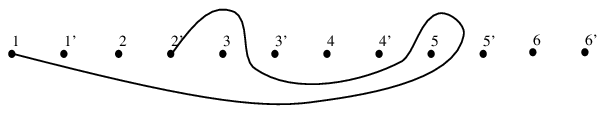}}}$ \hspace{.5cm} corresponding to 
$Z_{12'}$

\medskip
\noindent$\vcenter{\hbox{\epsfbox{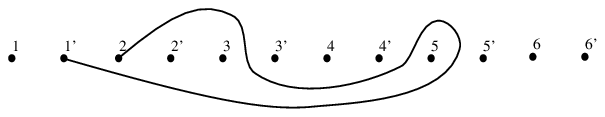}}}$ \hspace{.5cm} corresponding to  $Z_{1'2}$

\medskip

\noindent
(b) \underline{The paths corresponding to the factors in $(\hF_1)^{\rho^{-1}}$ for
$\rho^{-1} = Z^{-1}_{55'} Z^{-1}_{33'}$ }\ : \\

\noindent$\vcenter{\hbox{\epsfbox{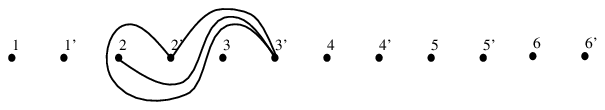}}}$ \hspace{.5cm} corresponding to $(\uZ^3_{22',3})^{\rho^{-1}}$

\medskip
\noindent$\vcenter{\hbox{\epsfbox{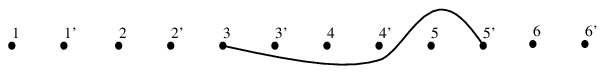}}}$ \hspace{.5cm} corresponding to $(Z^2_{3'5})^{\rho^{-1}}$

\medskip
\noindent$\vcenter{\hbox{\epsfbox{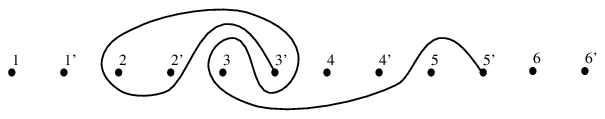}}}$ \hspace{.5cm} corresponding to $\left((\uZ^2_{35})^{Z^2_{22',3}}\right)^{\rho^{-1}}$

\medskip
\noindent  $\vcenter{\hbox{\epsfbox{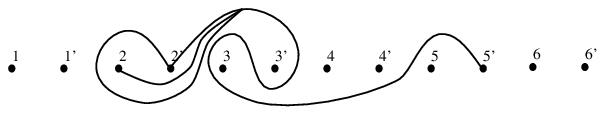}}}$ \hspace{.5cm}  corresponding to $\left((\uZ^3_{22',5})^{Z^2_{22',3}}\right)^{\rho^{-1}}$

\medskip
\noindent $\vcenter{\hbox{\epsfbox{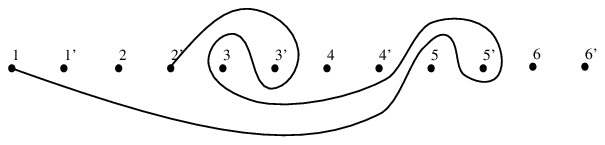}}}$ \hspace{.5cm} corresponding to $Z_{12'}^{\rho^{-1}}$

\medskip
\noindent $\vcenter{\hbox{\epsfbox{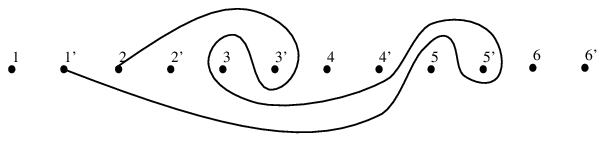}}}$ \hspace{.5cm}
 corresponding to  $Z_{1'2}^{\rho^{-1}}$}.

{\bf Theorem 8.6.} {\em  The local braid monodromy of $S^{(0)}_1$ around $V_1$, denoted by $H_{V_1}$, equals:\\
\noindent
$H_{V_1} = (\uZ^2_{4'i})_{i = 22',55',6,6'} \cdot Z^3_{33',4} \cdot
(\uZ^2_{4i})_{i = 22',55',6,6'}^{Z^2_{33',4}} \cdot
\uZ^3_{11',4} \cdot
(Z_{44'})^{Z^2_{33',4}\uZ^2_{11',4}} \cdot
Z^3_{55',6} \cdot (\uZ^2_{i6})^{Z^2_{i,22'}Z^2_{33',4}}_{i=11',33'} \cdot
\uZ^3_{22',6} \cdot
(\uZ^2_{i6'})^{\uZ^2_{i6}\uZ^2_{i,55'}Z^2_{i,22'}Z^2_{33',4}}_{i = 11',33'} \cdot
(Z_{66'})^{Z^2_{55',6}\uZ^2_{22',6}} \cdot
\left(\hF_1(\hF_1)^{\rho^{-1}}\right)^{Z^2_{55',6}Z^2_{33',4}}$, where the paths corresponding to these braids are (the paths corresponding to $
\hF_1(\hF_1)^{\rho^{-1}}$ are above):  \\

\noindent
(1)\hspace{.5cm}$\vcenter{\hbox{\epsfbox{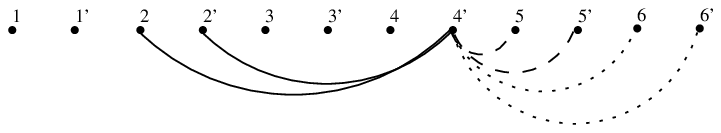}}}$ \hspace{.5cm} corresponding to $(\uZ^2_{4'i})_{i = 22',55',6,6'}$

\medskip
\noindent
(2)\hspace{.5cm}$\vcenter{\hbox{\epsfbox{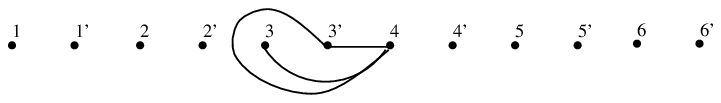}}}$ \hspace{.5cm}
 corresponding to $Z^3_{33',4}$

\medskip
\noindent
(3)\hspace{.5cm}$\vcenter{\hbox{\epsfbox{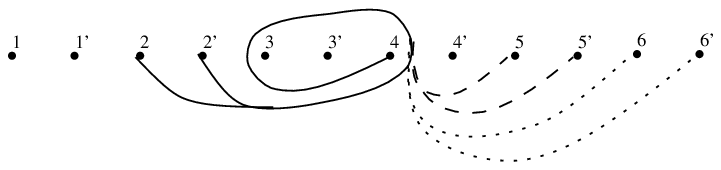}}}$ \hspace{.5cm}
corresponding to $(\uZ^2_{4i})_{i = 22',55',6,6'}^{Z^2_{33',4}}$

\medskip
\noindent(4)\hspace{.5cm} $\vcenter{\hbox{\epsfbox{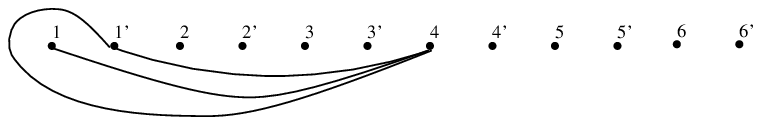}}}$ 
\hspace{.5cm} corresponding to $\uZ^3_{11',4}$

\medskip
\noindent(5)\hspace{.5cm}$\vcenter{\hbox{\epsfbox{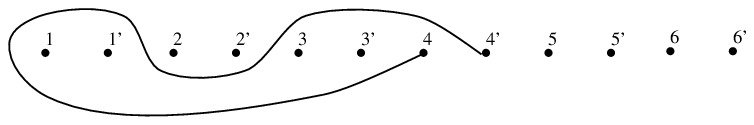}}}$ \hspace{.5cm}corresponding to $(Z_{44'})^{Z^2_{33',4}\uZ^2_{11',4}}$

\medskip
\noindent(6) \hspace{.5cm}$\vcenter{\hbox{\epsfbox{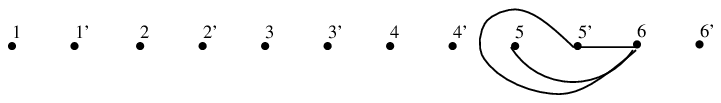}}}$ \hspace{.5cm}corresponding to $Z^3_{55',6}$

\medskip
\noindent(7)  \hspace{.5cm}$\vcenter{\hbox{\epsfbox{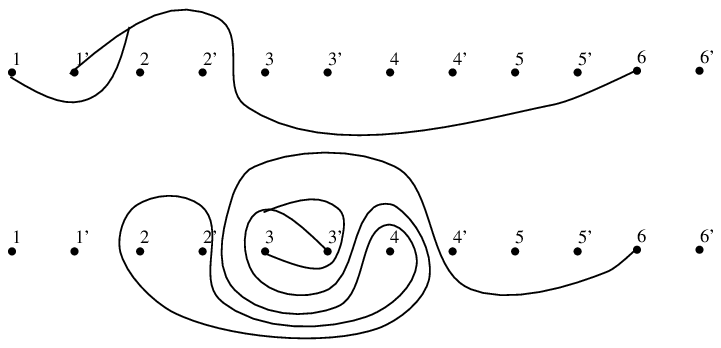}}}$ \hspace{.5cm}corresponding to $(\uZ^2_{i6})^{Z^2_{i,22'}Z^2_{33',4}}_{i=11',33'}$

\medskip

\medskip
\noindent(8) \hspace{.5cm}$\vcenter{\hbox{\epsfbox{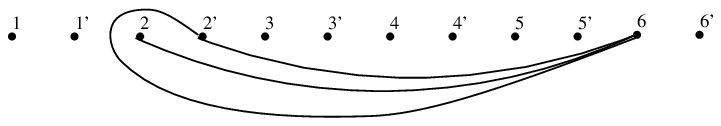}}}$ \hspace{.5cm}corresponding to $\uZ^3_{22',6}$

\medskip
\noindent(9) \hspace{.5cm}$\vcenter{\hbox{\epsfbox{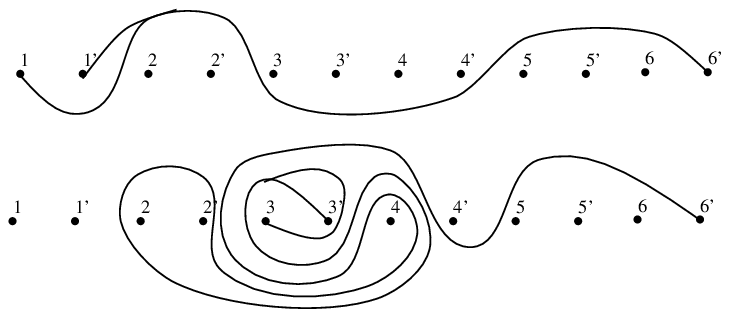}}}$ \hspace{.5cm}corresponding to  $(\uZ^2_{i6'})^{\uZ^2_{i6}\uZ^2_{i,55'}Z^2_{i,22'}Z^2_{33',4}}_{i = 11',33'}$

\medskip

\noindent(10) \hspace{.5cm}$\vcenter{\hbox{\epsfbox{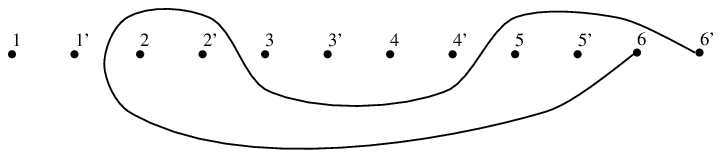}}}$ \hspace{.5cm}corresponding to  $(Z_{66'})^{Z^2_{55',6}\uZ^2_{22',6}}$}

\noindent
\underline{Proof:} In the last step of regeneration, $L_1$ and $L_2$
are replaced by conics and in the typical fiber 1 is replaced by 1 and 1', 2
is replaced by 2 and 2'.  To obtain a formula for a local braid monodromy of $S^{(0)}_1$
around $V_1$ we have to take formula (3)$_1$ and make the following changes: \\
(i) All the $(\;)^\ast$ are conjugations by the braids induced from motions
 as mentioned in the beginning of the theorem.  These conjugations affect the above paths.\\
(ii) Replace $F_1(F_1)^{\rho^{-1}}$ by $\hF_1 (\hF_1)^{\rho^{-1}}$. \\
(iii) Each of the degree 4 factors in (3)$_1$, that involves index 1 or 2, is replaced by 3
 cubes according to the third regeneration rule.\\
(iv) Each of the degree 2 factors outside of $F_1(F_1)^{\rho^{-1}}$ that involves index 1 or
 2 is replaced by 2 degree 2 factors as in the second regeneration rule.
\hfill
$\Box$

\noindent
{\bf Invariance Property 8.7.}  {\em $H_{V_1}$ is invariant under
$(Z_{44'} Z_{66'})^p (Z_{33'}Z_{55'})^q (Z_{11'} Z_{22'})^r \quad \forall p,q,r\in \Z$.}

\noindent
\underline{Proof:}  Recall Invariance Rules I, II, III, Invariance Remarks (i)-(v), the Chakiri 
Lemma, Conjugation Property   and Complex Conjugation in Section \ref{ssec:241}.\\
\underline{Case 1}: $p = q = r$.\\
$\Delta^2_{54} = \prodl^9_{i=1} C_i H_{V_i} \prod b_i$.  For each
$i$, we can choose a forgetful homomorphism, which neglects all
the indices which are
 not the indices in $H_{V_i}$.  For braids of the form $Z^2_{ii',j}$ and $Z_{ii'}$
we obtain $Z^2_{ii'}, Z_{ii'}$ when forgetting $j$.  For $Z^3_{ii',j}$
we can apply here Lemma 2(i) in [MoTe8] to obtain $Z_{ii'}$.  Applying this homomorphism
each time on $\Delta^2_{54}$ we obtain $\Delta^2_{12}$.  So each $H_{V_i}$ is a product
of $\Delta^2_{12}$ with $Z^m_{ii'}$, for some $m$.  Thus by Chakiri's Lemma, $H_{V_1}$
is invariant under this product with a power $-p$.  Since $\Delta^2_{12}$ is a central
 element, $H_{V_1}$ is invariant under $(Z_{44'} Z_{66'})^p (Z_{33'}Z_{55'})^p (Z_{11'}
  Z_{22'})^p$.

\noindent
\underline{Case 2}: $p = 0$.\\
Denote: $\epsilon = (Z_{33'} Z_{55'})^q (Z_{11'}Z_{22'})^r$.
We want to prove that  $H_{V_1}$ is invariant under $\epsilon $.

\noindent
\underline{Step 1:}  Product of the factors outside of $\hF_1(\hF_1)^{\rho^{-1}}$.
\\
$(Z_{66'})^{Z^{2}_{55',6}
\uZ^2_{22',6}}$ and
$(Z_{44'})^{ Z^2_{33',4}\uZ^2_{11',4}}$ commute with $\epsilon $. \\
$Z^3_{33',4}, \uZ^3_{11',4}, Z^3_{55',6} , \uZ^3_{22',6} $ are
invariant under $\epsilon $ by Invariance Rule III.

All conjugations are invariant under $\epsilon $ 
by Invariance Rule II and by Invariance 
Remark (iv).
The degree 2 factors are of the form $Z^2_{\alpha \alpha ', \beta
}$ where
 $\beta = 6,6',4,4'$ or $Z^2_{\alpha , \beta \beta '}$ where $\alpha = 4$ or 4'.
  By Invariance Rule II, they are invariant under $Z_{\alpha \alpha '}$ or $Z_{\beta
   \beta '}$ respectively, and since the other halftwists in $\epsilon $ commute with
   $Z_{\alpha \beta }, Z_{\alpha ' \beta }$ ( or $Z_{\alpha \beta }, Z_{\alpha \beta '})$
    we get that $Z^2_{\alpha \alpha ', \beta }$ and $Z^2_{\alpha , \beta \beta '}$ are
    invariant under $\epsilon $.

\noindent
\underline{Step 2}:   $\hF_1(\hF_1)^{\rho^{-1}}$.\\
In order to prove that  $\hF_1(\hF_1)^{\rho^{-1}}$ is invariant under $\epsilon $
we consider the following subcases:\\
\underline{Subcase 2.1}:  $q = 0 ; \epsilon = (Z_{11'}Z_{22'})^r$.\\
$\uZ^3_{22',3}, (\uZ^3_{22',5})^{Z^2_{22',3}}$ are invariant under $Z_{22'}$
(Invariance Rule III) and commute with $Z_{11'}$ (Invariance Remark (iv)). Thus
 they are invariant under $\epsilon $ (Invariance Remarks (i) and (v)).  $Z^2_{3'5},
 (\uZ^2_{35})^{Z^2_{22',3}}$ commute with $Z_{11'}$ and $Z_{22'}$ and thus with
 $\epsilon$.  $\alpha _1(m_1), \alpha _2(m_1)$ are invariant under
($Z_{11'}  Z_{22'})^r$ by Invariance Rule I.  Thus $\hF_1$ is invariant under $\epsilon $.
 Since $\rho^{-1} = Z^{-1}_{55'} Z^{-1}_{33'}, \rho^{-1}$ commutes with $\epsilon$ and
  thus $(\hF_1)^{\rho^{-1}}$ is also invariant under $\epsilon $.  Thus
   $\hF_1(\hF_1)^{\rho^{-1}}$ is invariant under \nolinebreak $\epsilon $.

\noindent
\underline{Subcase 2.2}:  $r=0 \ , \ q = 1 \ ; \ \epsilon = Z_{33'}Z_{55'} \ ,\
\epsilon = \rho$. \\
To prove that  $\hF_1(\hF_1)^{\rho^{-1}}$  is invariant under $\rho$, it is
necessary to prove that $\hF_1(\hF_1)^{\rho^{-1}}$  is Hurwitz equivalent $(He)$(see [Am2, Definition 1.55]) to
$(\hF_1)^{\rho}\hF_1$. Since $AB$ is Hurwitz equivalent to $BA^B$, it is enough
 to prove that  $\hF_1(\hF_1)^{\rho^{-1}}$ is Hurwitz equivalent to
 $\hF_1\left((\hF_1)^\rho\right)^{\hF_1}$.  Thus it is enough to prove that
$(\hF_1)^{\rho^{-1}}$ is Hurwitz equivalent to
  $\left((\hF_1)^{\rho}\right)^{\hF_1}$ or that
   $\left((\hF_1)^{\rho^{-1}}\right)^{\hF_1^{-1}}$ is Hurwitz equivalent to $(\hF_1)^\rho$.

By Theorem 8.4 (iv): $\hF_1 = \Delta^2_8 \rho^{-2} Z^{-2}_{11'} Z^{-2}_{22'}
(\hF_1^{-1})^{\rho^{-1}}$.  Thus $(\hF_1)^{-1} = (\hF_1)^{\rho^{-1}} (Z_{11'}
 Z_{22'})^2 \rho^2 \Delta^2_8$.  Thus:\\
$\left((\hF_1)^{\rho^{-1}}\right)^{\hF_1^{-1}}  =  \left((\hF_1)^{\rho^{-1}}\right)^{
(\hF_1)^{\rho^{-1}} (Z_{11'}Z_{22'})^2 \rho^2 \ \Delta^2_8}
\left . \begin{array}{c}
= \\ [-.5cm]{\scriptstyle as \ factorized} \\ [-.5cm]{\scriptstyle expression}
\end{array}  \right .
 \left((\hF_1)^{\rho^{-1}}
\right)^{(\hF_1)^{\rho^{-1}} (Z_{11'}Z_{22'})^2 \rho^2}\\
\left . \begin{array}{c}
\simeq \\ [-.5cm]{\scriptstyle He \ by \ Chakiri's} \\ [-.5cm]{\scriptstyle
Lemma} \end{array}  \right . $
$\left((\hF_1)^{\rho^{-1}}\right)^{
(Z_{11'}Z_{22'})^2 \rho^2} =  \left((\hF_1)^{
(Z_{11'}Z_{22'})^2} \right)^\rho \; \left . \begin{array}{c}
\simeq \\ [-.5cm]{\scriptstyle He \ by  \ Subcase \ 2.1}
\end{array}  \right .(\hF_1)^{\rho} .$

\noindent
\underline{Subcase 2.3}:  $q=2q' \ ; \ \epsilon = (Z_{33'}Z_{55'})^{2q'}
(Z_{11'}Z_{22'})^r$. \\
By Theorem 8.4 (iv): $\hF_1(\hF_1)^{\rho^{-1}} = \Delta^2_8(Z_{33'} Z_{55'})^{-2}
(Z_{11'}Z_{22'})^{-2}$.  By Chakiri's Lemma,
$\hF_1(\hF_1)^{\rho^{-1}}$  is invariant under $(\Delta^2_8 (Z_{33'}Z_{55'})^{-2}
(Z_{11'}Z_{22'})^{-2})^{-q'}$ and thus under $(Z_{33'}Z_{55'})^{2q'}
(Z_{11'}Z_{22'})^{2q'}$.  By Subcase 2.1, $\hF_1(\hF_1)^{\rho^{-1}}$  is invariant
 under  $(Z_{11'}Z_{22'})^{r-2q'}$.  By Invariance Remark (v):
  $\hF_1(\hF_1)^{\rho^{-1}}$ is invariant under $(Z_{33'}Z_{55'})^{2q'} (Z_{11'}Z_{22'})^{2q'}
(Z_{11'}Z_{22'})^{r-2q'} = \epsilon $.\\
\noindent
\underline{Subcase 2.4}:  $q=2q'+1 \ ; \  \epsilon = (Z_{33'}Z_{55'})^{2q'+1}
(Z_{11'}Z_{22'})^r$. \\
By Subcases 2.2, 2.3 and Invariance Remark (v). Thus $\hF_1(\hF_1)^{\rho^{-1}}$
is invariant under $(Z_{33'}Z_{55'})^{q} \\
(Z_{11'}Z_{22'})^r \; \forall q, r \in \Z$.  This finishes Step 2 of Case 2
and thus Case 2 is complete.\\
\noindent
\underline{Case 3}:  $p,q,r$ arbitrary; $\epsilon = (Z_{44'}Z_{66'})^{p} (Z_{33'}Z_{55'})^{q}
(Z_{11'}Z_{22'})^{r}$. \\
By Case 1: $H_{V_1}$ is invariant under
$(Z_{44'}Z_{66'})^{p} (Z_{33'}Z_{55'})^{p}
(Z_{11'}Z_{22'})^{p}$. By Case 2:$H_{V_1}$ is invariant under
$(Z_{33'}Z_{55'})^{q-p} (Z_{11'}Z_{22'})^{r-p}$. By Invariance Remark (v),
 $H_{V_1}$ is invariant \nolinebreak  under $\epsilon $.
\hfill
$\Box$
\noindent
\underline{Notation}: $Z_{ii'} = \rho_i$.\\
\noindent
{\bf Theorem 8.8.}  {\em Consider $H_{V_1}$ from Theorem 8.6.
The paths corresponding to the factors in $H_{V_1}$ ,
 are shown below considering the Invariance Property 8.7 and the Conjugation Property.}\\
\noindent
\underline{Remark}:  By abuse of notation, the simple braids are denoted by $\underline{z}_{mk}$ and  the more complicated paths are denoted by $\tilde{z}_{mk}$.\\
\underline{Proof}:  By the Invariance Property 8.7,
 $H_{V_1}$ is invariant under $\epsilon =(\rho_4 \rho_6)^p (\rho_3 \rho_5)^q (\rho_1 \rho_2)^r \quad\\
 \forall p, q, r \in \Z$.
There are two possible applications for presenting the braids:
 (a) if $m,k$ are two indices in two different parts of $\epsilon $,
  then the braid is
$\rho^j_k \rho^i_m z_{mk} \rho^{-i}_m \rho^{-j}_k, \ i, j \in \Z, \ i \neq j$.
(b) If $m,k$ are in the same part of $\epsilon $, then the braid
is $\rho^i_k \rho^i_m z_{mk} \rho^{-i}_m \rho^{-i}_k$, i.e.: when conjugating
 $z_{mk}$ by $\rho^i_m$ for some $i$, we must conjugate $z_{mk}$ also
 by $\rho^i_k$ for the same $i$.

By Conjugation Property,  the first figure of (9) is
conjugated by (10), in order to simplify (9).
\begin{center}
\begin{tabular}[hb]{|l|p{3in}|l|c|}\hline
& The paths corresponding & The braids & The exponent \\
[-.2cm] & to the braids and their  & &  (according to  \\
[-.2cm] & complex conjugates & & singularity type)\\ \hline
& & & \\ [-.5cm]
(1) & $\vcenter{\hbox{\epsfbox{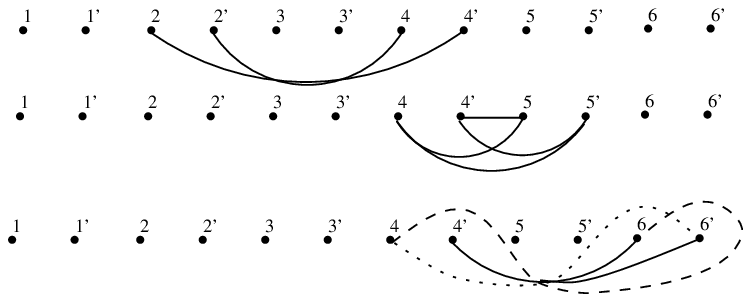}}}$ & $\begin{array}{ll}\rho^j_4 \rho^i_m \underline{z}_{m4} \rho^{-i}_m \rho^{-j}_4 \\ [-.2cm]
m = 2,5,6  
\end{array}$ & 2 \\  \hline
& & & \\ [-.5cm]
(2) & $\vcenter{\hbox{\epsfbox{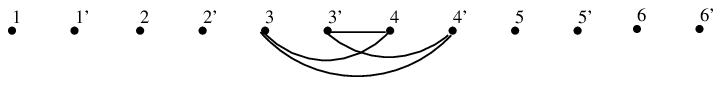}}}$ & $\rho^j_4 \rho^i_3 z_{34} \rho^{-i}_3 \rho^{-j}_4$ & 3  \\[.2cm]  \hline
(3) & $\vcenter{\hbox{\epsfbox{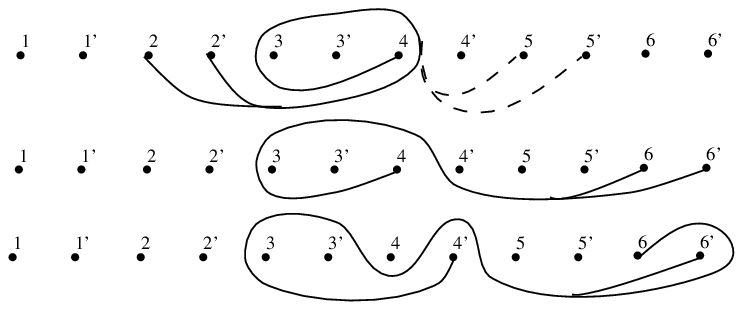}}}$ & $\begin{array}{ll}\rho^j_4 \rho^i_m\tilde{z}_{4m}\rho^{-i}_m \rho^{-j}_4  
 \\ [-.2cm] m = 2,5,6  \end{array}$ &  2  \\  \hline
&&&\\ [-.4cm]
(4) & $\vcenter{\hbox{\epsfbox{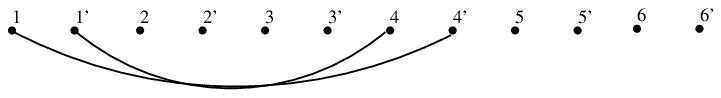}}}$ & $\rho^j_4 \rho^i_1 \underline{z}_{14} \rho^{-i}_1 \rho^{-j}_4$ & 3 \\  \hline
& & & \\ [-.4cm]
(5) & $\vcenter{\hbox{\epsfbox{1v071.ps}}}$ & $\tilde{z}_{44'}$  & 1 \\
[.3cm] \hline
& & & \\[-.4cm]
(6) &   $\vcenter{\hbox{\epsfbox{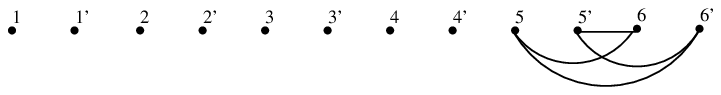}}}$ &  $\rho^j_6 \rho^i_5 z_{56} \rho^{-i}_5 \rho^{-j}_6$ & 3  \\ \hline
(7) & $\vcenter{\hbox{\epsfbox{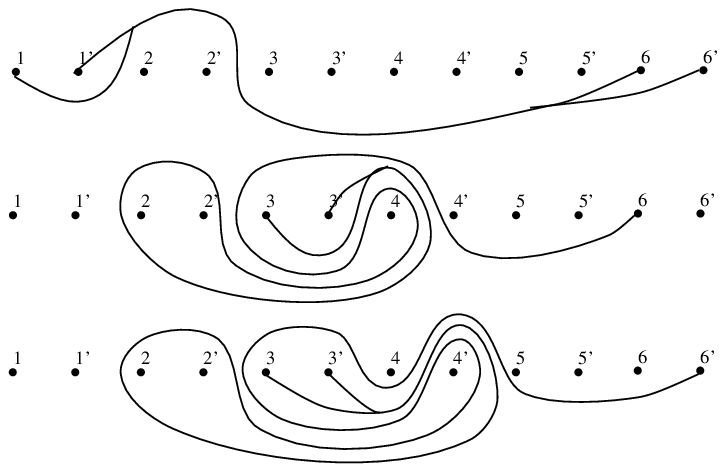}}}$ & $\begin{array}{ll}
\rho^j_6 \rho^i_m \tilde{z}_{m6}\rho^{-i}_m \rho^{-j}_6   \\[-.2cm]
m = 1,3\end{array}$ & 2 \\ \hline
& & & \\ [-.5cm](8) &   $\vcenter{\hbox{\epsfbox{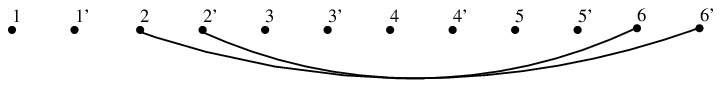}}}$ & $\rho^j_6 \rho^i_2 \underline{ z}_{26} \rho^{-i}_2 \rho^{-j}_6$ & 3 \\  \hline
& & & \\ [-.5cm](9) &  $\vcenter{\hbox{\epsfbox{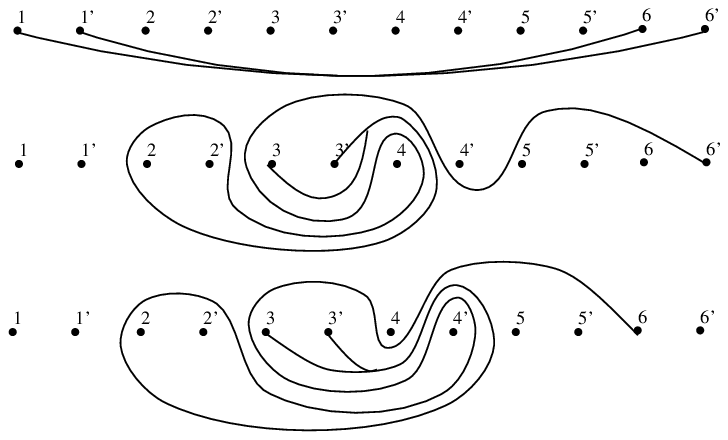}}}$ & $
\begin{array}{ll}
\rho^j_6 \rho^i_m \tilde{z}_{m6}\rho^{-i}_m \rho^{-j}_6\\ [-.2cm]
m = 1,3\end{array}$ & 2 \\ \hline
(10) &   $\vcenter{\hbox{\epsfbox{1v077.ps}}}$ &  $\tilde{z}_{66'}$ & 1 \\
   \hline
\end{tabular}
\end{center}
Consider now $\hF_1 (\hF_1)^{\rho^{-1}} (\rho^{-1} = \rho_5^{-1} \rho_3^{-1})$ in $H_{V_1}$.  The factors in $\hF_1$ and in $(\hF_1)^{\rho^{-1}}$ and their corresponding paths appear in Corollary 8.5.  Considering the Invariance Property 8.7, we obtain the following table of braids from $\hF_1$ and from $(\hF_1)
^{\rho^{-1}}$,
in which all paths are considered here and induce relations.
\begin{center}
\begin{tabular}[H]{|p{3in}|c|c|}\hline
$\vcenter{\hbox{\epsfbox{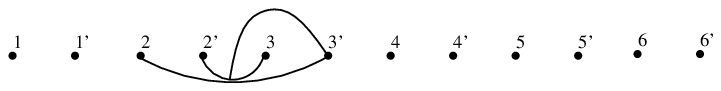}}}$  & $\rho^j_3 \rho^i_2 \underline{z}_{23} \rho^{-i}_2 \rho^{-j}_3$ & 3\\
   \hline
$\vcenter{\hbox{\epsfbox{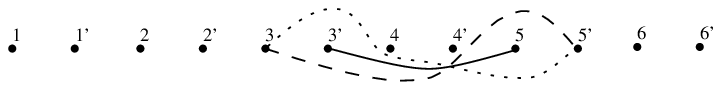}}}$  & $\rho^i_5 \rho^i_3  \underline{z}_{3'5} \rho^{-i}_3 \rho^{-i}_5$ & 2\\
    \hline
$\vcenter{\hbox{\epsfbox{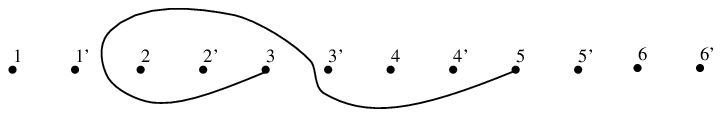}}}$ & $\rho^i_5 \rho^i_3 \tilde{z}_{35} \rho^{-i}_3 \rho^{-i}_5$ & 2\\
$\vcenter{\hbox{\epsfbox{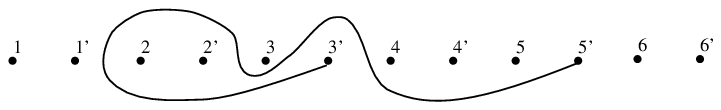}}}$ & & \\
$\vcenter{\hbox{\epsfbox{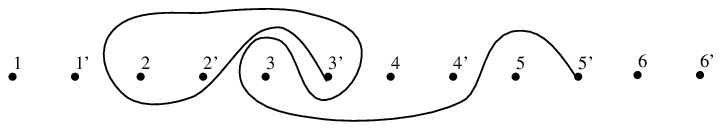}}}$ & & \\ \hline
$\vcenter{\hbox{\epsfbox{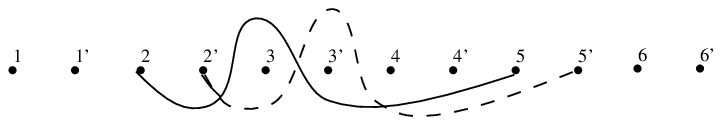}}}$ & $\rho^j_5 \rho^i_2 \tilde{z}_{25} \rho^{-i}_2 \rho^{-j}_5$ & 3\\
$\vcenter{\hbox{\epsfbox{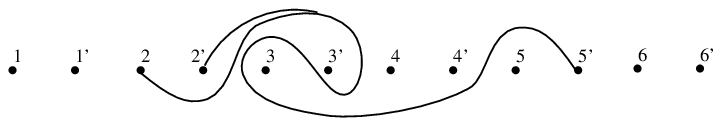}}}$ & & \\ \hline
$\vcenter{\hbox{\epsfbox{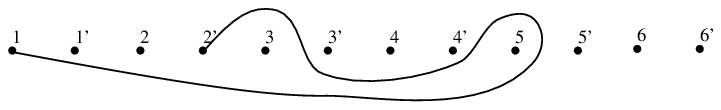}}}$  & $\alpha _1(m_1)$ & 1\\    \hline
$\vcenter{\hbox{\epsfbox{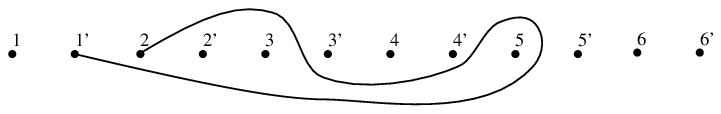}}}$   & $\alpha _2(m_1)$ & 1\\
    \hline
$\vcenter{\hbox{\epsfbox{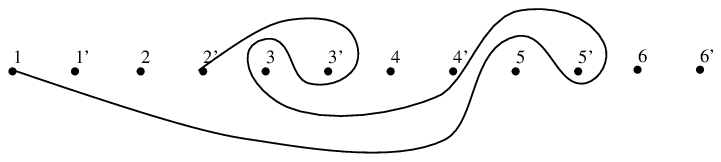}}}$   & $(\alpha _1(m_1))^{\rho^{-1}}$ & 1\\
   \hline
$\vcenter{\hbox{\epsfbox{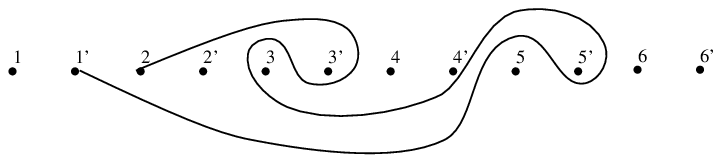}}}$     & $(\alpha _2(m_1))^{\rho^{-1}}$ & 1
\\   \hline
\end{tabular}
\end{center}
\hfill
$\Box$

\section{The computation of $H_{V_4}$}\label{sec:34}
\indent
We follow Figure 6 when computing the local braid monodromy 
for $V_4$, denoted by
$\vp^{(0)}_4$, but first we need to compute its factors which come from singularities in the neighbourhood of $V_4$.    We start with $\vp^{(6)}_4$.
 Then we compute $S^{(5)}_4$ and $\vp^{(5)}_4$ in the neighbourhood of $V_4$.  Notice that $\vp^{(5)}_4 = \vp^{(4)}_4 = \vp^{(3)}_4 = \vp^{(2)}_4$
since the fourth, fifth and sixth regenerations do not change the neighbourhood of $V_4$.
 Finally we compute $\vp^{(1)}_4$.  $\vp^{(1)}_4 =  \vp^{(0)}_4$ since the
eighth regeneration acts similarly on the neighbourhood of $V_4$.  Hence we get a factorized expression for $\vp^{(0)}_4$.

\noindent
{\bf Theorem 9.1.}  {\em In the neighbourhood of $V_4$, the local braid monodromy of $S^{(6)}_4$ around $V_4$ is given by:\\
$\vp^{(6)}_4 = (\uZ^2_{2'i})_{i=3,5,6,6'} \cdot
Z^4_{12} \cdot
\bZ^4_{2'4} \cdot
(\uZ^2_{2'i})_{i=3,5,6,6'}^{\uZ^{-2}_{i4}} \cdot
(Z_{22'})^{Z^2_{12}\uZ^2_{2'4}Z^2_{2'3}} \cdot
(Z_{66'})^{\uZ^{-2}_{36}Z^{-2}_{56}} \cdot
(\uZ^2_{i6'})^{Z^2_{12}}_{i = 1,4} \cdot
(\uZ^4_{36})^{Z^{-2}_{56}} \cdot
(\uZ^2_{i6})^{Z^{-2}_{56}Z^2_{12}}_{i = 1,4} \cdot
Z^4_{56} \cdot  (\Delta^2<1,3,4,5>)^{Z^2_{12}Z^{-2}_{56}}$, where $(Z_{22'})^{Z^2_{12}\uZ^2_{2'4}Z^2_{2'3}}$
is determined by \\ $\vcenter{\hbox{\epsfbox{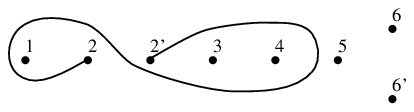}}}$ ,
$(Z_{66'})^{\uZumt_{36} Z^{-2}_{56}}$
is determined
by  $\vcenter{\hbox{\epsfbox{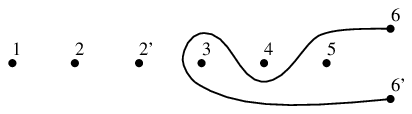}}}$
and the}  {\rm L.V.C.} {\em corresponding to
$\Delta^2<1,3,4,5>$ is given by}
$\vcenter{\hbox{\epsfbox{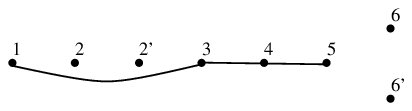}}}$ .

\noindent
\underline{Proof:} \ Let $\hat{L}_i, 1 \le i \le 6$,
 be the lines in $Z^{(8)}$ and in $Z^{(7)}$ intersecting in a point 4.
  Let $L_i = \pi^{(8)} (\hat{L}_i)$.  Then $L_i, 1 \le i \le 6$
, intersect in $V_4$.  $S^{(8)}$ (and $S^{(7)}$) around $V_4$ is
 $\bigcup\limits^6_{i=1} L_i$.  After the second regeneration, $S^{(6)}_4$ in a neighbourhood of $V_4$ is as follows:  the lines $L_2$ and $L_6$ are replaced by two conics $Q_2$ and $Q_6$, where $Q_2$ is tangent to the lines $L_1$ and $L_4$, $Q_6$ is tangent to $L_3$ and $L_5$ .
 This follows from the fact that $P_1 \cup P_5$ (resp. $P_{10} \cup P_{12})$ is regenerated to $R_2$ (resp. $R_5)$, and the tangency follows from Lemma 
\ref{lm:24}.  Thus, in a neighbourhood of $V_4$, $S^{(6)}_4$ is of the form: $S^{(6)}_4 =
L_1 \cup Q_2 \cup L_3 \cup L_4 \cup L_5 \cup Q_6$.

$S^{(6)}_4$ is shown in Figure 19.
\begin{figure}[h]
\epsfxsize=11cm 
\epsfysize=9cm 
\begin{minipage}{\textwidth}
\begin{center}
\epsfbox {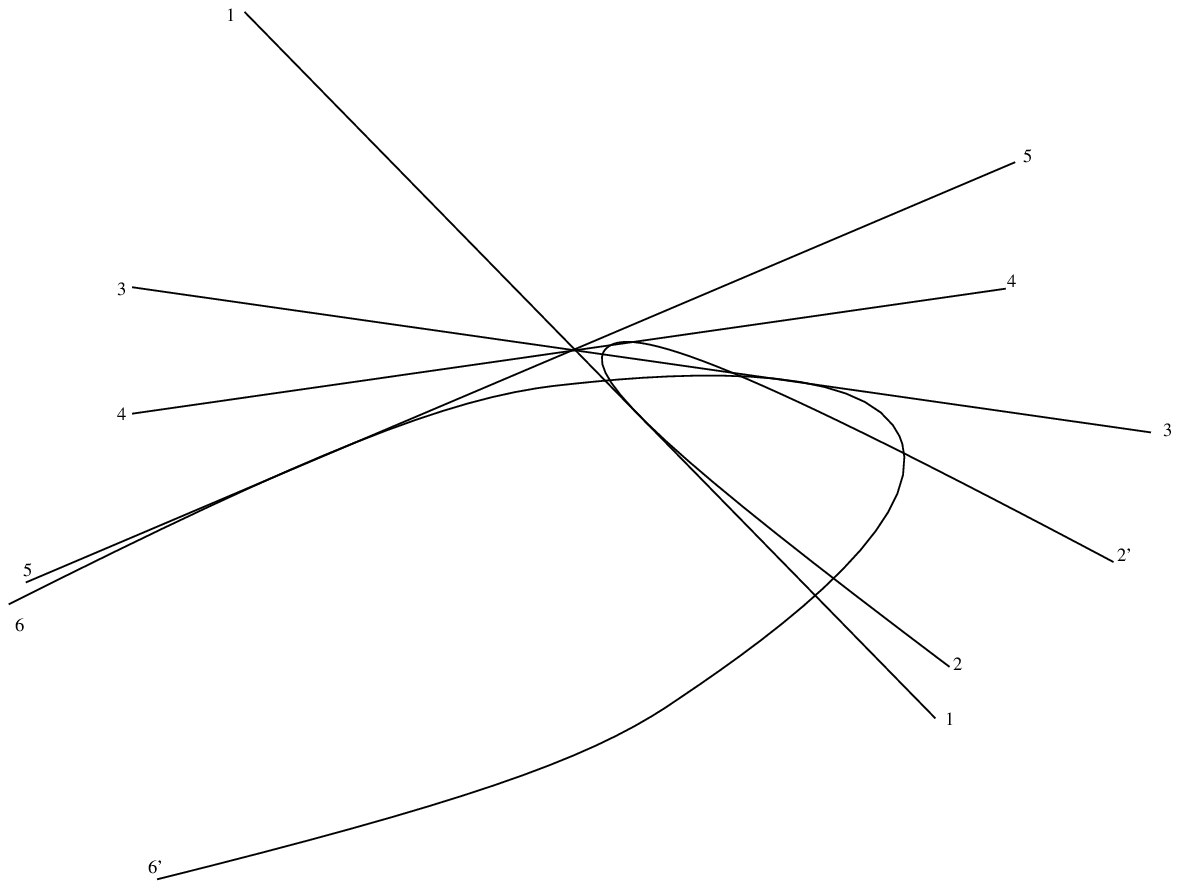}
\end{center}
\end{minipage}
\caption{}
\end{figure}

We devide the proof into parts:\\
\noindent
{\bf Proposition 9.1.1.} {\em Consider in $\C^2$ the following figure:  $C_1 = L_1 \cup Q_2 \cup L_3 \cup L_4$, where $L_i$ are lines, $i = 1,3,4, \ Q_2$ is a conic tangent to $L_1$ and $L_4$, it intersects $L_3$ and the 3 lines meet at one point (see Figure 20).

Then the braid monodromy of $C_1$ is given by:\\
$\vp_{C_1} = Z^2_{2'3} \cdot Z^4_{12} \cdot
\bZ^4_{2'4} \cdot
(Z^2_{2'3})^{Z^{-2}_{34}} \cdot
(Z_{22'})^{Z^2_{12}\uZ^2_{2'4}Z^2_{2'3}}
\cdot (\Delta^2<1,3,4>)^{Z^2_{12}}$,
 where
$(Z_{22'})^{Z^2_{12}\uZ^2_{2'4}Z^2_{2'3}}$
is determined by
$\vcenter{\hbox{\epsfbox{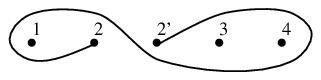}}}$
and the L.V.C. corresponding to $\Delta^2<1,3,4>$ is given by 
$\vcenter{\hbox{\epsfbox{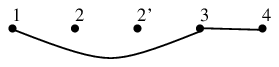}}}$}.

\begin{figure}[!h]
\epsfxsize=11cm 
\epsfysize=9cm 
\begin{minipage}{\textwidth}
\begin{center}
\epsfbox {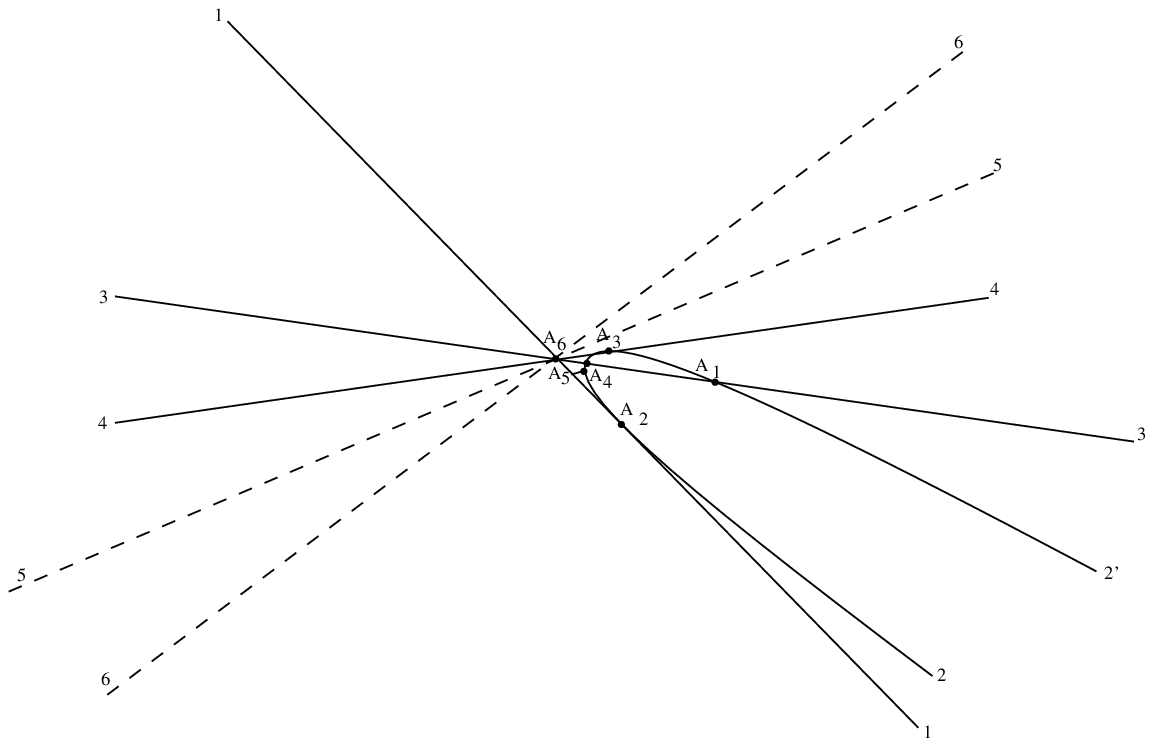}
\end{center}
\end{minipage}
\caption{}
\end{figure}

\underline{Proof}: \\
Let $\pi_1 : E \times D \rightarrow E$ be the projection to $E$.

Let $\{A_j\}^{6}_{j=1}$ be  singular points of $\pi_1$ as follows:\\
$A_1, A_4$ are intersection points of $Q_2$ with the line $L_3$. \\
$A_2, A_3$ are tangent points of $Q_2$
with the lines $L_1, L_4$ respectively.\\
$A_6$ is an intersection point of the lines $L_1, L_3, L_4$.\\
$A_5$ is a point of the type $a_1$ in $Q_2$.

Let $N,M, \{\ell(\g_j)\}^{6}_{j=1}$
as in Proposition 8.1.1.
Let $i = L_i \cap K$ for $i = 1,3,4$.
 Let $\{2,2'\} = Q_2 \cap K$.
 $K = \{1,2,2',3,4\}$, \st 1,2,2',3,4 are  real points.
Moreover: $1 < 2 < 2'< 3 < 4$.
Let $\beta_M$ be the diffeomorphism defined by: $\beta_M(i) = i$ for $i = 1,2, \ \beta_M(2') = 3, \ \beta_M (3) = 4, \  \beta_M(4) = 5$.  Moreover: deg$C_1 = 5$ and $\# K_\R (x) \geq 3 \ \ \forall x$.

We are looking for $\vp_M(\ell (\g_j))$ for $j = 1, \cdots , 6$.  We choose a $g$-base $\{\ell(\g_j)\}^{6}_{j=1}$ of $\pi_1(E - N,u)$, \st each path $\g_j$ is below the real line.

Put all data in the following table:

\vspace{-1cm}

\begin{center}
\begin{tabular}{cccc} \\
$j$ & $\lambda_{x_j}$ & $\epsilon_{x_j}
$ & $\delta_{x_j}$ \\ \hline
1 & $<3,4>$ & 2 & $\Delta<3,4>$\\
2 & $<1,2>$ & 4 & $\Delta^2<1,2>$\\
3 & $<4,5>$ & 4 & $\Delta^2<4,5>$\\
4 & $<3,4>$ & 2 & $\Delta<3,4>$\\
5 & $<2,3>$ & 1 & $\Delta^{\frac{1}{2}}_{I_2 \R}<2>$\\
6 & $<1,3>$ & 2 & $\Delta<1,3>$\\
\end{tabular}
\end{center}

For computations, we use the formulas in [Am2, Theorems 1.41, 1.44].

\noindent
\underline{Remark}:  $\beta_{x'_j} (K(x'_j)) = \{1,2,3,4,5\}$ for $
1 \leq j \leq 5$.\\
$\beta_{x'_6} (K(x'_6)) = \{1,2,3,4+i, 4-i\}$.\\

\medskip
\noindent
$(\xi_{x'_1}) \Psi_{\g'_1} = <3,4> \beta^{-1}_M = z_{2'3}$

\medskip
\noindent
$<3,4>$ \hspace{.3cm} $\vcenter{\hbox{\epsfbox{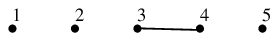}}}$ \hspace{.3cm} $\beta^{-1}_M$
 $\vcenter{\hbox{\epsfbox{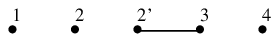}}}$
 \hspace{.3cm}  $\vp_M(\ell(\g_1)) = Z^2_{2'3}$ \\

\medskip

\noindent
$(\xi_{x'_2}) \Psi_{\g'_2} = <1,2> \Delta<3,4>\beta^{-1}_M = z_{12}$

\medskip
\noindent
$<1,2> \Delta <3,4>$ \ \  $\vcenter{\hbox{\epsfbox{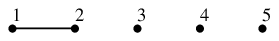}}}$ \ \ $\beta^{-1}_M$ \ \  $\vcenter{\hbox{\epsfbox{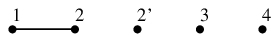}}}$ \ \  $\vp_M(\ell(\g_2)) = Z^4_{12}$ \\

\medskip

\noindent
$(\xi_{x'_3}) \Psi_{\g'_3}
 = <4,5>\Delta^2<1,2> \Delta<3,4> \beta^{-1}_M = \bar{z}_{2'4}$

\medskip
\noindent
$<4,5> \Delta^2 <1,2>$  \ \ $\vcenter{\hbox{\epsfbox{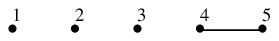}}}$ \ \ $\Delta<3,4>$ \ \ $\vcenter{\hbox{\epsfbox{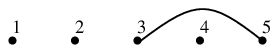}}}$ \
$\beta^{-1}_M$  \ \ $\vcenter{\hbox{\epsfbox{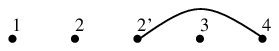}}}$ \\
$\vp_M(\ell(\g_3)) = \bZ^4_{2'4}$\\

\medskip

\noindent
$(\xi_{x'_4}) \Psi_{\g'_4}  = <3,4> \Delta^2<4,5>\Delta^2<1,2> \Delta<3,4>
\beta^{-1}_M = z_{2'3}^{Z^{-2}_{34}}$

\medskip

\noindent
$<3,4>$\ \ $\vcenter{\hbox{\epsfbox{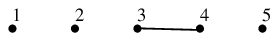}}}$ \ \
$\Delta^2<4,5>\Delta^2<1,2>$ \ \ $\vcenter{\hbox{\epsfbox{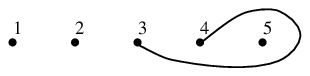}}}$ \\
$\Delta<3,4>$ \ \  $\vcenter{\hbox{\epsfbox{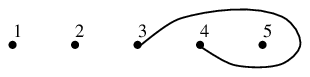}}}$ \ \ $\beta^{-1}_M$ \ \
$\vcenter{\hbox{\epsfbox{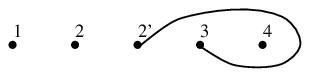}}}$ \ \
$\vp_M(\ell(\g_4)) = (Z^2_{2'3})^{Z^{-2}_{34}}$\\

\medskip

\noindent
$(\xi_{x'_5}) \Psi_{\g'_5}
 = <2,3> \Delta<3,4>\Delta^2<4,5> \Delta^2<1,2> \Delta<3,4> \beta^{-1}_M =
z_{22'}^{Z^2_{12} \uZ^2_{2'4} Z^2_{2'3}}$

\medskip

\noindent
$<2,3>$  \quad    $\vcenter{\hbox{\epsfbox{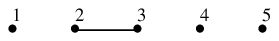}}}$ \quad   $\Delta<3,4>$
\quad   $\vcenter{\hbox{\epsfbox{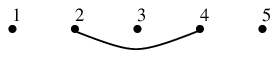}}}$ \ \
$\Delta^2<4,5>$  \quad    $\vcenter{\hbox{\epsfbox{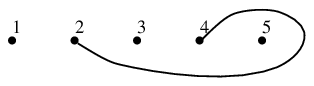}}}$ \\
$\Delta^2<1,2>$
\quad    $\vcenter{\hbox{\epsfbox{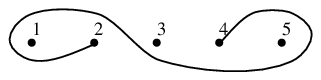}}}$ \ \
$\Delta<3,4>$ \quad   $\vcenter{\hbox{\epsfbox{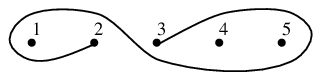}}}$ \ \
$\beta^{-1}_M$ \ \
 $\vcenter{\hbox{\epsfbox{v4019.ps}}}$ \\
$\vp_ M(\ell(\g_5)) = (Z_{22'})^{Z^2_{12}\uZ^2_{2'4}Z^2_{2'3}}$\\

\medskip

\noindent
$(\xi_{x'_6}) \Psi_{\g'_6}=  <1,3> \Delta^{\frac{1}{2}}_{I_2 \R}<2>
  \Delta<3,4> \Delta^2<4,5>\Delta^2<1,2> \Delta<3,4>  \beta^{-1}_M
= \\
(\Delta<1,3,4>)^{Z^2_{12}}$

\medskip

\noindent
$<1,3>$ \ $\vcenter{\hbox{\epsfbox{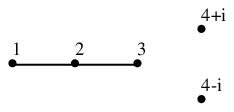}}}$ \ $\Delta^{\frac{1}{2}}_{I_2
\R}<2> $ \ $\vcenter{\hbox{\epsfbox{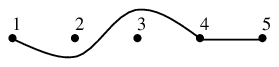}}}$ \
$\Delta<3,4>$ \  $\vcenter{\hbox{\epsfbox{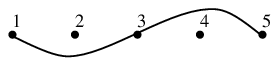}}}$ \\
\noindent
$\Delta^2<4,5>$   $\vcenter{\hbox{\epsfbox{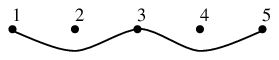}}}$ $\Delta^2<1,2>$
 $\vcenter{\hbox{\epsfbox{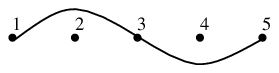}}}$    \quad
$\Delta<3,4>$  \quad  $\vcenter{\hbox{\epsfbox{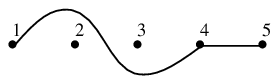}}}$
\\
$\beta^{-1}_M$  \quad $\vcenter{\hbox{\epsfbox{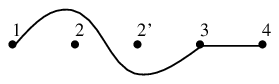}}}$ \quad
$\vp_M(\ell(\g_6)) = (\Delta^2<1,3,4>)^{Z^2_{12}}$

\medskip
\noindent
{\bf Remark 9.1.2.}
{\em In a similar way as in  Remark 8.1.2, the changes are:\\
(a) $(\Delta^2<1,3,4>)^{Z^2_{12}}$ is replaced by
$(\Delta^2<1,3,4,5,6>)^{Z^2_{12}}$ with the corresponding L.V.C.
which appears as follows:}
{\em $\vcenter{\hbox{\epsfbox{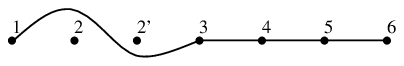}}}$ .\\
(b) $Z^2_{2'3}$ is replaced by $(\uZ^2_{2'i})_
{i = 3,5,6}$ as follows:
 $\vcenter{\hbox{\epsfbox{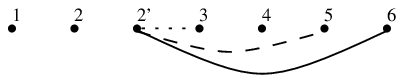}}}$ .\\
(c) $ (Z^2_{2'3})^{Z^{-2}_{34}}$ is replaced by
$(\uZ^2_{2'i})_
{i = 3,5,6}^{\uZ^{-2}_{i4}}$ as follows:
 $\vcenter{\hbox{\epsfbox{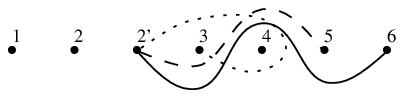}}}$ .}

\noindent
{\bf Proposition 9.1.3.} {\em Consider in $\C^2$ the following figure: $C_2 = L_1 \cup L_3 \cup L_5 \cup  Q_6$, where $L_i$ are lines, $i = 1,3,5, \ Q_6$ is a conic tangent to $L_3$ and $L_5$, it intersects $L_1$  and the 3 lines meet at one point (see Figure 21).

Then the braid monodromy of $C_2$ is given by:\\
$\vp_{C_2} =
( Z_{66'})^{\uZ^{-2}_{36}Z^{-2}_{56}}
\cdot \uZ^2_{16'} \cdot
(\uZ^4_{36})^{Z^{-2}_{56}} \cdot
(\uZ^2_{16})^{Z^{-2}_{56}} \cdot
(\Delta^2<1,3,5>)^{Z^{-2}_{56}} \cdot Z^4_{56}$ ,
 where
$(Z_{66'})^{\uZ^{-2}_{36}Z^{-2}_{56}}$  is determined by
$\vcenter{\hbox{\epsfbox{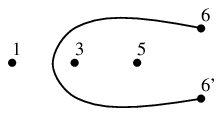}}}$
  and the {\rm   L.V.C.} corresponding to $\Delta^2<1,3,5>$ is given by
$\vcenter{\hbox{\epsfbox{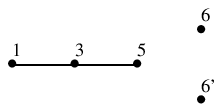}}}$}.

\begin{figure}[!h]
\epsfxsize=11cm 
\epsfysize=9cm 
\begin{minipage}{\textwidth}
\begin{center}
\epsfbox {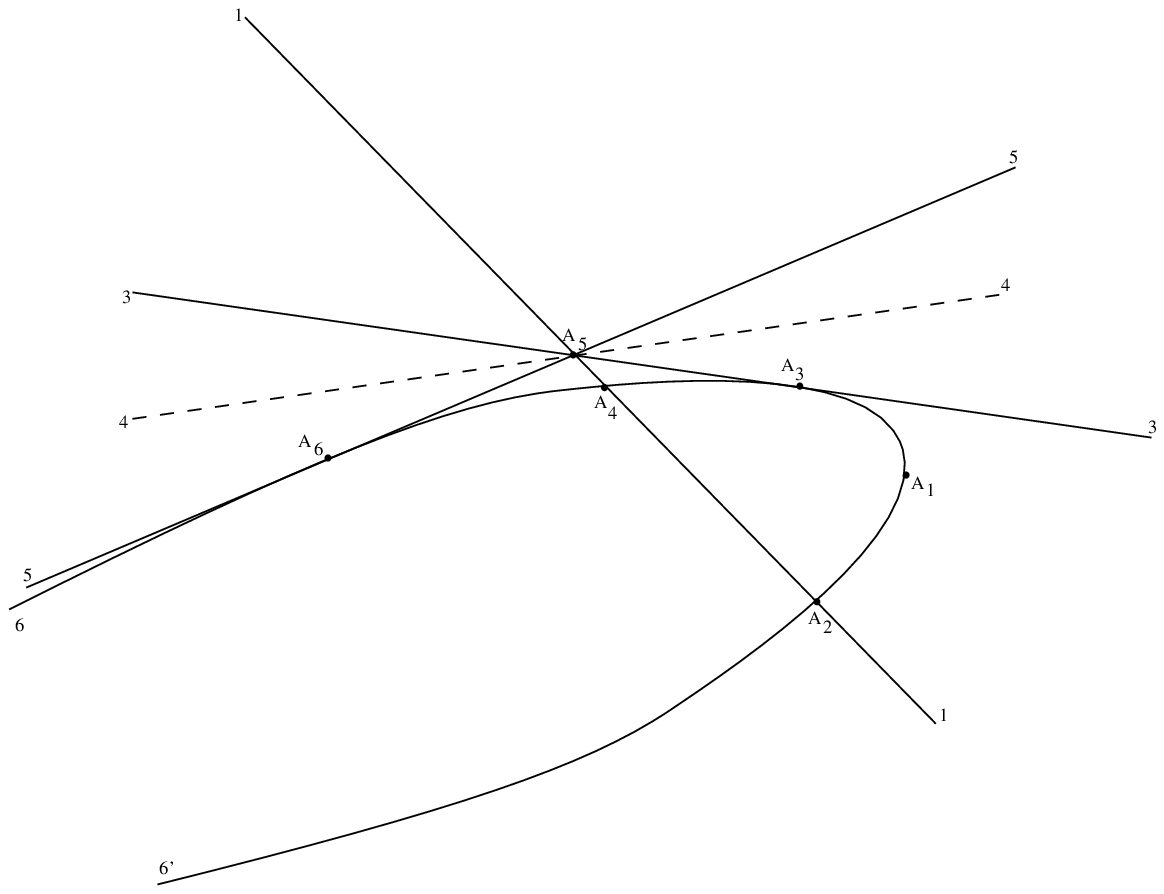}
\end{center}
\end{minipage}
\caption{}\end{figure}

\underline{Proof:} Let $\pi_1: E \times D \rightarrow E$ be the projection to $E$.

Let $\{A_j\}^{6}_{j=1}$ be singular points of $\pi_1$ as follows:\\
$A_2, A_4$ are intersection points of $Q_6$ with the line $L_1$.\\
$A_3, A_6$  are tangent points of $Q_6$ with the lines $L_3,L_5$ respectively. \\
$A_5$ is an intersection point of the lines $L_1, L_3, L_5$.\\
$A_1$ is a point of the type $a_2$ in $Q_6$.

Let $K,N,M, \{\ell (\g_j)\}^{6}_{j=1}$ be as in Proposition 8.1.3.
Let $i = L_i \cap K$, $i = 1,3,5, \ \{6,6'\} = Q_6 \cap K$.
$K = \{1,3,5,6,6'\}$, \st $1,3,5$ are real points and
$6, 6'$ are complex ones and $1 < 3 < 5 < Re(6), \ Re(6) = Re(6')$.
Let $\beta_M$ be the diffeomorphism defined by
$\beta_M(1) = 1 , \   \beta_M(3) = 2, \ \beta_M(5) = 3, \  \beta_M(6) = 4+i, \  \beta_M(6') = 4-i$.  deg$C_2 = 5$ and
 $\# K_\R (x) \geq 3 \ \forall x$.

Put all data in the following table:

\vspace{-1cm}
\begin{center}
\begin{tabular}{cccc} \\
$j$ & $\lambda_{x_j}$ & $\epsilon_{x_j}
$ & $\delta_{x_j}$ \\ \hline
1 & $P_2$ & 1 & $\Delta^{\frac{1}{2}}_{\R I_2}<2>$\\
2 & $<1,2>$ & 2 & $\Delta<1,2>$\\
3 & $<3,4>$ & 4 & $\Delta^2<3,4>$\\
4 & $<2,3>$ & 2 & $\Delta<2,3>$\\
5 & $<3,5>$ & 2 & $\Delta<3,5>$\\
6 & $<2,3>$ & 4 & $\Delta^2<2,3>$\\
\end{tabular}
\end{center}

\noindent
\underline{Remark}:  $\beta_{x'_j} (K(x'_j)) = \{1,2,3,4,5\}$ for $2 \leq j \leq 6$.\\
$\beta_{x'_1}(K(x'_1)) = \{1,2,3,4+i, 4-i\}$.\\

\medskip
\noindent
$(\xi_{x'_1}) \Psi_{\g'_1} = P_2 \beta^{-1}_M = z_{66'}^{\uZ^{-2}_{36}Z^{-2}_{56}}$

\noindent
$P_2$ $\vcenter{\hbox{\epsfbox{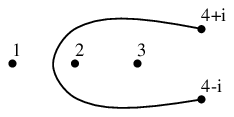}}}$ \ \
$\beta^{-1}_M$ $\vcenter{\hbox{\epsfbox{v4032.ps}}}$ \ \
$\vp_M(\ell(\g_1)) = (Z_{66'})^{\uZ^{-2}_{36}Z^{-2}_{56}}$\\

\noindent
$(\xi_{x'_2}) \Psi_{\g'_2} = <1,2> \Delta^{\frac{1}{2}}_{\R I_2}<2>\beta^{-1}_M = \underline{ z}_{16'}$

\noindent
$<1,2> \vcenter{\hbox{\epsfbox{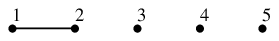}}}$ \  $\Delta^{\frac{1}{2}}_{\R I_2}<2> \ \vcenter{\hbox{\epsfbox{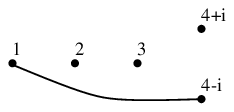}}}$ \ \
$\beta^{-1}_M$ \ \ $\vcenter{\hbox{\epsfbox{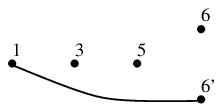}}}$\ \
$\vp_M(\ell(\g_2)) = \uZ^2_{16'}$\\

\noindent
$(\xi_{x'_3}) \Psi_{\g'_3}
 = <3,4> \Delta<1,2> \Delta^{\frac{1}{2}}_{\R I_2}<2> \beta^{-1}_M =
\underline{z}_{36}^{Z^{-2}_{56}}$

\noindent
$<3,4> \Delta<1,2>
 \vcenter{\hbox{\epsfbox{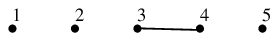}}}$ \ \
$\Delta^{\frac{1}{2}}_{\R I_2}<2> \vcenter{\hbox{\epsfbox{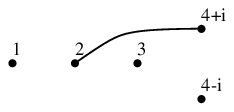}}}$ \ \
$\beta^{-1}_M$ $\vcenter{\hbox{\epsfbox{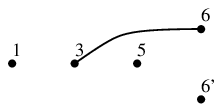}}}$ \ \
$\vp_M(\ell(\g_3)) = (\uZ^4_{36})^{Z^{-2}_{56}}$\\

\noindent
$(\xi_{x'_4}) \Psi_{\g'_4}  =  <2,3>
 \Delta^2<3,4> \Delta<1,2>\Delta^{\frac{1}{2}}_{\R I_2}<2> \beta^{-1}_M =
\underline{z}_{16}^{Z^{-2}_{56}}$

\noindent
$<2,3> \vcenter{\hbox{\epsfbox{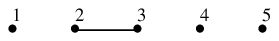}}}$ \ \
$ \Delta^2<3,4> \ \vcenter{\hbox{\epsfbox{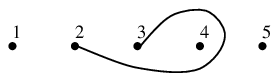}}}$ \ \
$ \Delta<1,2> \ \vcenter{\hbox{\epsfbox{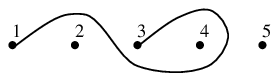}}}$ \ \
$\Delta^{\frac{1}{2}}_{\R I_2}<2>$ \ $\vcenter{\hbox{\epsfbox{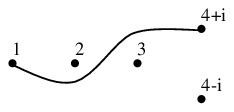}}}$ \ \
$\beta^{-1}_M$ $\vcenter{\hbox{\epsfbox{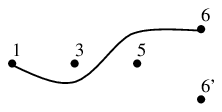}}}$ \ \
$\vp_M(\ell(\g_4)) = (\uZ^2_{16})^{Z^{-2}_{56}}$\\

\noindent
$(\xi_{x'_5}) \Psi_{\g'_5}  =   <3,5> \Delta<2,3>
\Delta^2<3,4>\Delta<1,2> \Delta^{\frac{1}{2}}_{\R I_2}<2> \beta^{-1}_M = \\
(\Delta<1,3,5>)^{Z^{-2}_{56}}$

\noindent
$<3,5>$  $\vcenter{\hbox{\epsfbox{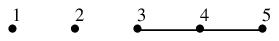}}}$ \ \
$\Delta<2,3>$  $\vcenter{\hbox{\epsfbox{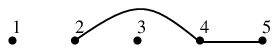}}}$ \ \
$\Delta^2<3,4>$
 $\vcenter{\hbox{\epsfbox{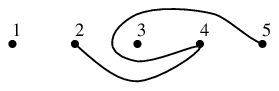}}}$ \\
$\Delta<1,2>$
 $\vcenter{\hbox{\epsfbox{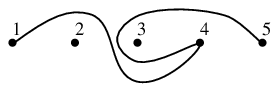}}}$ \ \
$ \Delta^{\frac{1}{2}}_{\R I_2}<2>$
 $\vcenter{\hbox{\epsfbox{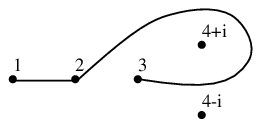}}}$ \ \
$\beta^{-1}_M$ $\vcenter{\hbox{\epsfbox{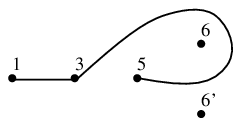}}}$ \\
$\vp_M(\ell(\g_5)) = (\Delta^2<1,3,5>)^{Z^{-2}_{56}}$\\

\noindent
$(\xi_{x'_6}) \Psi_{\g'_6}  =  <2,3> \Delta<3,5>\Delta<2,3>
 \Delta^2<3,4> \Delta<1,2> \Delta^{\frac{1}{2}}_{\R I_2}<2>
\beta^{-1}_M = z_{56}$

\noindent
$  <2,3>  \vcenter{\hbox{\epsfbox{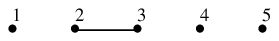}}}$ \ \
$ \Delta<3,5> \ \vcenter{\hbox{\epsfbox{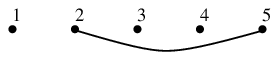}}}$ \ \
$ \Delta<2,3> \ \vcenter{\hbox{\epsfbox{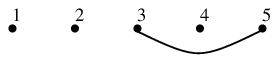}}}$ \ \
$ \Delta^2<3,4>\Delta<1,2> \ \vcenter{\hbox{\epsfbox{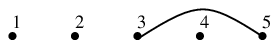}}}$ \ \
$\Delta^{\frac{1}{2}}_{\R I_2}<2> \vcenter{\hbox{\epsfbox{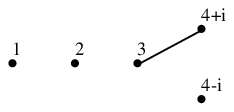}}}$ \ \
$\beta^{-1}_M$ $\vcenter{\hbox{\epsfbox{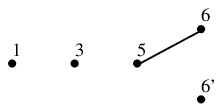}}}$ \\
$\vp_M(\ell(\g_6)) = Z^4_{56}$

\medskip
\noindent
{\bf Remark 9.1.4.}
{\em Similarly to Remark 9.1.2, the changes are:\\
(a) $(\Delta^2<1,3,5>)^{Z^{-2}_{56}}$ is replaced by $(\Delta^2<1,3,4,5>)^{Z^{-2}_{56}}$ with the corresponding {\rm L.V.C.}:}
{\em $\vcenter{\hbox{\epsfbox{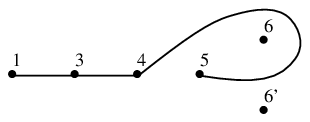}}}$ .\\
(b) $\uZ^2_{16'}$ is replaced by $(\uZ^2_{i6'})_
{i = 1,4}$ with the following {\rm L.V.C}.:
 $\vcenter{\hbox{\epsfbox{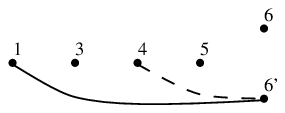}}}$ .\\
(c) $(\uZ^2_{16})^{Z^{-2}_{56}}$ is replaced by
$(\uZ^2_{i6})_{i = 1,4}^{Z^{-2}_{56}}$ with the following {\rm L.V.C.}:
 $\vcenter{\hbox{\epsfbox{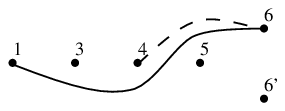}}}$ .}

\noindent
\underline{Proof of Theorem 9.1:}

Until now we computed the braid monodromy of all singularities.  Each one of the intersection points of $L_6 \cap Q_2$ is replaced by 2 intersection points $(\subseteq Q_2 \cap Q_6)$ which are close to each other.  So $\pi^{-1} (M) \cap S^{(6)}_4 =
\{1,2,2',3,4,5,6,6'\}$.

The changes are:\\
(I) Similar to (I) in the proof of Theorem 8.1.\\
(II) $\uZ^2_{2'6}$ is replaced by
$(\uZ^2_{2'i})_{i = 6,6'}$.  These braids correspond to the paths
$(\underline{z}_{2'i})_{i = 6,6'}$.
$(\uZ^2_{2'6})^{\uZ^{-2}_{46}}$ is replaced by $(\uZ^2_{2'i})_{i=6,6'}^{\uZ^{-2}_{4i}}$.
These braids correspond to the paths $(\underline{z}_{2'i})_{i = 6,6'}^{\uZ^{-2}_{4i}}$.\\
(III) In a similar proof as in Theorem 8.1, we have to conjugate all braids from Proposition 9.1.3 and Remark 9.1.4 by $Z^2_{12}$.

According to these changes, we present the list of braids:


\begin{tabbing}
$\vcenter{\hbox{\epsfbox{v4051.ps}}}$vvdddvv \= dd corresponding to
$ Z^4_{34})^2_{i = 2,5,6,6'}$ \kill
$\vcenter{\hbox{\epsfbox{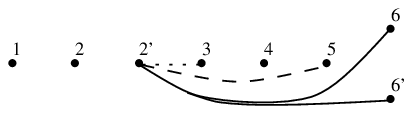}}}$ \> corresponding to
$ (\uZ^2_{2'i})_{i=3,5,6,6'}        $ \\ [.5cm]
$\vcenter{\hbox{\epsfbox{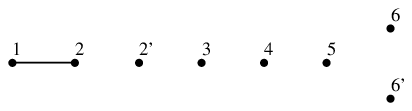}}}$ \> corresponding to
$Z^4_{12}$\\[.5cm]
$\vcenter{\hbox{\epsfbox{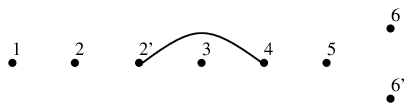}}}$ \> corresponding to
$\bZ^4_{2'4}$\\[.5cm]
$\vcenter{\hbox{\epsfbox{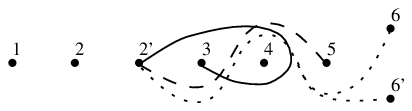}}}$ \> corresponding to
$(\uZ^2_{2'i})_{i = 3,5,6,6'}^{\uZ^{-2}_{i4}} $\\[.5cm]
$\vcenter{\hbox{\epsfbox{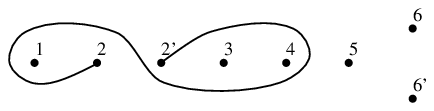}}}$ \> corresponding to
$(Z_{22'})^{Z^2_{12}\uZ^2_{2'4}Z^2_{2'3}} $\\[.5cm]
$\vcenter{\hbox{\epsfbox{v4064.ps}}}$ \> corresponding to
$(Z_{66'})^{\uZ^{-2}_{36}Z^{-2}_{56}}$\\[.5cm]
$\vcenter{\hbox{\epsfbox{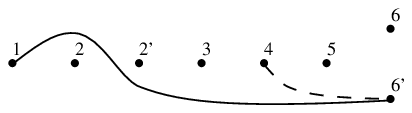}}}$ \> corresponding to
$(\uZ^2_{i6'})^{Z^{2}_{12}}_{i = 1,4}$\\[.5cm]
$\vcenter{\hbox{\epsfbox{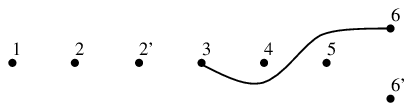}}}$ \> corresponding to $
(\uZ^4_{36})^{Z^{-2}_{56}}$\\[.5cm]
$\vcenter{\hbox{\epsfbox{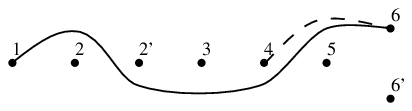}}}$ \> corresponding to
$(\uZ^2_{i6})^{Z^{-2}_{56}Z^2_{12}}_{i=1,4}$ \\[.5cm]
$\vcenter{\hbox{\epsfbox{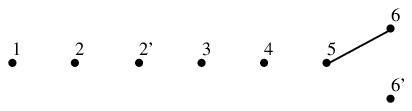}}}$ \> corresponding to $
Z^4_{56}$\\[.5cm]
$\vcenter{\hbox{\epsfbox{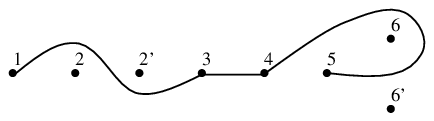}}}$ \> corresponding to $(\Delta^2<1,3,4,5>)^{Z^2_{12}Z^{-2}_{56}}$ \hspace{3.5cm} $\Box$
\end{tabbing}

Now we compute $\vp^{(5)}_4$.  Recall that $\vp^{(5)}_4 = \vp^{(4)}_4 =\vp^{(3)}_4 = \vp^{(2)}_4$, so in fact we computed $\vp^{(2)}_4$.

\noindent
{\bf Theorem 9.2.} {\em Let $\omega = L_1 \cup L_3 \cup L_4 \cup L_5$ in $S^{(6)}_4$. \\
(a) In a small neighbourhood $U_2$ of $\omega, T^{(2)} = S^{(2)}_4 \cap \ U_2$
resembles Figure 23, i.e., the singularities of $T^{(2)}$ are 4
tangency points, 4 nodes  and 2 branch points.\\
(b) The local braid monodromy of $S^{(2)}_4$ in that neighbourhood
is presented by $F_1(F_1)^{\rho^{-1}}$, where:\\
 $F_1 = Z^4_{1'3} \cdot Z^4_{45} \cdot (Z_{34})^{Z^2_{45}Z^2_{1'3}} \cdot
(\uZ^2_{1'5})^{Z^2_{1'3}}
\cdot \uZ^2_{15} \ .\\
\rho = Z_{11'} Z_{55'}$ and $(z_{34})^{Z^2_{45} Z^2_{1'3}} =
\vcenter{\hbox{\epsfbox{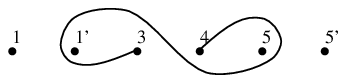}}}$ } .

\noindent
\underline{Proof}:

(a) Let $\hat{\omega} = R_1 \cap R_2 \cap R_5 \cap R_7$ in $Z^{(6)}$.

For local analysis we can use holomorphic coordinates in a
neighbourhood $\hat{U}$ of $\hat{\omega}$ in $\C \P^n$.  This allows us to
consider $\hat{U}$ as a subvariety of a neighbourhood of the origin in $\C^4$
with coordinates $X, Y, Z, T$ defined by the following system of equations
$
\hat{U} : \left\{ \begin{array}{l}
XT = 0 \\
YZ = 0 \end{array} \right .$.

By abuse of notation, $\hat{L}_1$  is called the $x$-axis (see Figure 22), 
and
then: $R_1 = \{T=0 \ ,\ \ Y=0\} \ , \ R_2 = \{Z=0, T=0\}\ , \
 R_5 = \{X=0\ , \ Y=0 \} \ , \
R_7 = \{Z=0, X=0\}$.

\begin{figure}[h]
\epsfxsize=5cm 
\epsfysize=5cm 
\begin{minipage}{\textwidth}
\begin{center}
\epsfbox {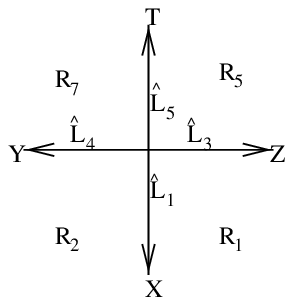}
\end{center}
\end{minipage}
\caption{}
\end{figure}

Now, we consider Lemma 6 in  [MoTe8] to obtain Figure 23.
%
\begin{figure}[h]
\epsfxsize=11cm 
\epsfysize=9cm 
\begin{minipage}{\textwidth}
\begin{center}
\epsfbox {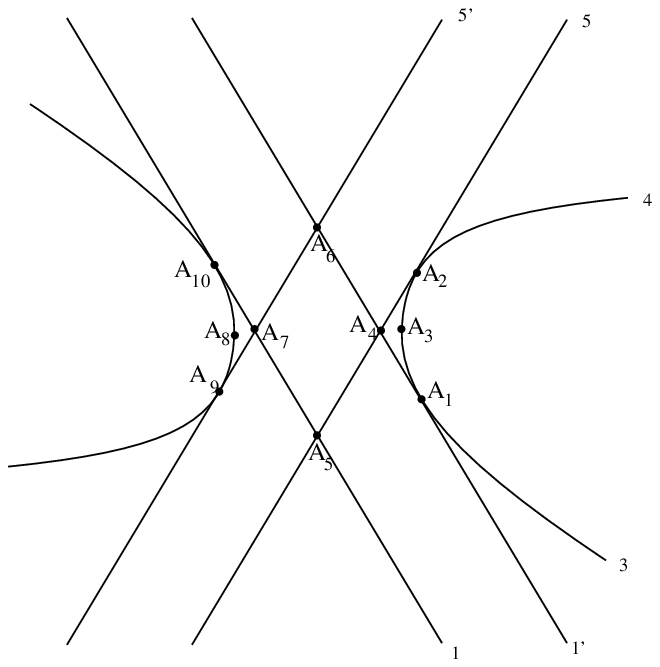}
\end{center}
\end{minipage}
\caption{}\end{figure}

(b) $S^{(2)}_4 = L_1 \cup L_{1'} \cup h_{34} \cup L_5 \cup L_{5'}$, see
Figure 23.
Here $K = \{1,1', 3,4, 5,5'\}$.  We apply the proof of Proposition 2.5 in [MoTe5]
(where $K = \{1,2, \cdots , 6\})$ to obtain the above $F_1 (F_1)^{\rho^{-1}}$ for
 $\rho^{-1} = Z^{-1}_{55'} Z^{-1}_{11'}$.\\
\noindent
The paths corresponding to the factors in $F_1$:\\
\noindent
$\vcenter{\hbox{\epsfbox{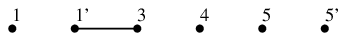}}}$ \ corresponding to $Z^4_{1'3}$ \ ; \
$\vcenter{\hbox{\epsfbox{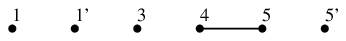}}}$ \
corresponding to $Z^4_{45}$\\
\noindent
$\vcenter{\hbox{\epsfbox{4v070.ps}}}$
 corresponding to $(Z_{34})^{Z^2_{45}Z^2_{1'3}}$  ;
$\vcenter{\hbox{\epsfbox{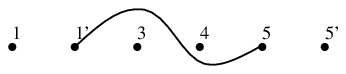}}}$  corresponding to $(\uZ^2_{1'5})^{Z^2_{1'3}}$\\
$\vcenter{\hbox{\epsfbox{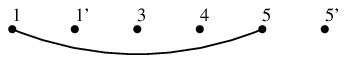}}}$ \quad  corresponding to $\uZ^2_{15}$ \ .
\hfill
$\Box$

\vspace{.5cm}
\noindent
{\bf Proposition 9.3.} {\em The local braid monodromy for $S^{(2)}_4$ around $V_4$ is obtained from formula $\vp^{(6)}_4$ in Theorem 9.1 by the following replacements:\\
(i) Consider the following $( \ )^\ast$ as conjugations by the braids, induced from the following motions: \\
$(\quad \quad)^{Z^{-2}_{4,55'}}$ \ $\vcenter{\hbox{\epsfbox{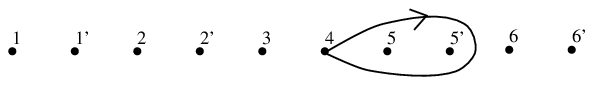}}}$ ;\\
$(\quad \quad)^{Z^{2}_{11',2}}$ ; $\vcenter{\hbox{\epsfbox{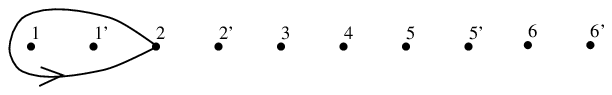}}}$ ;\\
$(\quad \quad)^{Z^{-2}_{55',6}}$ \ $\vcenter{\hbox{\epsfbox{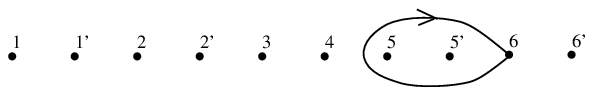}}}$ .

All the other conjugations do not change, since indices 1, 5 are not
involved.

\noindent
(ii) $\Delta^2<1,3,4,5>$ is replaced by  $F_1(F_1)^{\rho^{-1}}$.\\
(iii) Each of the degree 4 factors in $\vp^{(6)}_4$ that involves index 1 or 5 is replaced by 3 cubes as in the third regeneration rule.\\
(iv) Each of the degree 2 factors in $\vp^{(6)}_4$ that involves index 1 or 5 is replaced by 2 degree 2 factors, where : 1 and 1' are replacing 1; 5 and 5' are replacing 5.

We call the formula that was obtained: (2)$_4$.}

\noindent
\underline{Proof}: A similar proof as the proof of Proposition 8.3, but here changes are applied on indices 1 and 5.
\hfill
$\Box$

Now we compute $\vp^{(1)}_4$.  Recall that $\vp^{(1)}_4 = \vp^{(0)}_4$, so in fact we compute  $\vp^{(0)}_4$.

\noindent
{\bf Theorem 9.4.}  {\em In the notation of Theorem 9.2:  let $T^{(0)}$ be the curve
 obtained from $T^{(2)}$ in the regeneration $Z^{(2)} \leadsto Z^{(0)}$. \\
(i) Then the local braid monodromy of $T^{(0)}$ is $\hat{F}_1 (\hat{F}_2)$, where:\\
$\hat{F}_1 = \uZ^3_{1',33'} \cdot \uZ^3_{44',5}\cdot Z_{34'} \cdot Z_{3'4}
\cdot (\uZ^2_{1'5})^{\uZ^2_{1',33'}}
 \cdot \uZ^2_{15} $, \\
$\hat{F}_2 = (\uZ^3_{1',33'})^{\rho^{-1}} \cdot(\uZ^3_{44',5})^{\rho^{-1}}
\cdot Z_{34'}^{\rho^{-1}} \cdot Z_{3'4}^{\rho^{-1}}
\cdot\left((\uZ^2_{1'5})^{\uZ^2_{1',33'}}\right)^{\rho^{-1}} \; \; \cdot
(\uZ^2_{15})^{\rho^{-1}}$. \\
(ii) The singularities of the $x$-projection of \ $T^{(0)}$ are those arising from $T^{(2)}$ by the regeneration rules, namely: 4 nodes that exist in $T^{(2)}
\ , \ 3 \times 4$ cusps arising from 4 tangency points in $T^{(2)}\ , \ 2 \times 2$ branch points from 2 branch points of $T^{(2)}$.\\
(iii) The braid monodromy of $T^{(0)}$ is $\hat{F}_1 (\hat{F}_2)$. \\
(iv) $\hat{F}_1 (\hat{F}_2) = \Delta^{-2}_8 Z^{-2}_{11'} Z^{-2}_{33'} Z^{-2}_
{44'} Z^{-2}_{55'}$}.

\underline{Proof}:  A similar proof as the proof of Theorem 8.4 but replace $i = 1,2,3,5$ with $i = 1,3,4,5$.
\vspace{-1.2cm}
\begin{flushright}
$\Box$
\end{flushright}
\vspace{-.4cm}

The local braid monodromy of $S^{(0)}_4$ in $U_2$ is $\hat{F}_1 (\hat{F}_1)^{\rho^{-1}}$.  Applying regeneration rules on $F_1(F_1)^{\rho^{-1}}$, we obtain 
$\hat{F}_1 (\hat{F}_1)^{\rho^{-1}}$.

{\bf Corollary 9.5.} {\em The paths corresponding to $\hF_1(\hF_1)^{\rho^{-1}}$ (without their conjugation). \\
(a) \underline{ The paths corresponding to the factors in $\hF_1$}:

\noindent
$\vcenter{\hbox{\epsfbox{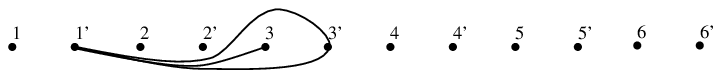}}}$ \quad
corresponding to $\uZ^3_{1',33'}$\\
\noindent
$\vcenter{\hbox{\epsfbox{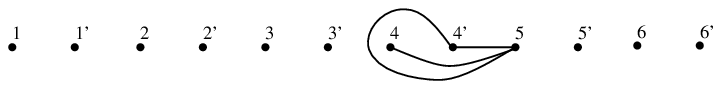}}}$ \quad  corresponding to $\uZ^3_{44',5}$\\
\noindent
$\vcenter{\hbox{\epsfbox{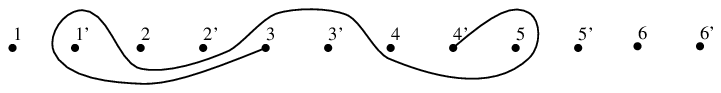}}}$ \quad   corresponding to $Z_{34'}$\\
\noindent
$\vcenter{\hbox{\epsfbox{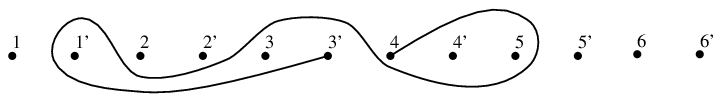}}}$ \quad  corresponding to  $Z_{3'4}$\\
\noindent
$\vcenter{\hbox{\epsfbox{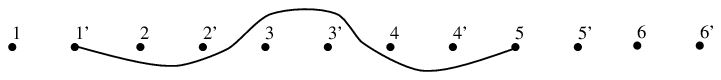}}}$ \quad
corresponding to $(\uZ^2_{1'5})^{\uZ^2_{1',33'}}$\\
\noindent
$\vcenter{\hbox{\epsfbox{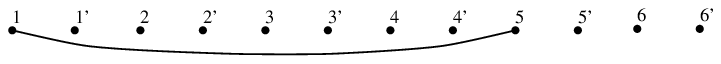}}}$ \quad  corresponding to $\uZ^2_{15}$

\medskip

\noindent
(b) \underline{ The paths corresponding to the factors in $(\hF_1)^{\rho^{-1}}$ for
$\rho^{-1} = Z^{-1}_{55'} Z^{-1}_{11'}$}: \\
\noindent
$\vcenter{\hbox{\epsfbox{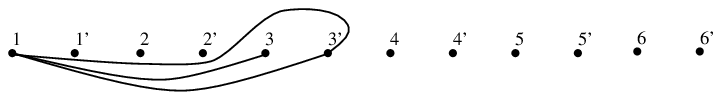}}}$ \quad  corresponding to $(\uZ^3_{1',33'})^{\rho^{-1}}$\\
\noindent
$\vcenter{\hbox{\epsfbox{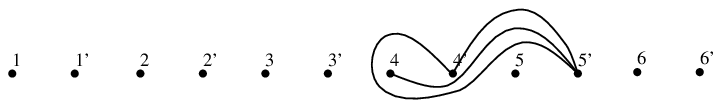}}}$ \quad  corresponding to $(\uZ^3_{44',5})^{\rho^{-1}}$\\
\noindent
$\vcenter{\hbox{\epsfbox{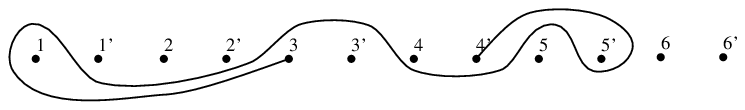}}}$ \quad  corresponding to $Z_{34'}^{\rho^{-1}}$\\
\noindent
$\vcenter{\hbox{\epsfbox{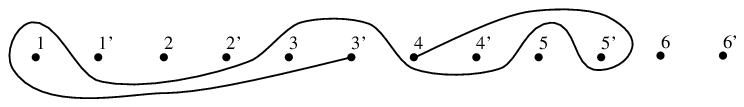}}}$ \quad  corresponding to  $Z_{3'4}^{\rho^{-1}}$\\
\noindent
$\vcenter{\hbox{\epsfbox{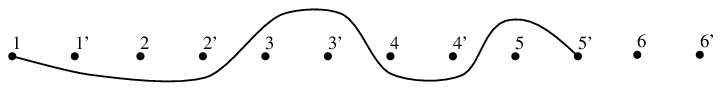}}}$ \quad corresponding to $\left((\uZ^2_{1'5})^{\uZ^2_{1',33'}}\right)^{\rho^{-1}}$\\
\noindent
$\vcenter{\hbox{\epsfbox{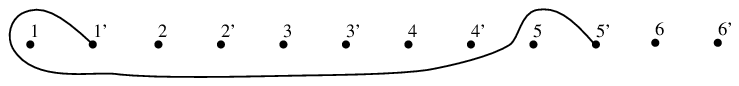}}}$ \quad  corresponding to $(\uZ^2_{15})^{\rho^{-1}}$

}
\noindent
{\bf Theorem 9.6.} {\em  The local braid monodromy of $S^{(0)}_4$ around $V_4$, denoted by $H_{V_4}$, equals:\\
$H_{V_4}$ = $(\uZ^2_{2'i})_{i=33',55',6,6'} \cdot
Z^3_{11',2} \cdot
\bZ^3_{2',44'}
\cdot
(\uZ^2_{2'i})^{\uZ^{-2}_{i,44'}}_{i=33',55',6,6'}
\cdot
(Z_{22'})^{Z^2_{11',2}\uZut_{2',44'}Z^2_{2',33'}}\cdot
(Z_{66'})^{\uZumt_{33',6}  Z^{-2}_{55',6}}\cdot
(\uZ^2_{i6'})^{Z^2_{11',2} }_{i = 11',44'}
\cdot
(\uZ^3_{33',6})^{Z^{-2}_{55',6}}
\cdot
(\uZ^2_{i6})^{Z^{-2}_{55',6}Z^2_{11',2}}_{i = 11',44'}
\cdot Z^3_{55',6}
\cdot (\hF_1(\hF_1)^{\rho^{-1}})^{Z^2_{11',2}Z^{-2}_{55',6}}$ ,
 where the paths corresponding to these braids are (the paths corresponding to $\hF_1(\hF_1)^{\rho^{-1}}$ are above):

\begin{tabbing}
(1) \  $\vcenter{\hbox{\epsfbox{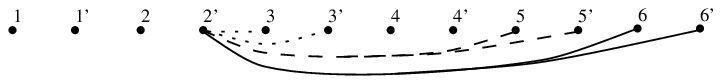}}}$ vvv \=   vvv corresponding to $ (\uZ^2_{2'i})_{33',55',6,6'}$ \kill
(1) \   $\vcenter{\hbox{\epsfbox{v4072.ps}}}$  \>
corresponding to
$(\uZ^2_{2'i})_{i=33',55',6,6'}$ \\[.5cm]
(2) \   $\vcenter{\hbox{\epsfbox{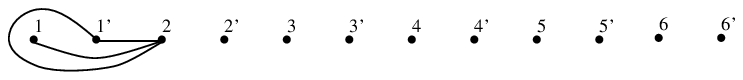}}}$  \> corresponding to $Z^3_{11',2}$ \\ [.5cm]
(3)  \  $\vcenter{\hbox{\epsfbox{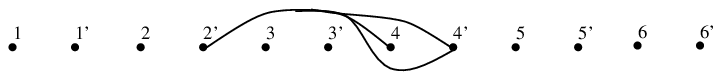}}}$  \> corresponding to $\bZ^3_{2',44'}$ \\ [.5cm]
(4) \   $\vcenter{\hbox{\epsfbox{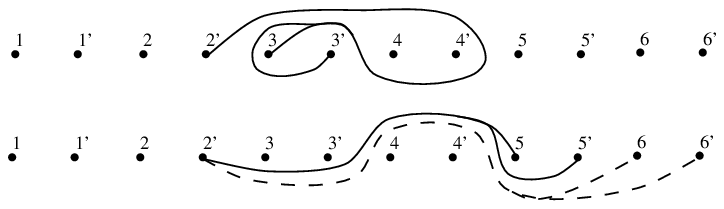}}}$  \> corresponding to
 $(\uZ^2_{2'i})^{\uZ^{-2}_{i,44'}}_{i=33',55',6,6'}$
 \\ [.5cm]
(5)  \  $\vcenter{\hbox{\epsfbox{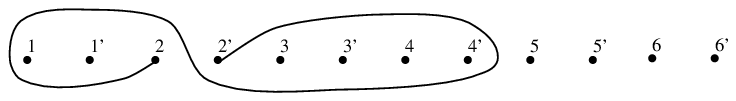}}}$  \> corresponding to $(Z_{22'})^{Z^2_{11',2}\uZut_{2',44'}Z^2_{2',33'}}$ \\ [.5cm]
(6)  \  $\vcenter{\hbox{\epsfbox{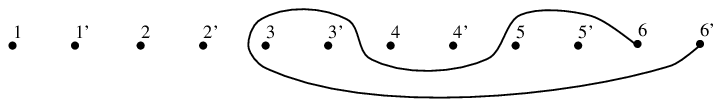}}}$  \> corresponding to $(Z_{66'})^{\uZumt_{33',6}  Z^{-2}_{55',6}}$ \\ [.5cm]
(7)   \ $\vcenter{\hbox{\epsfbox{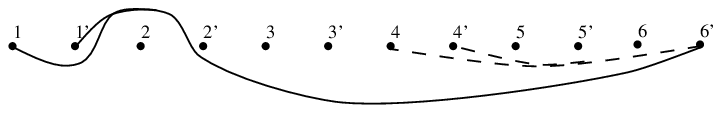}}}$  \> corresponding to $(\uZ^2_{i6'})^{Z^2_{11',2} }_{i = 11',44'}$\\ [.5cm]
(8)   $\vcenter{\hbox{\epsfbox{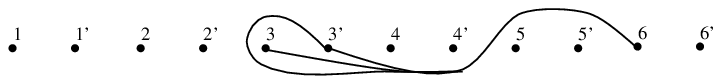}}}$  \>
 corresponding to $(\uZ^3_{33',6})^{Z^{-2}_{55',6}}$\\ [.5cm]
(9)   $\vcenter{\hbox{\epsfbox{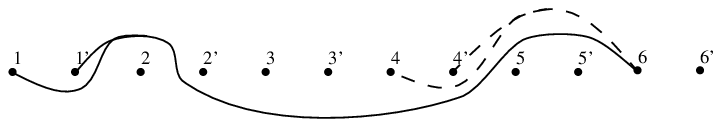}}}$  \> corresponding to $(\uZ^2_{i6})^{Z^{-2}_{55',6}Z^2_{11',2}}_{i = 11',44'}$\\ [.5cm]
(10)   $\vcenter{\hbox{\epsfbox{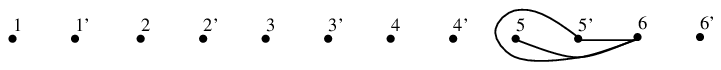}}}$  \> corresponding to $Z^3_{55',6}$
\end{tabbing}
}
\noindent
\underline{Proof}: A similar proof as the proof of Theorem 8.6 but changes are applied on (2)$_4$, to obtain a formula for a local braid monodromy of
$S^{(0)}_4$.  Moreover, all changes were applied on indices 3,4 to obtain
$H_{V_4}$.
\hfill
$\Box$

\noindent
{\bf Invariance Property 9.7.}  {\em $H_{V_4}$ is invariant under
$(Z_{22'} Z_{66'})^p (Z_{11'}Z_{55'})^q (Z_{33'} Z_{44'})^r \\
\forall p,q,r\in \Z$.}

\noindent
\underline{Proof}: 
\underline{Case 1}: $p = q = r$.\\
A similar proof as Case 1 in Invariance Property 8.7. \\
\underline{Case 2}: $p = 0$.\\
Denote: $\epsilon = (Z_{11'} Z_{55'})^q (Z_{33'}Z_{44'})^r$.
We want to prove that  $H_{V_4}$ is invariant under $\epsilon $.

\noindent
\underline{Step 1:}  Product of the factors outside of $\hF_1(\hF_1)^{\rho^{-1}}$.
\\
$(Z_{22'})^{Z^2_{11',2}\uZut_{2',44'}Z^2_{2',33'}}$
 and
$(Z_{66'})^{\uZumt_{33',6}  Z^{-2}_{55',6}} $
commute with $\epsilon $. \\
The degree 3 factors are invariant under $\epsilon $ by Invariance Rule III.
All conjugations are invariant under $\epsilon$ too.
The
degree 2 factors are of the forms $Z^2_{\alpha,\beta \beta ' }$ where $\alpha
= 2'$, $Z^2_{\alpha  \alpha ', \beta }$ where $\beta = 6,6'$.  So by Invariance Rule II, they are invariant under
$\epsilon $.

\noindent
\underline{Step 2:}   $\hF_1(\hF_1)^{\rho^{-1}}$.\\
A similar proof as step 2 in Invariance Property 8.7, but for:\\
Case 2.1:  $q = 0 \ ; \ \epsilon = (Z_{33'}Z_{44'})^r$.\\
Case 2.2:  $r = 0 \ , \ q=1; \epsilon = Z_{11'} Z_{55'}$.\\
Case 2.3:  $q = 2q' \ ; \ \epsilon = (Z_{11'} Z_{55'})^{2q'}(Z_{33'}
Z_{44'})^r$.\\
Case 2.4:  $q = 2q'+1 \ ; \ \epsilon = (Z_{11'} Z_{55'})^{2q'+1}(Z_{33'}Z_{44'})^r$.\\
\noindent
\underline{Case 3}:  $p,q,r$ arbitrary ; \ $\epsilon = (Z_{22'}Z_{66'})^p(Z_{11'}Z_{55'})^q(Z_{33'}Z_{44'})^r$.\\
By Case 1: $H_{V_4}$ is invariant under
$(Z_{22'}Z_{66'})^{p} (Z_{11'}Z_{55'})^{p}
(Z_{33'}Z_{44'})^{p}$.
By Case 2: $H_{V_4}$ is invariant under
$(Z_{11'}Z_{55'})^{q-p} (Z_{33'}Z_{44'})^{r-p}$.
By Invariance Remark (v),  $H_{V_4}$ is invariant  under \nolinebreak
$\epsilon \ . \Box$

\noindent
{\bf Theorem 9.8.}  {\em Consider $H_{V_4}$ from Theorem 9.6.  The paths corresponding to the factors in $H_{V_4}$ , are shown below considering the Invariance Property 9.7.  Moreover, below to (3),(4),(9) their complex conjugates appear too.}

\noindent
\underline{Remark}:  By abuse of notation, the simple braids are denoted by $
\underline{z}_{mk}$ and the more complicated paths are denoted by $\tilde{z}_{mk}$.

\noindent
\underline{Proof:}  By  Invariance Property 9.7,  $H_{V_4}$ is invariant under $\epsilon =(\rho_2 \rho_6)^p (\rho_1 \rho_5)^q (\rho_3 \rho_4)^r \quad\\
 \forall p, q, r \in \Z$.
There are two possible applications for presenting the braids: (a) if $m,k$ are two indices in two different parts of $\epsilon $, then the braid is
$\rho^j_k \rho^i_m z_{mk} \rho^{-i}_m \rho^{-j}_k, \ i, j \in \Z, \ i \neq j$.
(b) If $m,k$ are in the same part of $\epsilon $, then the braid is $\rho^i_k \rho^i_m z_{mk} \rho^{-i}_m \rho^{-i}_k$, i.e.: when conjugating $z_{mk}$ by $\rho^i_m$ for some $i$, we must conjugate $z_{mk}$ also by $\rho^i_k$ for the same $i$.

We apply Complex Conjugation  on (3), (4), (9).  Their
complex conjugates appear below them.

\begin{center}
\begin{tabular}[hb]{|l|p{3in}|l|c|}\hline
& The paths corresponding & The braids & The exponent  \\
[-.2cm] & to the braids and their  & &  (according to  \\
[-.2cm] & complex conjugates & & singularity type)\\  \hline (1) &
$\vcenter{\hbox{\epsfbox{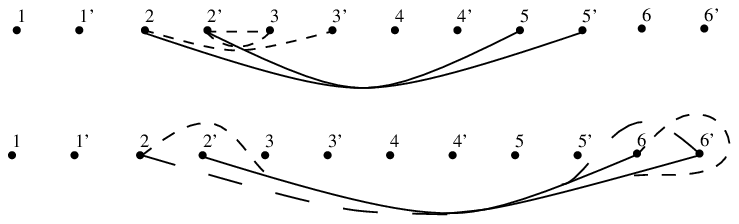}}}$
 & $\begin{array}{lll}\rho^j_m \rho^i_2 \underline{z}_{2'm} \rho^{-i}_2 \rho^{-j}_m\\ [-.2cm]
m = 3,5 \\ [-.2cm]
\rho^i_6 \rho^i_2 \underline{z}_{2'6} \rho^{-i}_2 \rho^{-i}_6
\end{array}$ & 2 \\
    \hline
(2) &  $\vcenter{\hbox{\epsfbox{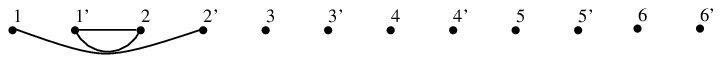}}}$ & $\rho^j_2 \rho^i_1 \underline{z}_{12} \rho^{-i}_1 \rho^{-j}_2$ & 3  \\
   \hline
(3) &  $\vcenter{\hbox{\epsfbox{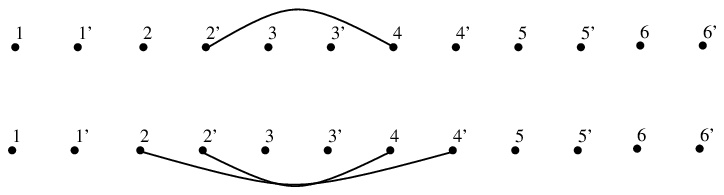}}}$ & $ \rho^j_4 \rho^i_2 \underline{z}_{2'4} \rho^{-i}_2 \rho^{-j}_4$ & 3
\\
   \hline
(4) &  $\vcenter{\hbox{\epsfbox{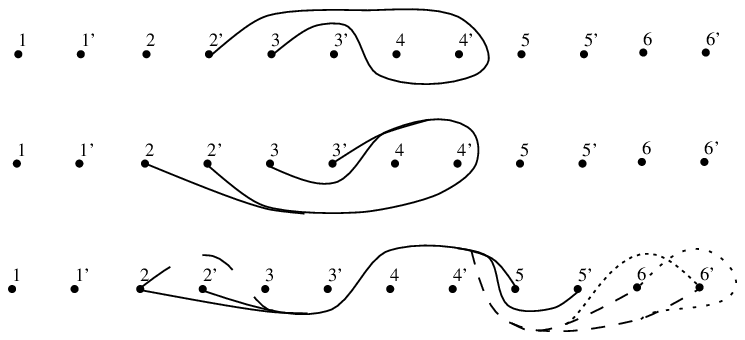}}}$ & $\begin{array}{lll}
\rho^j_m \rho^i_2 \tilde{z}_{2'm} \rho^{-i}_2 \rho^{-j}_m \\ [-.2cm]
m = 3,5\\ [-.2cm]
\rho^i_6 \rho^i_2 \tilde{z}_{2'6} \rho^{-i}_2 \rho^{-i}_6 \end{array}$ & 2 \\
    \hline
(5) & $\vcenter{\hbox{\epsfbox{v4076.ps}}}$ & $\tilde{z}_{22'}$ & 1\\
   \hline
(6) & $\vcenter{\hbox{\epsfbox{v4077.ps}}}$ & $\tilde{z}_{66'}$ & 1  \\
   \hline
(7) & $\vcenter{\hbox{\epsfbox{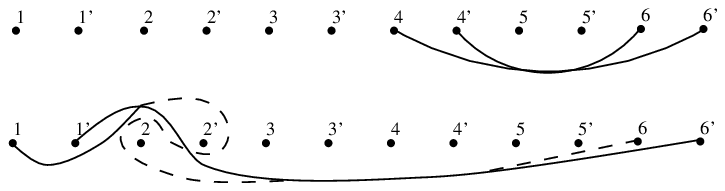}}}$ & $\begin{array}{ll}\rho^j_6 \rho^i_m \tilde{z}_{m6'}\rho^{-i}_m \rho^{-j}_6 \\ [-.2cm]
m = 1,4\end{array}$ & 2\\
   \hline
(8) & $\vcenter{\hbox{\epsfbox{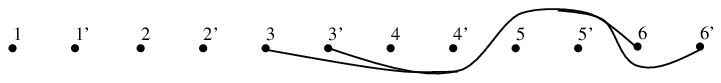}}}$ & $\rho^j_6 \rho^i_3 \tilde{z}_{36} \rho^{-i}_3 \rho^{-j}_6$ & 3 \\
    \hline
(9) & $\vcenter{\hbox{\epsfbox{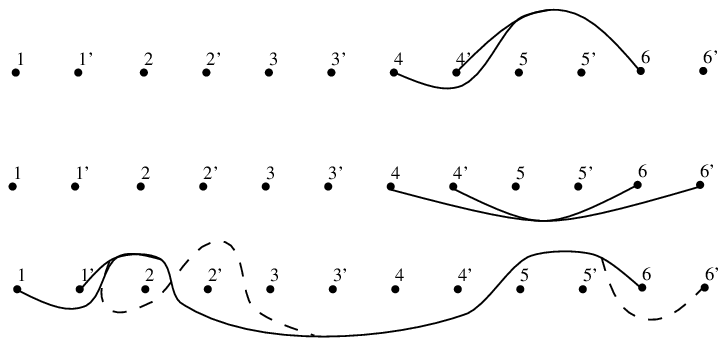}}}$ &
$\begin{array}{ll}\rho^j_6 \rho^i_m \tilde{z}_{m6} \rho^{-i}_m \rho^{-j}_6
\\ [-.2cm]
m = 1,4 \end{array}$ & 2 \\
   \hline
(10) & $\vcenter{\hbox{\epsfbox{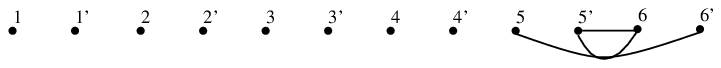}}}$ &
$\rho^j_6 \rho^i_5 z_{56} \rho^{-i}_5 \rho^{-j}_6$ & 3 \\
   \hline
\end{tabular}
\end{center}

Consider now $\hF_1 (\hF_1)^{\rho^{-1}} (\rho^{-1} = \rho^{-1}_5 \rho_1^{-1})$ in $H_{V_4}$.
In a similar proof of Theorem 8.8 we obtain the following table:

\begin{center}
\begin{tabular}[H]{|l|c|c|} \hline
 $\vcenter{\hbox{\epsfbox{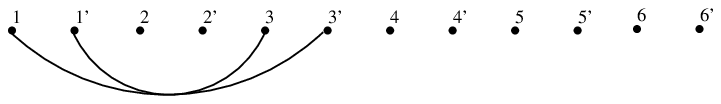}}}$ & $\rho^j_3 \rho^i_1 \underline{z}_{1'3} \rho^{-i}_1 \rho^{-j}_3$ & 3\\ [.3cm]
    \hline
 $\vcenter{\hbox{\epsfbox{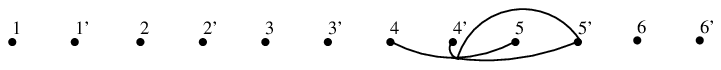}}}$  & $\rho^j_5 \rho^i_4  \underline{z}_{45} \rho^{-i}_4 \rho^{-j}_5$ & 3\\
    \hline
 $\vcenter{\hbox{\epsfbox{4v079.ps}}}$ & $\alpha _1(m_1)$ & 1\\
    \hline
 $\vcenter{\hbox{\epsfbox{4v080.ps}}}$  & $\alpha _2(m_1)$ & 1\\
    \hline
 $\vcenter{\hbox{\epsfbox{4v081.ps}}}$  & $(\alpha _1(m_1))^{\rho^{-1}}$ & 1\\
    \hline
 $\vcenter{\hbox{\epsfbox{4v082.ps}}}$   & $(\alpha _2(m_1))^{\rho^{-1}}$ & 1\\   \hline
 $\vcenter{\hbox{\epsfbox{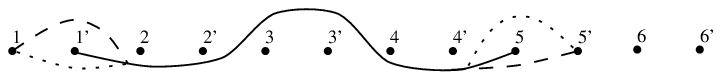}}}$ & $\rho^i_5 \rho^i_1 \tilde{z}_{1'5} \rho^{-i}_1 \rho^{-i}_5$ & 2\\
    \hline
 $\vcenter{\hbox{\epsfbox{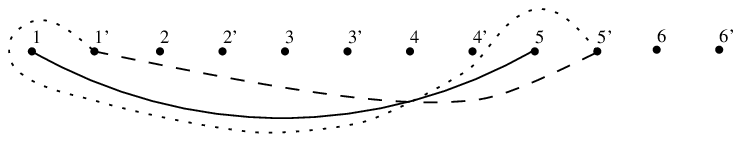}}}$  & $\rho^i_5 \rho^i_1 \underline{z}_{15} \rho^{-i}_1 \rho^{-i}_5$ & 2\\
   \hline
\end{tabular}
\end{center}
\hfill
$\Box$

\section{The computation of $H_{V_7}$}\label{sec:37}

\indent
We follow Figure 6 when computing the local braid monodromy 
for $V_7$, denoted by
$\vp^{(0)}_7$, but first we need to compute its factors which come from singularities in the neighbourhood of $V_7$.
We start with $\vp^{(6)}_7$.
Then we compute $S^{(5)}_7$ and $\vp^{(5)}_7$ in the neighbourhood of $V_7$.
Notice that $\vp^{(5)}_7 = \vp^{(4)}_7  = \vp^{(3)}_7= \vp^{(2)}_7 =\vp^{(1)}_7$ since the fourth, fifth, sixth, seventh regenerations act similarly on the
 the neighbourhood of $V_7$.  Finally we compute $\vp^{(0)}_7$.
Hence we get a factorized expression for $\vp^{(0)}_7$.

\noindent
{\bf Theorem 10.1.}  {\em In a neighbourhood of $V_7$, the local braid monodromy of $S^{(6)}_7$ around $V_7$ is given by:\\
$\vp^{(6)}_7 = (\uZ^2_{2'i})_{i=3,4,4',6} \cdot
Z^4_{12} \cdot
(\uZ^4_{2'5})^{Z^2_{2'3}} \cdot
(\uZ^2_{2i})_{i=3,4,4',6}^{Z^{2}_{12}} \cdot
(Z_{22'})^{Z^2_{12}\bZ^2_{2'5}Z^2_{44',5}}
\cdot
Z^4_{34} \cdot
(\uZ^2_{i4'})^{Z^2_{12}}_{i = 1,5} \cdot
\bZ^4_{4'6}  \cdot
(\uZ^2_{i4})^{Z^{2}_{34}Z^2_{12}}_{i = 1,5} \cdot
(Z_{44'})^{\bZ^2_{4'6}Z^2_{34}} \cdot
(\Delta^2<1,3,5,6>)^{\uZ^{-2}_{45}\uZ^{-2}_{46}Z^2_{12}}$,
  where $(Z_{22'})^{Z^2_{12}\bZ^2_{2'5}Z^2_{44',5}}$
is determined by  $\vcenter{\hbox{\epsfbox{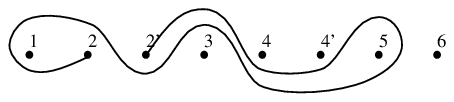}}}$ ,
$(Z_{44'})^{\bZut_{4'6} Z^{2}_{34}}$
is determined
by  $\vcenter{\hbox{\epsfbox{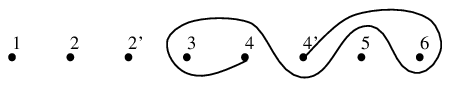}}}$
and the}  {\rm L.V.C.} {\em corresponding to
$\Delta^2<1,3,5,6>$ is given by}
$\vcenter{\hbox{\epsfbox{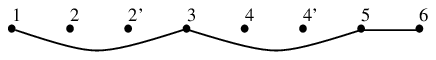}}}$ .

\noindent
\underline{Proof}: \ Let $\hat{L}_i, 1 \le i \le 6$,
 be the lines in $Z^{(8)}$ and in $Z^{(7)}$ intersecting in a point 7.
  Let $L_i = \pi^{(8)} (\hat{L}_i)$.  Then $L_i, 1 \le i \le 6$
, intersect in $V_7$.  $S^{(8)}$ (and $S^{(7)})$ around $V_7$ is
 $\bigcup\limits^6_{i=1} L_i$.  After the second regeneration, $S^{(6)}_7$ in a neighbourhood of $V_7$ is as follows:  the lines $L_2$ and $L_4$ are replaced by two conics $Q_2$ and $Q_4$, where $Q_2$ is tangent to the lines $L_1$ and $L_5$ and $Q_4$ is tangent to $L_3$ and $L_6$ .
 This follows from the fact that $P_7 \cup P_{11}$ (resp. $P_{8} \cup P_{16})$ is regenerated to $R_4$ (resp. $R_7)$, and the tangency follows from Lemma \ref{lm:24}.  Thus, in a neighbourhood of $V_7$, $S^{(6)}_7$ is of the form: $S^{(6)}_7 =
L_1 \cup Q_2 \cup L_3 \cup Q_4 \cup L_5 \cup L_6$.

$S^{(6)}_7$ is shown in Figure 24.
\begin{figure}[h]
\epsfxsize=11cm 
\epsfysize=9cm 
\begin{minipage}{\textwidth}
\begin{center}
\epsfbox {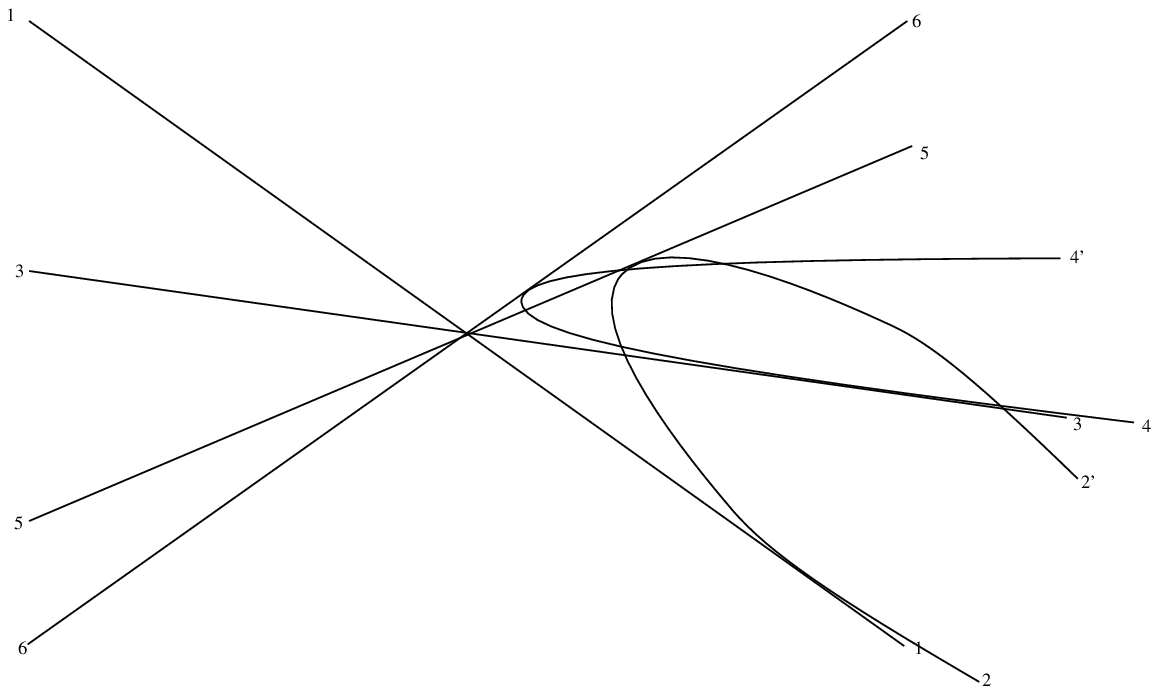}
\end{center}
\end{minipage}
\caption{}
\end{figure}

We divide the proof into parts:\\
\noindent
{\bf Proposition 10.1.1.} {\em Consider in $\C^2$ the following figure:  $C_1 = L_1 \cup Q_2 \cup L_3 \cup L_5$, where $L_i$ are lines, $i = 1,3,5, \ Q_2$ is a conic tangent to $L_1$ and $L_5$, it intersects $L_3$ and the 3 lines meet at one point (see Figure 25).

Then the braid monodromy of $C_1$ is given by:\\
$\vp_{C_1} = Z^2_{2'3}
\cdot Z^4_{12} \cdot (\uZ^4_{2'5})^{Z^2_{2'3}} \cdot (\uZ^2_{23})^{Z^2_{12}} \cdot
(Z_{22'})^{Z^2_{12}\bZ^2_{2'5}} \cdot
(\Delta^2<1,3,5>)^{Z^2_{12}}$,
 where
$(Z_{22'})^{Z^2_{12}\bZ^2_{2'5}}$
is determined by
$\vcenter{\hbox{\epsfbox{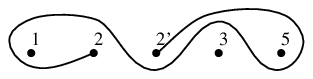}}}$ and the L.V.C. corresponding to $\Delta^2<1,3,5>$ is given by 
$\vcenter{\hbox{\epsfbox{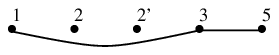}}}$}.

\begin{figure}[!h]
\epsfxsize=11cm 
\epsfysize=9cm 
\begin{minipage}{\textwidth}
\begin{center}
\epsfbox {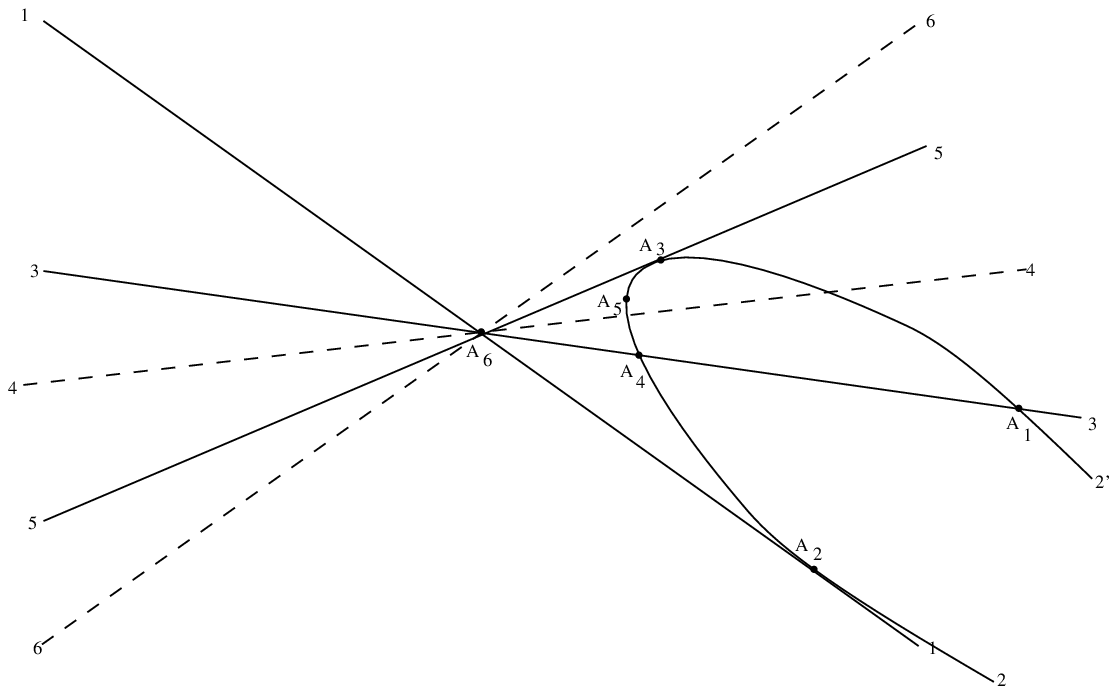}
\end{center}
\end{minipage}
\caption{}
\end{figure}

\underline{Proof}:
Let $\pi_1 : E \times D \rightarrow E$ be the projection to $E$.

Let $\{A_j\}^{6}_{j=1}$ be  singular points of $\pi_1$ as follows:\\
$A_1, A_4 $ are  intersection points of $Q_2$ with the line $L_3$.\\
$A_2, A_3 $ are tangent  points of $Q_2$ with the lines $L_1, L_5$ respectively.\\
$A_{6}$ is an intersection  point of the  lines $L_1, L_3 ,L_5$.\\
$A_5$ is a  point of the type $a_1$ in $Q_2$.

Let $N,M, \{\ell(\g_j)\}^{6}_{j=1}$
as in Proposition 8.1.1.
Let $i = L_i \cap K$ for $i = 1,3,5$.
 Let $\{2,2'\} = Q_2 \cap K$.
 $K = \{1,2,2',3,5\}$, \st 1,2,2',3,5 are  real points.
Moreover: $1 < 2 < 2'< 3 < 5 $.
Let $\beta_M$ be the diffeomorphism defined by: $\beta_M(i) = i$ \ for \ $i = 1,2,5,  \ \beta_M(2') = 3, \ \beta_M (3) = 4$.  Moreover: deg$C_1 = 5$ and $\# K_\R (x) \geq 3 \ \ \forall x$.

We are looking for $\vp_M(\ell (\g_j))$ for $j = 1, \cdots , 6$.  We choose a $g$-base $\{\ell(\g_j)\}^{6}_{j=1}$ of $\pi_1(E - N,u)$, \st each path $\g_j$ is below the real line.

Put all data in the following table:
\vspace{-1cm}
\begin{center}
\begin{tabular}[H]{cccc} \\
$j$ & $\lambda_{x_j}$ & $\epsilon_{x_j}
$ & $\delta_{x_j}$ \\ \hline
1 & $<3,4>$ & 2 & $\Delta<3,4>$\\
2 & $<1,2>$ & 4 & $\Delta^2<1,2>$\\
3 & $<4,5>$ & 4 & $\Delta^2<4,5>$\\
4 & $<2,3>$ & 2 & $\Delta<2,3>$\\
5 & $<3,4>$ & 1 & $\Delta^{\frac{1}{2}}_{I_2 \R}<3>$\\
6 & $<1,3>$ & 2 & $\Delta<1,3>$\\
\end{tabular}\end{center}

For computations, we use the formulas in [Am2, Theorems 1.41, 1.44].

\noindent
\underline{Remark}:  $\beta_{x'_j} (K(x'_j)) = \{1,2,3,4,5\}$ for $1 \leq j \leq 5$.\\
$\beta_{x'_6} (K(x'_6)) = \{1,2,3,4+i, 4-i \}$.\\
\\
\noindent
$(\xi_{x'_1}) \Psi_{\g'_1} = <3,4> \beta ^{-1}_M  = z_{2'3}$

\noindent
$<3,4>$ \ \ $\vcenter{\hbox{\epsfbox{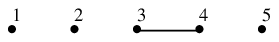}}}$ \ \ $\beta ^{-1}_M$
 \ $\vcenter{\hbox{\epsfbox{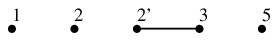}}}$ \ \
$\vp_M(\ell(\g_1)) = Z^2_{2'3}$ \\

\noindent
$(\xi_{x'_2}) \Psi_{\g'_2} = <1,2> \Delta<3,4>  \beta ^{-1}_M
= z_{12}$\\
\noindent
$<1,2>  \Delta<3,4>$ \ $\vcenter{\hbox{\epsfbox{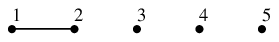}}}$
\ $\beta ^{-1}_M$ \ $\vcenter{\hbox{\epsfbox{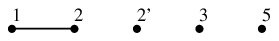}}}$  \ $ \vp_M(\ell(\g_2)) = Z^4_{12}$\\

\noindent
$(\xi_{x'_3}) \Psi_{\g'_3} =
<4,5> \Delta^2<1,2> \Delta<3,4> \beta ^{-1}_M  =
\underline{z}_{2'5}^{Z^2_{2'3}} $

\noindent
$<4,5> \Delta^2<1,2>$   \ $\vcenter{\hbox{\epsfbox{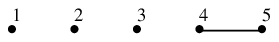}}}$ \ $ \Delta<3,4>$ \ $\vcenter{\hbox{\epsfbox{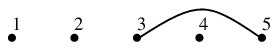}}}$\ \
$\beta ^{-1}_M$  \ $\vcenter{\hbox{\epsfbox{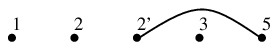}}}$  \ \
$\vp_M(\ell(\g_3)) = (\uZ^4_{2'5})^{Z^2_{2'3}} $\\

\noindent
$(\xi_{x'_4}) \Psi_{\g'_4} = <2,3> \Delta^2<4,5> \Delta^2<1,2> \Delta<3,4> \beta ^{-1}_M   =\underline{z}_{23}^{Z^2_{12}}  $

\noindent
$<2,3> \Delta^2<4,5>$  $\vcenter{\hbox{\epsfbox{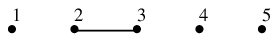}}}$  \ \
$\Delta^2<1,2>$ $\vcenter{\hbox{\epsfbox{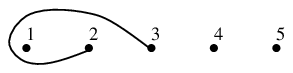}}}$  \\
$\Delta<3,4>$ $\vcenter{\hbox{\epsfbox{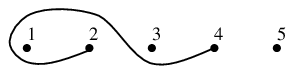}}}$  \ \
$\beta^{-1}_M $
 $\vcenter{\hbox{\epsfbox{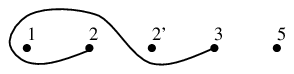}}}$ \ \
$ \vp_M(\ell(\g_4)) = (\uZ^2_{23})^{Z^2_{12}}$\\

\noindent
$(\xi_{x'_5}) \Psi_{\g'_5}   = <3,4> \Delta<2,3> \Delta^2<4,5> \Delta^2<1,2>
\Delta<3,4>  \beta ^{-1}_M = z_{22'}^{Z^{2}_{12}\bZ^2_{2'5}}$

\noindent
$<3,4> $  $\vcenter{\hbox{\epsfbox{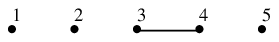}}}$  \   $\Delta<2,3>$\
\ $\vcenter{\hbox{\epsfbox{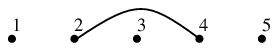}}}$  \ \
$\Delta^2<4,5>$   $\vcenter{\hbox{\epsfbox{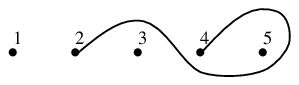}}}$ \\
$\Delta^2<1,2>$ \
$\vcenter{\hbox{\epsfbox{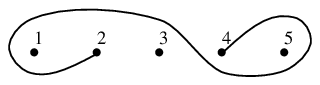}}}$  \ \
$\Delta<3,4>$ $\vcenter{\hbox{\epsfbox{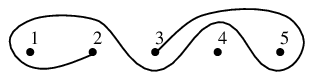}}}$ \
$\beta ^{-1}_M $ \ $\vcenter{\hbox{\epsfbox{v7019.ps}}}$ \
$\vp_M(\ell(\g_5)) = (Z_{22'})^{Z^{2}_{12}\bZ^2_{2'5}} $\\

\noindent
$(\xi_{x'_6}) \Psi_{\g'_6}= <1,3> \Delta^{\frac{1}{2}}_{I_2 \R}<3>
\Delta<2,3> \Delta^2<4,5> \Delta^2<1,2>\Delta<3,4>\beta ^{-1}_M = \\
(\Delta<1,3,5>)^{
Z^2_{12}}$

\noindent
$<1,3>$ $\vcenter{\hbox{\epsfbox{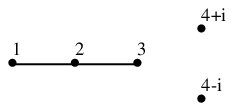}}}$ \  $\Delta^{\frac{1}{2}}_{I_2 \R}<3> $ \ $\vcenter{\hbox{\epsfbox{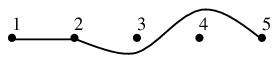}}}$ \ \
$\Delta<2,3>$ $\vcenter{\hbox{\epsfbox{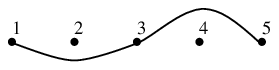}}}$ \\
$\Delta^2<4,5> $ \ $\vcenter{\hbox{\epsfbox{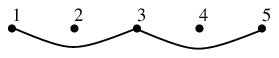}}}$ \ \
$\Delta^2<1,2>$ $\vcenter{\hbox{\epsfbox{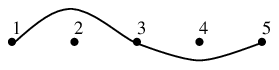}}}$ \ \
 $\Delta<3,4>$ \
$\vcenter{\hbox{\epsfbox{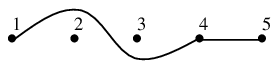}}}$ \\
$\beta ^{-1}_M $
\ $\vcenter{\hbox{\epsfbox{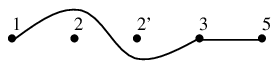}}}$\ \
$\vp_M(\ell(\g_6)) =
(\Delta^2<1,3,5>)^{Z^2_{12}}$\\

\noindent
{\bf Remark 10.1.2.} {\em In a similar way as in  Remark 8.1.2, the changes are:\\
(a) $(\Delta^2<1,3,5>)^{Z^2_{12}}$ is replaced by
$(\Delta^2<1,3,4,5,6>)^{Z^2_{12}}$ with the corresponding {\rm L.V.C.}
which appears as follows:}
{\em $\vcenter{\hbox{\epsfbox{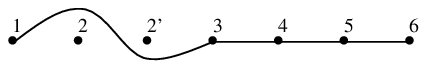}}}$ .\\
(b) $Z^2_{2'3}$ is replaced by $(\uZ^2_{2'i})_
{i = 3,4,6}$ as follows:
 $\vcenter{\hbox{\epsfbox{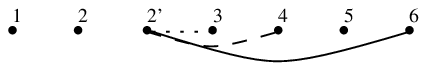}}}$ .\\
(c) $(\uZ^2_{23})^{Z^2_{12}}$ is replaced by
$(\uZ^2_{2i})_{i = 3,4,6}^{Z^2_{12}}$ as follows:
 $\vcenter{\hbox{\epsfbox{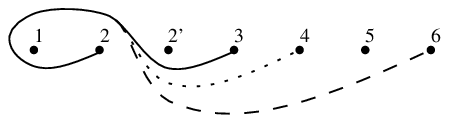}}}$ .}

\noindent
{\bf Proposition 10.1.3.} {\em Consider in $\C^2$ the following figure: $C_2 = L_3 \cup Q_4 \cup L_5 \cup  L_6$, where $L_i$ are lines, $i = 3,5,6, \ Q_4$ is a conic tangent to $L_3$ and $L_6$, it intersects $L_5$  and the 3 lines meet at one point (see Figure 26).

Then the braid monodromy of $C_2$ is given by:\\
$\vp_{C_2} = Z^4_{34} \cdot Z^2_{4'5} \cdot \bZ^4_{4'6} \cdot
(\uZ^2_{45})^{Z^{2}_{34}} \cdot
(Z_{44'})^{\bZ^2_{4'6}Z^2_{34}} \cdot  (\Delta^2<3,5,6>)^{\uZ^{-2}_{45}\uZ^{-2}_{46}}$,
 where
$(Z_{44'})^{\bZ^{2}_{4'6}Z^2_{34}}$  is determined by
$\vcenter{\hbox{\epsfbox{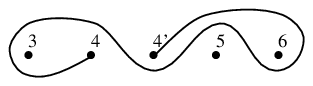}}}$
  and the {\rm   L.V.C.} corresponding to $\Delta^2<3,5,6>$ is given by
$\vcenter{\hbox{\epsfbox{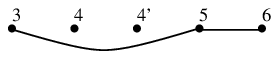}}}$.}

\begin{figure}[h]
\epsfxsize=11cm 
\epsfysize=9cm 
\begin{minipage}{\textwidth}
\begin{center}
\epsfbox {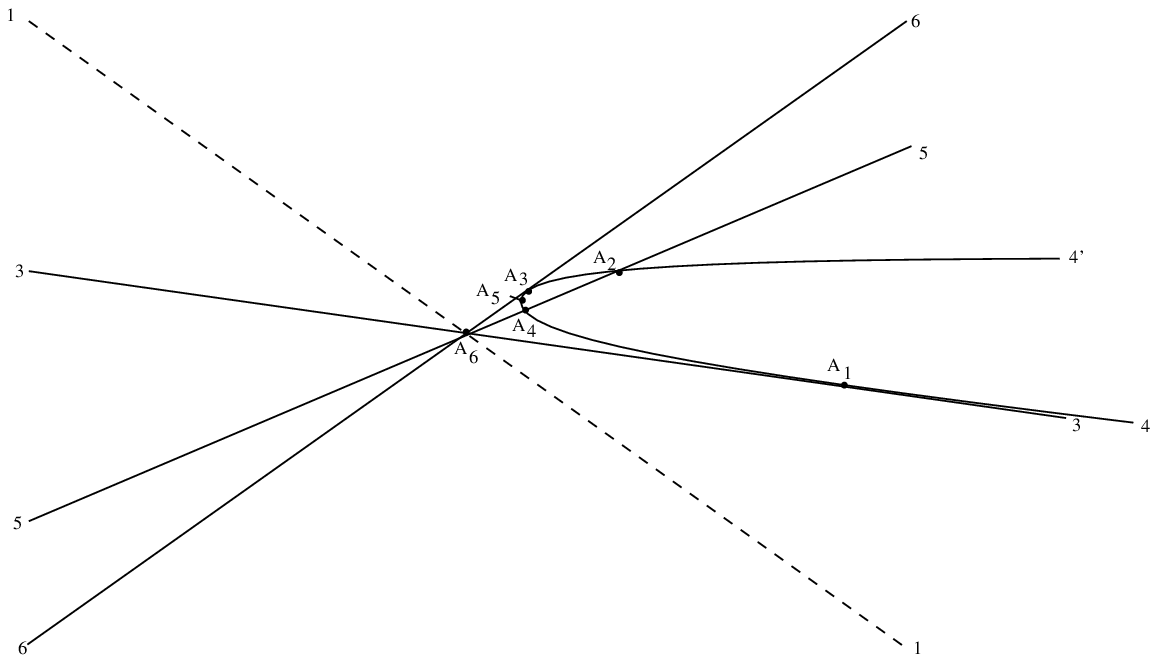}
\end{center}
\end{minipage}
\caption{}
\end{figure}

\underline{Proof:} Let $\pi_1: E \times D \rightarrow E$ be the projection to $E$.

Let $\{A_j\}^{6}_{j=1}$ be singular points of $\pi_1$ as follows:\\
$A_1, A_3$ are  tangent points of $Q_4$ with the lines $L_3, L_6$ respectively.\\
$A_2, A_4$  are intersection points of $Q_4$ with the line $L_5$. \\
$A_6$ is an intersection point of the lines $L_3, L_5, L_6$.\\
$A_5$ is a point of the type $a_1$ in $Q_4$.

Let $K,N,M, \{\ell (\g_j)\}^{6}_{j=1}$ be as in Proposition 8.1.3.

Let $i = L_i \cap K$ for $i = 3,5,6, \ \{4,4'\} = Q_4 \cap K$.
$K = \{3,4,4',5,6\}$, \st $3,4,4',5,6$ are real points and
$3 < 4 < 4'< 5 < 6$.
Let $\beta_M$ be the diffeomorphism defined by
$\beta_M(3) = 1 , \   \beta_M(4) = 2, \ \beta_M(4') = 3, \  \beta_M(5) = 4, \  \beta_M(6) = 5$.  deg$C_2 = 5$ and
 $\# K_\R (x) \geq 3 \ \forall x$.

Put all data in the following table:
\begin{center}
\begin{tabular}[!H]{cccc} \\
$j$ & $\lambda_{x_j}$ & $\epsilon_{x_j}
$ & $\delta_{x_j}$ \\ \hline
1 & $<1,2>$ & 4 & $\Delta^2<1,2>$\\
2 & $<3,4>$ & 2 & $\Delta<3,4>$\\
3 & $<4,5>$ & 4 & $\Delta^2<4,5>$\\
4 & $<2,3>$ & 2 & $\Delta<2,3>$\\
5 & $<3,4>$ & 1 & $\Delta^{\frac{1}{2}}_{I_2 \R}<3>$\\
6 & $<1,3>$ & 2 & $\Delta<1,3>$
\end{tabular}
\end{center}

\noindent
\underline{Remark}:  $\beta_{x'_j} (K(x'_j)) = \{1,2,3,4,5\}$ for $1 \leq j \leq 5$.\\
$\beta_{x'_6}(K(x'_6)) = \{1,2,3,4+i, 4-i\}$.\\

\medskip
\noindent
$(\xi_{x'_1}) \Psi_{\g'_1} = <1,2> \beta^{-1}_M = z_{34}$

\noindent
$<1,2>$ $\vcenter{\hbox{\epsfbox{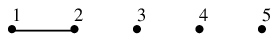}}}$ \ \
$\beta^{-1}_M$ $\vcenter{\hbox{\epsfbox{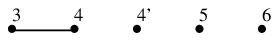}}}$ \ \
$\vp_M(\ell(\g_1)) = Z^4_{34}$\\

\noindent
$(\xi_{x'_2}) \Psi_{\g'_2} = <3,4> \Delta^2<1,2> \beta^{-1}_M = z_{4'5}$

\noindent

\noindent
$<3,4> \Delta^2<1,2>  \vcenter{\hbox{\epsfbox{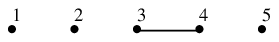}}}$ \
$\beta^{-1}_M$ \ \ $\vcenter{\hbox{\epsfbox{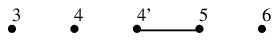}}}$\ \
$\vp_M(\ell(\g_2)) = Z^2_{4'5}$\\

\noindent
$(\xi_{x'_3}) \Psi_{\g'_3} = <4,5> \Delta<3,4> \Delta^2<1,2>
  \beta^{-1}_M =
\bar{z}_{{4'6}}$

\medskip

\noindent
$<4,5>  \vcenter{\hbox{\epsfbox{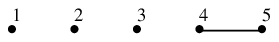}}}$ \ \
$\Delta<3,4> \Delta^2<1,2> $ $\vcenter{\hbox{\epsfbox{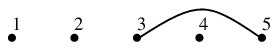}}}$ \ \
$\beta^{-1}_M$ $\vcenter{\hbox{\epsfbox{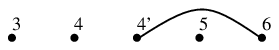}}}$ \ \
$\vp_M(\ell(\g_3)) = \bZ^4_{4'6}$\\

\noindent
$(\xi_{x'_4}) \Psi_{\g'_4}  =  <2,3> \Delta^2<4,5>
 \Delta<3,4> \Delta^2<1,2> \beta^{-1}_M = \underline{ z}_{45}^{Z^2_{34}}$

\noindent
$<2,3>\Delta^2<4,5>  \vcenter{\hbox{\epsfbox{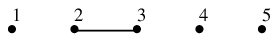}}}$ \ \
$\Delta<3,4>$
 $\vcenter{\hbox{\epsfbox{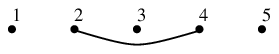}}}$ \\
$\Delta^2<1,2>$
 $\vcenter{\hbox{\epsfbox{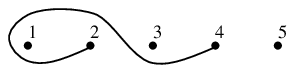}}}$ \ \
$\beta^{-1}_M$ $\vcenter{\hbox{\epsfbox{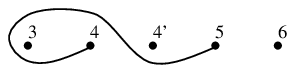}}}$ \ \
$\vp_M(\ell(\g_4)) = (\uZ^2_{45})^{Z^2_{34}}$\\

\noindent
$(\xi_{x'_5}) \Psi_{\g'_5}  =  <3,4>  \Delta<2,3> \Delta^2<4,5>
 \Delta<3,4> \Delta^2<1,2> \beta^{-1}_M = z_{44'}^{\bZ^2_{4'6}Z^2_{34}}$

\noindent
$<3,4>$  $\vcenter{\hbox{\epsfbox{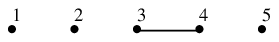}}}$ \ \
$\Delta<2,3>$  $\vcenter{\hbox{\epsfbox{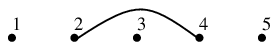}}}$ \ \
$\Delta^2<4,5>$
 $\vcenter{\hbox{\epsfbox{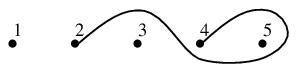}}}$ \\
\noindent
$\Delta<3,4>$
 $\vcenter{\hbox{\epsfbox{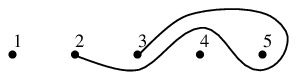}}}$ \ \
$\Delta^2<1,2>$ \
 $\vcenter{\hbox{\epsfbox{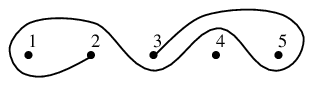}}}$ \ \
$\beta^{-1}_M$ $\vcenter{\hbox{\epsfbox{v7047.ps}}}$ \ \
$\vp_M(\ell(\g_5)) = (Z_{44'})^{\bZ^2_{4'6}Z^2_{34}}$\\

\medskip
\noindent
$(\xi_{x'_6}) \Psi_{\g'_6}  =  <1,3>   \Delta^{\frac{1}{2}}_{I_2 \R}<3>
\Delta<2,3>\Delta^2<4,5>  \Delta<3,4>\Delta^2<1,2>\beta^{-1}_M =
(\Delta<3,5,6>)^{\uZ^{-2}_{45}\uZ^{-2}_{46}}$

\noindent
$ <1,3>  \vcenter{\hbox{\epsfbox{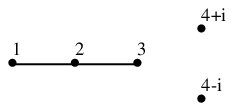}}}$ \ \
$\Delta^{\frac{1}{2}}_{I_2 \R}<3>
\vcenter{\hbox{\epsfbox{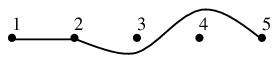}}}$ \ \
$ \Delta<2,3> \ \vcenter{\hbox{\epsfbox{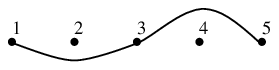}}}$ \\
$\Delta^2<4,5>$
 $\vcenter{\hbox{\epsfbox{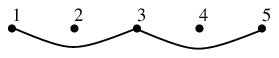}}}$ \ \
$\Delta<3,4>$
 $\vcenter{\hbox{\epsfbox{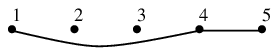}}}$ \ \
$\Delta^2<1,2>$
 $\vcenter{\hbox{\epsfbox{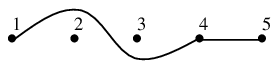}}}$ \\
$\beta^{-1}_M$ $\vcenter{\hbox{\epsfbox{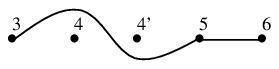}}}$ \quad
$\vp_M(\ell(\g_6)) = (\Delta^2<3,5,6>)^{\uZ^{-2}_{45}\uZ^{-2}_{46}}$\\
{\bf Remark 10.1.4.} {\em 
Similarly to Remark 10.1.2, the changes are:\\
(a) $(\Delta^2<3,5,6>)^{\uZ^{-2}_{45}\uZ^{-2}_{46}}$ is replaced by $(\Delta^2<1,3,5,6>)^{\uZ^{-2}_{45}\uZ^{-2}_{46}} $ with the corresponding 
{\rm L.V.C.}:
$\vcenter{\hbox{\epsfbox{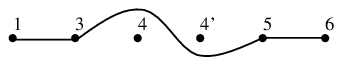}}}$ .\\
(b) $Z^2_{4'5}$ is replaced by $(\uZ^2_{4'i})_
{i = 1,5}$ with the following {\rm L.V.C.}:
 $\vcenter{\hbox{\epsfbox{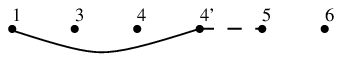}}}$ .\\
(c) $(\uZ^2_{45})^{Z^2_{34}}$ is replaced by
$(\uZ^2_{4i})_{i = 1,5}^{Z^2_{34}}$ with the following {\rm L.V.C.}:
 $\vcenter{\hbox{\epsfbox{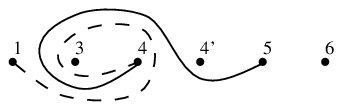}}}$.}

\noindent
\underline{Proof of Theorem 10.1:}

Until now we computed the braid monodromy of all singularities.  Each one of the intersection points of $L_4 \cap Q_2$ is replaced by 2 intersection points $(\subseteq Q_2 \cap Q_4)$ which are close to each other.  So $\pi^{-1} (M) \cap S^{(6)}_7 =
\{1,2,2',3,4,4',5,6\}$.

The changes are:\\
(I) Similar to (I) in the proof of Theorem 8.1.\\
(II) $\uZ^2_{2'4}$ is replaced by
$(\uZ^2_{2'i})_{i = 4,4'}$.  These braids correspond to the paths
$(\underline{z}_{2'i})_{i = 4,4'}$.  $(\uZ^2_{24})^{Z^2_{12}}$
is replaced by $(\uZ^2_{2i})_{i=4,4'}^{Z^2_{12}}$.
These braids correspond to the paths $(\underline{z}_{2i})_{i = 4,4'}^{Z^2_{12}}$.\\
(III) In a similar proof as in Theorem 8.1, we have to conjugate all braids from Proposition 10.1.3 and Remark 10.1.4 by $Z^2_{12}$.

According to the above changes, we present here the list of braids:


\begin{tabbing}
$\vcenter{\hbox{\epsfbox{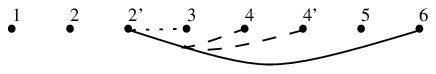}}}$vvvvvvv \= dd corresponding to
$ Z^4_{34})^2_{i = 2,5,6,6'}$ \kill
$\vcenter{\hbox{\epsfbox{v7058.ps}}}$ \> corresponding to
$ (\uZ^2_{2'i})_{i=3,4,4',6}$ \\ [.5cm]
$\vcenter{\hbox{\epsfbox{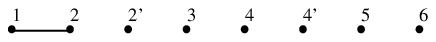}}}$ \> corresponding to
$Z^4_{12}$\\[.5cm]
$\vcenter{\hbox{\epsfbox{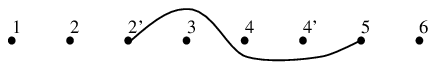}}}$ \> corresponding to
$(\uZ^4_{2'5})^{Z^2_{2'3}}$\\[.5cm]
$\vcenter{\hbox{\epsfbox{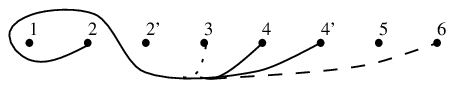}}}$ \> corresponding to
$(\uZ^2_{2i})_{i = 3,4,4',6}^{Z^2_{12}} $\\[.5cm]
$\vcenter{\hbox{\epsfbox{v7062.ps}}}$ \> corresponding to
$(Z_{22'})^{Z^2_{12}\bZ^2_{2'5}Z^2_{44',5}} $\\[.5cm]
$\vcenter{\hbox{\epsfbox{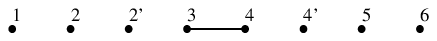}}}$ \> corresponding to
$Z_{34}^4$\\[.5cm]
$\vcenter{\hbox{\epsfbox{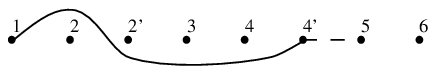}}}$ \> corresponding to
$(\uZ^2_{i4'})^{Z^{2}_{12}}_{i = 1,5}$\\[.5cm]
$\vcenter{\hbox{\epsfbox{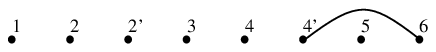}}}$ \> corresponding to $
\bZ^4_{4'6}$\\[.5cm]
$\vcenter{\hbox{\epsfbox{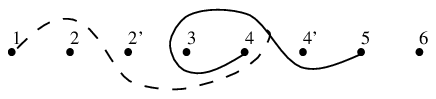}}}$ \> corresponding to
$(\uZ^2_{i4})^{Z^{2}_{34}Z^2_{12}}_{i=1,5}$ \\[.5cm]
$\vcenter{\hbox{\epsfbox{v7067.ps}}}$ \> corresponding to $
(Z_{44'})^{\bZ^2_{4'6}Z^2_{34}}$\\[.5cm]
$\vcenter{\hbox{\epsfbox{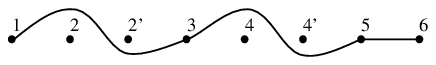}}}$ \> corresponding to
$(\Delta^2<1,3,5,6>)^{\uZ^{-2}_{45}\uZ^{-2}_{46}Z^2_{12}}$ \quad \quad \quad $\Box$
\end{tabbing}

Now we compute $\vp^{(5)}_7$.
Recall that $\vp^{(5)}_7 = \vp^{(i)}_7$ for $i = 1,2,3,4$,
  so in fact we compute $\vp^{(1)}_7$.

\noindent
{\bf Theorem 10.2.} {\em Let $\omega = L_1 \cup L_3 \cup L_5 \cup L_6$ in $S^{(6)}_7$. \\
(a) In a small neighbourhood $U_1$ of $\omega, T^{(1)} = S^{(1)}_7
\cap \ U_1$ resembles Figure 28, i.e., the singularities of
$T^{(1)}$ are 4 nodes,
4   tangency points and 2 branch points.\\
(b) The local braid monodromy of $S^{(1)}_7$ in that neighbourhood is presented by $F_1(F_1)^{\rho^{-1}}$, where:\\
 $F_1 = Z^2_{1'3}\cdot  Z^4_{3'5}  \cdot (\bZ^{2}_{1'3'})^{Z^2_{3'5}}
\cdot (\uZ^4_{1'5})^{Z^2_{1'3}} \cdot (Z_{56})^{\bZ^2_{1'5}Z^2_{3'5}}.\\
\rho = Z_{11'}   Z_{33'}$ and $(z_{56})^{\bZ^2_{1'5} Z^2_{3'5}} = $}
$\vcenter{\hbox{\epsfbox{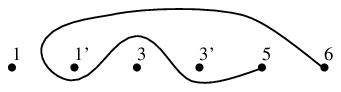}}}$ .

\noindent
\underline{Proof}:
(a) Let $\hat{\omega} = R_4 \cap R_5 \cap R_6 \cap R_7$ in $Z^{(6)}$.

For local analysis we can use holomorphic coordinates in a neighbourhood
$\hat{U}$ of $\hat{\omega}$ in $\C \P^n$.  This allows us to consider
 $\hat{U}$ as a neighbourhood of the origin in $\C^4$ with coordinates $X, Y, Z, T$
defined by the following system of equations
$
\hat{U} : \left\{ \begin{array}{l}
XT = 0 \\
YZ = 0 \ . \end{array} \right .$

By abuse of notation, $\hat{L}_3$  is called the $x$-axis (see Figure 27), and then: $R_4 = \{X=0, Y=0\}, \ R_5 = \{Y=0, T=0\}, \ R_6 = \{X=0, Z=0\}, \ R_7 = \{T=0, \ Z=0 \}$.
%
\begin{figure}[!h]
\epsfxsize=5cm 
\epsfysize=5cm 
\begin{minipage}{\textwidth}
\begin{center}
\epsfbox {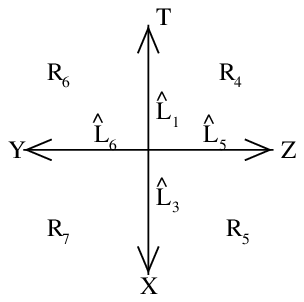}
\end{center}
\end{minipage}
\caption{}
\end{figure}

Consider the curve we obtain in Lemma 6, [MoTe8]. 
 As explained in [Am2, Subsection 2.4.1],  we can spin this curve to obtain Figure 28.

\begin{figure}[h]
\epsfxsize=11cm 
\epsfysize=9cm 
\begin{minipage}{\textwidth}
\begin{center}
\epsfbox{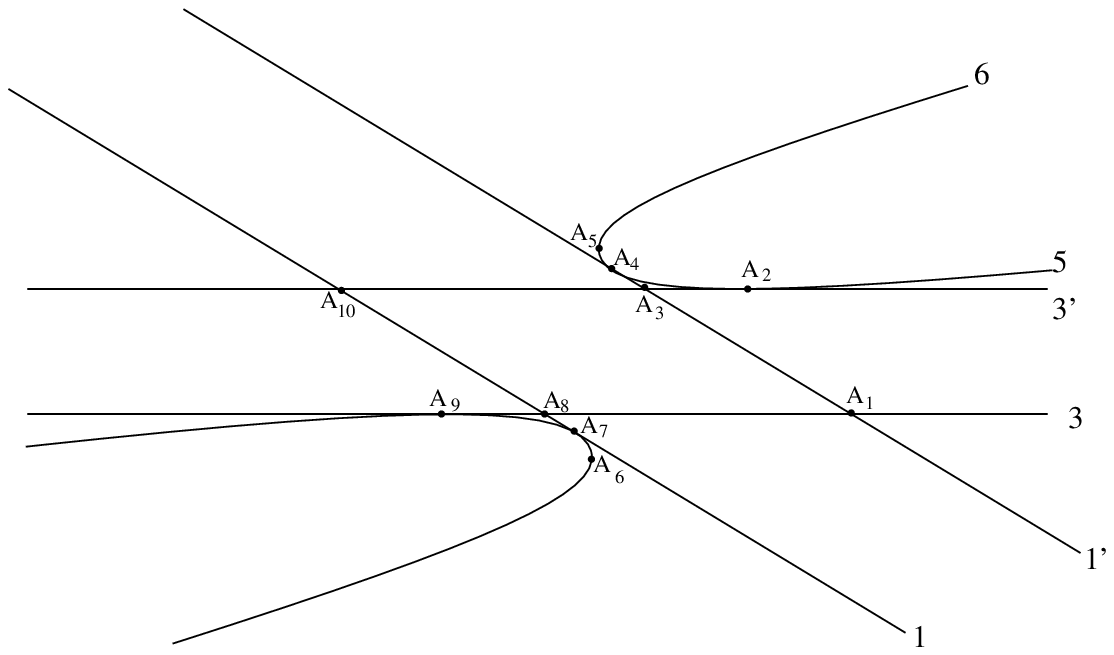}
\end{center}
\end{minipage}
\caption{}
\end{figure}

Let $M$ be a real point, $x << M \;\; \forall x \in N, K(M) = \{1,1',3,3',5,6\}$.\\
(b) We use the proof  of Theorem 8.2(b) here in a similar way.
$S^{(1)}_7 = L_1 \cup  L_{1'} \cup  L_3 \cup  L_{3'} \cup  h_{56}$,
see Figure 28.

Let $N, M, E, D, \{\ell(\g_j)\}^{10}_{j=1}$ be as in Theorem 8.2(b) and let $\{A_j\}^{10}_{j=1}$ be the singular points as shown in Figure 28.

Let $-P, P_1, P, T, T_1, T_2$ be as in Figure 18 and let
$x(L_1 \cap L_{1'}), x (L_3 \cap L_{3'})$, \st  
$M << x(L_1 \cap L_{1'}) < P_1 < x (L_3 \cap L_{3'}) < P \ , \ -P < x_{10}. \
K(M)=\{1,1',3,3',5,6\}$.

Choose a  diffeomorphism $\beta _M$ which satisfies:
$\beta_M(i) = i$ for $i = 1,3,5,6, \ \beta_M(1') = 2, \  \beta_M(3') = 4.$

We have the following table for the first five points:
\vspace{-1cm}
\begin{center}
\begin{tabular}{cccc} \\
$j$ & $\lambda_{x_j}$ & $\epsilon_{x_j}
$ & $\delta_{x_j}$ \\ \hline
1 & $<2,3>$ & 2 & $\Delta<2,3>$\\
2 & $<4,5>$ & 4 & $\Delta^2<4,5>$\\
3 & $<3,4>$ & 2 & $\Delta<3,4>$\\
4 & $<4,5>$ & 4 & $\Delta^2<4,5>$\\
5 & $<5,6>$ & 1 & $\Delta^{\frac{1}{2}}_{I_2 \R}<5>$\\
\end{tabular}
\end{center}

For computations we use the formulas from [Am2, 
Theorems 1.41, 1.44].

\medskip
\noindent
$(\xi_{x'_1}) \Psi_{\g'_1} = <2,3> \beta^{-1}_M = z_{1'3}$

\noindent
$<2,3>$ $\vcenter{\hbox{\epsfbox{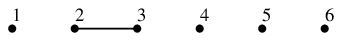}}}$ \ $\beta^{-1}_M$ \
$\vcenter{\hbox{\epsfbox{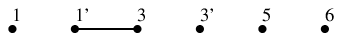}}}$ \ $\vp_M(\ell(\g_1)) = Z^2_{1'3}$ \\

\noindent
$(\xi_{x'_2}) \Psi_{\g'_2} = <4,5> \Delta<2,3>\beta^{-1}_M = z_{3'5}$

\medskip
\noindent
$<4,5> \Delta <2,3>$  $\vcenter{\hbox{\epsfbox{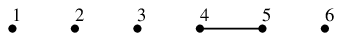}}}$ \
$\beta^{-1}_M$  $\vcenter{\hbox{\epsfbox{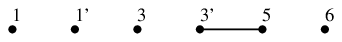}}}$ \ $\vp_M(\ell(\g_2)) = Z^4_{3'5}$ \\

\noindent
$(\xi_{x'_3}) \Psi_{\g'_3}
= <3,4>\Delta^2<4,5> \Delta<2,3>
\beta^{-1}_M =
\bar{z}_{1'3'}^{Z^2_{3'5}}
$

\noindent
$<3,4>$   $\vcenter{\hbox{\epsfbox{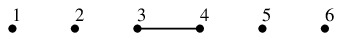}}}$ \ $\Delta^2<4,5>$ \
 $\vcenter{\hbox{\epsfbox{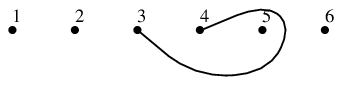}}}$ \ $\Delta<2,3>$  $\vcenter{\hbox{\epsfbox{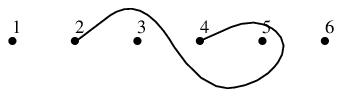}}}$ \\
\noindent
$\beta^{-1}_M$ \ $\vcenter{\hbox{\epsfbox{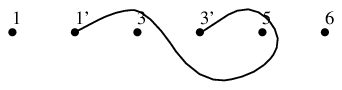}}}$ \
$\vp_M(\ell(\g_3)) = (\bZ^2_{1'3'})^{Z^2_{3'5}}$\\

\medskip

\noindent
$(\xi_{x'_4}) \Psi_{\g'_4}  = <4,5>   \Delta<3,4>\Delta^2<4,5>
\Delta<2,3> \beta^{-1}_M = \underline{ z}_{1'5}^{Z^2_{1'3}}$

\noindent
$<4,5>$ $\vcenter{\hbox{\epsfbox{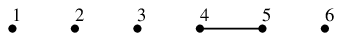}}}$ \  $\Delta<3,4>$ \
$\vcenter{\hbox{\epsfbox{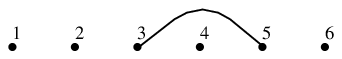}}}$ \
$\Delta^2<4,5>$ \ $\vcenter{\hbox{\epsfbox{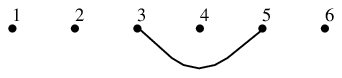}}}$ \\
$\Delta<2,3>$ \ $\vcenter{\hbox{\epsfbox{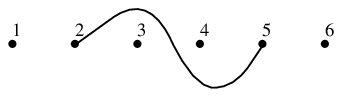}}}$ \
$\beta^{-1}_M$ $\vcenter{\hbox{\epsfbox{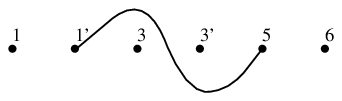}}}$ \
$\vp_M(\ell(\g_4)) = (\uZ^4_{1'5})^{Z^2_{1'3}}$\\

\noindent
$(\xi_{x'_5}) \Psi_{\g'_5}
 =  <5,6> \Delta^2<4,5>   \Delta<3,4>\Delta^2<4,5>
\Delta<2,3> \beta^{-1}_M =  z_{56}^{\bZ^2_{1'5}Z^2_{3'5}}$\\
\noindent
$<5,6>$ \ $\vcenter{\hbox{\epsfbox{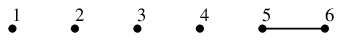}}}$ \ $\Delta^2<4,5>$ \
$\vcenter{\hbox{\epsfbox{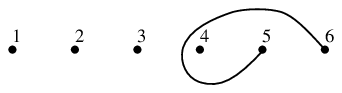}}}$ \
$\Delta<3,4>$   $\vcenter{\hbox{\epsfbox{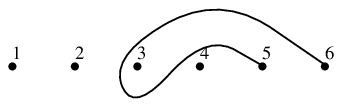}}}$ \\
$\Delta^2<4,5>$  \ $\vcenter{\hbox{\epsfbox{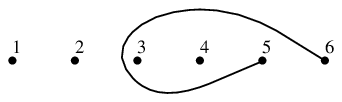}}}$ \
$\Delta<2,3>$  $\vcenter{\hbox{\epsfbox{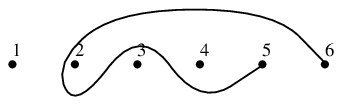}}}$ \
$\beta^{-1}_M$ \
$\vcenter{\hbox{\epsfbox{7v062.ps}}}$ \\
\noindent
$\vp_M(\ell(\g_5)) =   (Z_{56})^{\bZ^2_{1'5}Z^2_{3'5}}$

The sequence of the braids that we obtain is:\\
$Z^2_{1'3}, Z^4_{3'5}, (\bZ^2_{1'3'})^{Z^2_{3'5}}, (\uZ^4_{1'5})^{Z^2_{1'3}} ,
(Z_{56})^{\bZ^2_{1'5} Z^2_{3'5}}$ .
So we obtain $F_1$ as a factorized expression.

Now we want to compute $\beta^\vee_M (\vp_M(\ell(\g_j)))$  for $6 \le j \le 10$.
First, we have to compute $\beta^\vee_{-P} (\vp_{-P}(\ell(\tg_j)))$.
We have the following table:
\vspace{-1cm}
\begin{center}
\begin{tabular}{cccc} \\
$j$ & $\lambda_{x_j}$ & $\epsilon_{x_j}
$ & $\delta_{x_j}$ \\ \hline
10 & $<4,5>$ & 2 & $\Delta<4,5>$\\
9 & $<2,3>$ & 4 & $\Delta^2<2,3>$\\
8 & $<3,4>$ & 2 & $\Delta<3,4>$\\
7 & $<2,3>$ & 4 & $\Delta^2<2,3>$\\
6 & $<1,2>$ & 1 & $\Delta^{\frac{1}{2}}_{I_2 \R}<1>$\\
\end{tabular}
\end{center}

\medskip

Now we are computing $\beta^\vee_{-P} (\vp_{-P}(\ell(\tg_j )))$ for $j = 10, \cdots , 6$.\\
\noindent
L.V.C. $(\tg_{10}){\beta_{-P}} = z_{45}$\\
\noindent
$<4,5>$  \ $\vcenter{\hbox{\epsfbox{7v071.ps}}}$ \ $\beta^\vee_{-P}(\vp_{-P}(\ell(\tg_{10})))= Z^2_{45}$\\

\medskip
\noindent
L.V.C. $(\tg_{9}){\beta_{-P}} = <2,3> \Delta<4,5> = z_{23}$\\
\noindent
$<2,3>$ \ $\Delta<4,5>$ $\vcenter{\hbox{\epsfbox{7v063.ps}}}$ \ $\beta^\vee_{-P}(\vp_{-P}(\ell(\tg_{9}))) = Z^4_{23}$\\

\medskip
\noindent
L.V.C. $(\tg_{8}){\beta_{-P}} = <3,4> \Delta^2<2,3>\Delta<4,5> = \underline{
z}_{35}^{Z^2_{23}}$\\
\noindent
$<3,4>$  \ $\vcenter{\hbox{\epsfbox{7v067.ps}}}$ \ $\Delta^2<2,3>$ \
$\vcenter{\hbox{\epsfbox{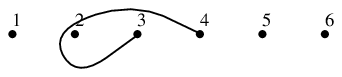}}}$ \
$\Delta<4,5>$ \ $\vcenter{\hbox{\epsfbox{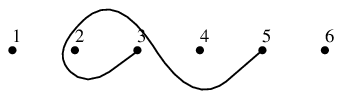}}}$ \\
\noindent
$\beta^\vee_{-P}( \vp_{-P}(\ell(\tg_{8}) ))= (\uZ^2_{35})^{Z^2_{23}}$\\

\medskip
\noindent
L.V.C.$(\tg_7)\beta_{-P} = <2,3>\Delta<3,4> \Delta^2<2,3>\Delta <4,5> =
 \underline{z}_{25}^{Z^2_{23}}$\\
\noindent
$<2,3>$  $\vcenter{\hbox{\epsfbox{7v063.ps}}}$ \
$\Delta<3,4>$ \ $\vcenter{\hbox{\epsfbox{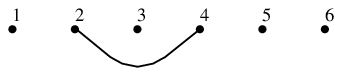}}}$ \
$\Delta^2<2,3>$   $\vcenter{\hbox{\epsfbox{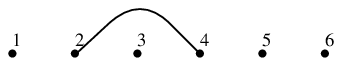}}}$ \\
\noindent
$\Delta<4,5>$ $\vcenter{\hbox{\epsfbox{7v074.ps}}}$ \
$\beta^\vee_{-P}(\vp_{-P}(\ell(\tg_{7}))) = (\uZ^4_{25})^{Z^2_{23}}$\\

\medskip
\noindent
L.V.C. $(\tg_6)\beta_{-P} = <1,2>\Delta^2<2,3>\Delta<3,4> \Delta^2<2,3>
\Delta<4,5> = z_{12}^{\uZ^2_{25} Z^2_{23}}$\\
\noindent
$<1,2>$ $\vcenter{\hbox{\epsfbox{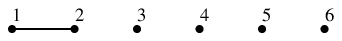}}}$ \ $\Delta^2<2,3>$ \
$\vcenter{\hbox{\epsfbox{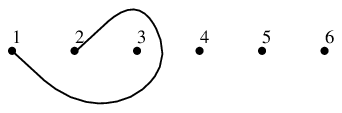}}}$ \
$\Delta<3,4>$ $\vcenter{\hbox{\epsfbox{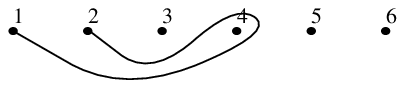}}}$ \\
\noindent
$\Delta^2<2,3>$ \ $\vcenter{\hbox{\epsfbox{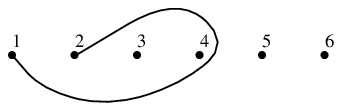}}}$ \
$\Delta<4,5>$ $\vcenter{\hbox{\epsfbox{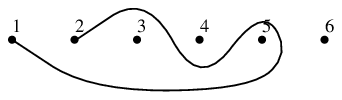}}}$ \\ $\beta^\vee_{-P}(\vp_{-P}(\ell(\tg_{6}))) = (Z_{12})^{\uZ^2_{25}Z^2_{23}}$\\

Apply $180^0$ rotation on L.V.C. $(\tg_j)\beta_{-P}$ to get L.V.C.$(\tg_jT)\beta_{P}$ for $j = 10, \cdots ,6$:\\
\noindent
L.V.C. $ (\tg_{10}T)\beta_P = \vcenter{\hbox{\epsfbox{7v063.ps}}}  = z_{23}$ ;
L.V.C. $ (\tg_{9}T)\beta_P = \vcenter{\hbox{\epsfbox{7v071.ps}}} = z_{45}$ \\
\noindent
L.V.C. $ (\tg_{8}T)\beta_P = \vcenter{\hbox{\epsfbox{7v069.ps}}} = \bar{z}_{24}^{Z^2_{45}}$  ;
L.V.C. $ (\tg_{7}T)\beta_P =  \vcenter{\hbox{\epsfbox{7v074.ps}}} = \underline{z}_{25}^{Z^2_{23}}$ \\
\noindent
L.V.C. $ (\tg_{6}T)\beta_P = \vcenter{\hbox{\epsfbox{7v080.ps}}} = z_{56}^{\bZ^2_{25}Z^2_{45}}$ .

Moreover $\beta ^\vee_P(\vp_P(\ell(\tg_jT))) =
\Delta < \ $L.V.C.$(\tg_jT)\beta_{P}>^{\epsilon_{x_j}}$ for $j = 10, \cdots ,6$.
Thus:\\
$\beta^\vee_{P}(\vp_{P}(\ell(\tg_{10}T)))= Z^2_{23}$ \ ; \
$\beta^\vee_{P}(\vp_{P}(\ell(\tg_{9}T)))= Z^4_{45}$\ ; \
$\beta^\vee_{P}(\vp_{P}(\ell(\tg_{8}T)))= (\bZ^2_{24})^{Z^2_{45}}$ \ ; \
$\beta^\vee_{P}(\vp_{P}(\ell(\tg_{7}T)))= (\uZ^4_{25})^{Z^2_{23}}$ \ ; \
$\beta^\vee_{P}(\vp_{P}(\ell(\tg_{6}T)))= (Z_{56})^{\bZ^2_{25}Z^2_{45}}$.

By Figure 28, L-pair $(x(L_1
\cap L_{1'})) = (1,1')$ and
 L-pair $(x(L_3
\cap L_{3'})) = (3,3')$.
$\rho = \Delta(1,1') \Delta(3,3')$ and:
$\vp_M(\ell(\g_{10})) = (Z^2_{1'3})^{\rho^{-1}}$ \ ; \
$\vp_M(\ell(\g_{9})) = (Z^4_{3'5})^{\rho^{-1}}$ \ ; \
$\vp_M(\ell(\g_{8})) = ((\bZ^2_{1'3'})^ {Z^2_{3'5}})^{\rho^{-1}}$\ ; \
$\vp_M(\ell(\g_{7})) = ((\uZ^4_{1'5})^{Z^2_{1'3}})^{\rho^{-1}}$ \ ; \
$\vp_M(\ell(\g_{6})) = ((Z_{56})^{\bZ^2_{1'5}Z^2_{3'5}})^{\rho^{-1}}$.

We obtain $(F_1)^{\rho^{-1}}$.

Therefore, the braid monodromy w.r.t. $E \times D$ is $\vp_M = F_1 (F_1)^{\rho^{-1}}$.
\hfill
$\Box$
\medskip
\noindent
{\bf Proposition 10.3.} {\em The local braid monodromy for $S^{(1)}_7$ around $V_7$ is obtained from formula $\vp^{(6)}_7$ in Theorem 10.1 by the following replacements:\\(i) Consider the following $( \ )^\ast$ as conjugations by the braids, induced from the following motions: \\
\noindent
$(\quad \quad)^{Z^2_{2',33'}} \ \vcenter{\hbox{\epsfbox{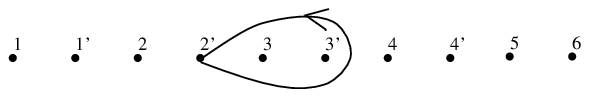}}}$ ;
 $(\quad \quad)^{Z^{2}_{11',2}} \ \vcenter{\hbox{\epsfbox{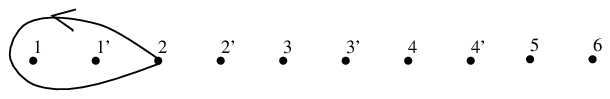}}}$\\
$(\quad \quad)^{Z^{2}_{33',4}} \
\vcenter{\hbox{\epsfbox{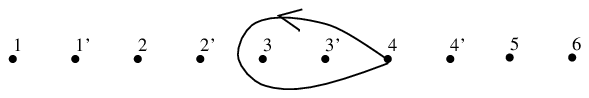}}}$

All the other conjugations do not change, since indices 1, 3 are not involved.

\noindent
(ii) $\Delta^2<1,3,5,6>$ is replaced by  $F_1(F_1)^{\rho^{-1}}$.\\
(iii) Each of the degree 4 factors in $\vp^{(6)}_7$ that involves index 1 or 3 is replaced by 3 cubes as in the third regeneration rule.\\
(iv) Each of the degree 2 factors in $\vp^{(6)}_7$ that involves index 1 or 3 is replaced by 2 degree 2 factors, where : 1 and 1' are replacing 1; 3 and 3' are replacing 3.

We call the formula that was obtained: (1)$_7$.}

\noindent
\underline{Proof}:  A similar proof as the proof of Proposition 8.3 but here changes are applied on indices 1 and 3.
\hfill
$\Box$

Now we compute $\vp^{(0)}_{7}$.

\noindent
{\bf Theorem 10.4.} {\em  In the notation of Theorem 10.2:  let $T^{(0)}$ be the curve
 obtained from $T^{(1)}$ in the regeneration $Z^{(1)} \leadsto Z^{(0)}$. \\
(i) Then the local braid monodromy of $T^{(0)}$ is $\hat{F}_1 (\hat{F}_2)$, where:\\
$\hat{F}_1 = Z^2_{1'3} \cdot \uZ^3_{3',55'}\cdot(\uZ^2_{1'3'})^{\uZ^2_{3',55'}\uZ^2_{1'3}}
 \cdot(\uZ^3_{1',55'})^{Z^2_{1'3}} \; \; \cdot Z_{56'} \cdot  Z_{5'6},$ \\
$\hat{F}_2 =  (Z^2_{1'3})^{\rho^{-1}} \cdot (\uZ^3_{3',55'})^{\rho^{-1}}
\cdot((\uZ^2_{1'3'})^{\uZ^2_{3',55'}\uZ^2_{1'3}})^{\rho^{-1}}
 \cdot((\uZ^3_{1',55'})^{Z^2_{1'3}})^{\rho^{-1}} \; \; \cdot
 Z_{56'}^{\rho^{-1}} \cdot  Z_{5'6}^{\rho^{-1}}.$

$m_1, m_2$ arise from the first regeneration rule applied on $\alpha ^{(1)}$ and $(\alpha ^{(1)})^{\rho^{-1}}$ respectively.
Recall that $\alpha = (Z_{56})^{\bZ^2_{1'5} Z^2_{3'5}}$ and $\rho = Z_{11'} Z_{33'}$.\\
\noindent
(ii) The singularities of the $x$-projection of \ $T^{(0)}$ are those arising from $T^{(1)}$ by the regeneration rules, namely: 4 nodes that exist in $T^{(1)},\ 3 \times 4$ cusps arising from 4 tangency points in $T^{(1)}, \
2 \times 2$ branch points for 2 branch points of $T^{(1)}$.\\
(iii) The braid monodromy of $T^{(0)}$ is $\hat{F}_1 (\hat{F}_2)$. \\
(iv) $\hat{F}_1 (\hat{F}_2) = \Delta^{-2}_8 Z^{-2}_{11'} Z^{-2}_{33'} Z^{-2}_{55'} Z^{-2}_{66'}$.}

\noindent
\underline{Proof}: A similar proof as the proof of Theorem 8.4 but replace $i = 1,2,3,5$ with $i = 1,3,5,6$.

\vspace{-.2cm}

\hfill
$\Box$

\def\hF{\hat{F}}
 The local braid monodromy of $S^{(0)}_7$ in $U_1$ is $\hat{F}_1 (\hat{F}_1)^{\rho^{-1}}$.  Applying regeneration rules on $F_1(F_1)^{\rho^{-1}}$, we obtain $\hat{F}_1 (\hat{F}_1)^{\rho^{-1}}$.\\
\noindent
{\bf Corollary 10.5.} {\em The paths corresponding to $\hF_1(\hF_1)^{\rho^{-1}}$ (without their conjugation). \\
(a) \underline{ The paths corresponding to the factors in $\hF_1$}:

\medskip\noindent $\vcenter{\hbox{\epsfbox{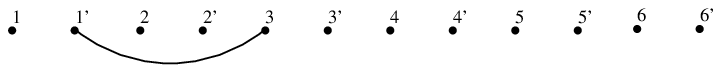}}}$ \quad corresponding to $Z^2_{1'3}$

\medskip

\noindent $\vcenter{\hbox{\epsfbox{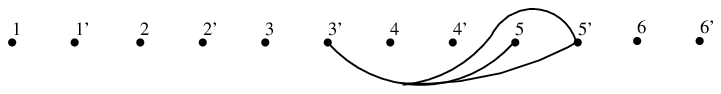}}}$ \quad  corresponding to
 $\uZ^3_{3',55'}$

\medskip

\noindent $\vcenter{\hbox{\epsfbox{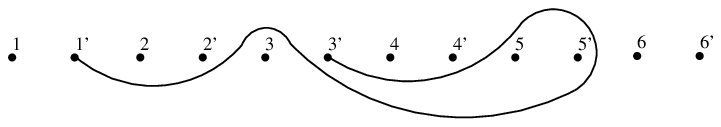}}}$ \quad  corresponding to $(\uZ^2_{1'3'})^{\uZ^2_{3',55'}\uZ^2_{1'3}}$

\medskip
\noindent $\vcenter{\hbox{\epsfbox{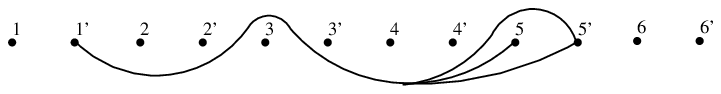}}}$ \quad  corresponding to $(\uZ^3_{1',55'})^{Z^2_{1'3}}$

\medskip

\noindent $\vcenter{\hbox{\epsfbox{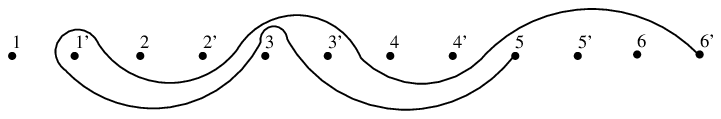}}}$ \quad corresponding to 
$Z_{56'}$
\medskip

\noindent $\vcenter{\hbox{\epsfbox{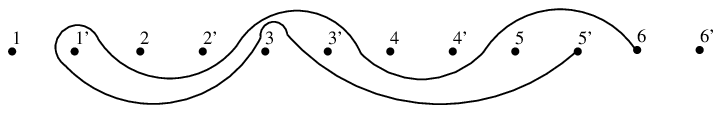}}}$ \quad corresponding to  
$Z_{5'6}$
\medskip

\noindent
(b) \underline{ The paths corresponding to $(\hF_1)^{\rho^{-1}}$ for
$\rho^{-1} = Z^{-1}_{33'} Z^{-1}_{11'}$} \ :

\medskip
\noindent $\vcenter{\hbox{\epsfbox{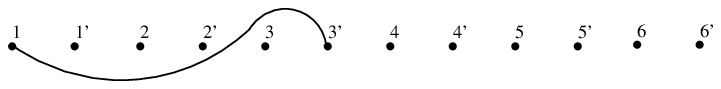}}}$ \quad  corresponding to $(Z^2_{1'3})^{\rho^{-1}}$

\medskip

\noindent $\vcenter{\hbox{\epsfbox{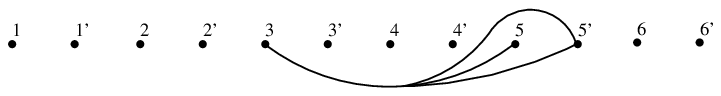}}}$ \quad  corresponding to $(\uZ^3_{3',55'})^{\rho^{-1}}$

\medskip

\noindent $\vcenter{\hbox{\epsfbox{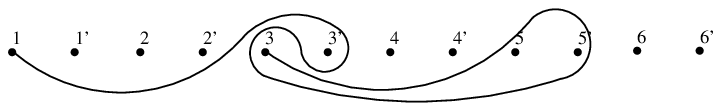}}}$ \quad corresponding to $\left((\uZ^2_{1'3'})^{\uZ^2_{3',55'}\uZ^2_{1'3}}\right)^{\rho^{-1}}$

\medskip

\noindent $\vcenter{\hbox{\epsfbox{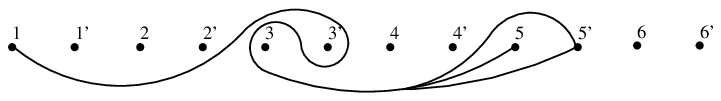}}}$ \quad corresponding to $\left((\uZ^3_{1',55'})^{Z^2_{1'3}}\right)^{\rho^{-1}}$

\medskip

\noindent $\vcenter{\hbox{\epsfbox{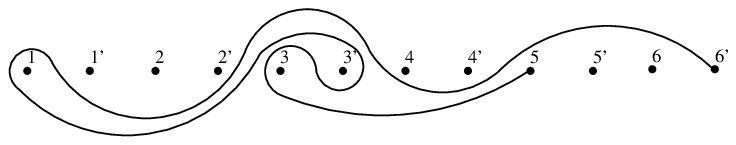}}}$ \quad corresponding to $Z_{56'}^{\rho^{-1}}$

\medskip

\noindent $\vcenter{\hbox{\epsfbox{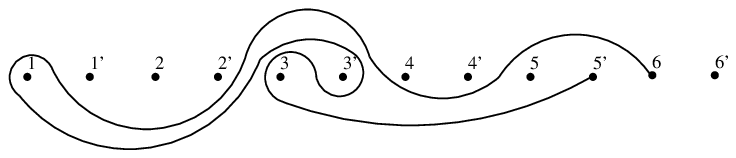}}}$ \quad  corresponding to  $Z_{5'6}^{\rho^{-1}}$}

\medskip

\noindent

\noindent
{\bf Theorem 10.6.} {\em  The local braid monodromy of $S^{(0)}_7$ around $V_7$, denoted by $H_{V_7}$, equals:\\
$H_{V_7} =
(\uZ^2_{2'i})_{i=33',4,4',66'} \cdot
Z^3_{11',2} \cdot
(\uZ^3_{2',55'})^{Z^2_{2',33'}} \cdot
(\uZ^{2}_{2i})^{Z_{11',2}}_{i = 33',4,4',66'}
 \cdot
(Z_{22'})^{Z^2_{11',2}\bZ^2_{2',55'}Z^2_{44',55'}} \cdot
Z^3_{33',4} \cdot
(\uZ^2_{i4'})^{Z^{2}_{11',2}}_{i=11',55'} \cdot
\bZ^3_{4',66'} \cdot
(\uZ^2_{i4})^{Z^2_{33',4}Z^2_{11',2}}_{i = 11',55'} \cdot
(Z_{44'})^{\bZ^2_{4',66'}Z^2_{33',4}}
\cdot
\left(\hF_1(\hF_1)^{\rho^{-1}}\right)^{\uZ^{-2}_{4,55'}
 \uZ^{-2}_{4,66'}Z^2_{11',2}}$ ,
where the paths corresponding to these braids are (the paths corresponding to
$\hF_1(\hF_1)^{\rho^{-1}}$ are above):

\begin{tabbing}
$\vcenter{\hbox{\epsfbox{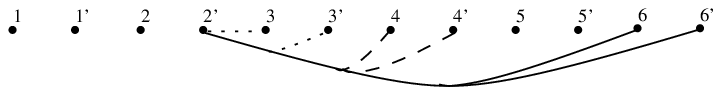}}}$vvvvvvv \= dd corresponding to
$ Z^4_{34})^2_{i = 2,5,6,6'}$ \kill
(1)\quad $\vcenter{\hbox{\epsfbox{v7069.ps}}}$ \> corresponding to $ (\uZ^2_{2'i})_{i=33',4,4',66'}$\\[.5cm]
(2)\quad $\vcenter{\hbox{\epsfbox{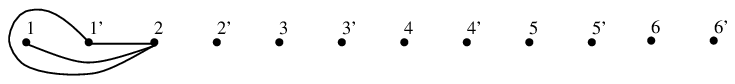}}}$ \>  corresponding to $Z^3_{11',2}   $\\[.5cm]
(3) \quad $\vcenter{\hbox{\epsfbox{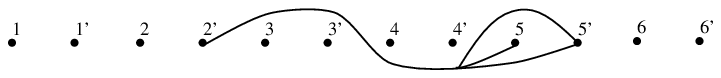}}}$ \> corresponding to $
(\uZ^3_{2',55'})^{Z^2_{2',33'}} $ \\[.5cm]
(4)\quad $\vcenter{\hbox{\epsfbox{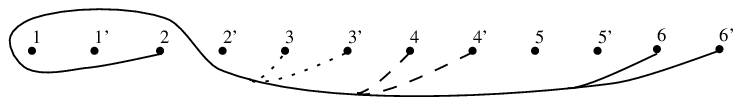}}}$ \>  corresponding to $
(\uZ^2_{2i})_{i=33',4,4',66'}^{Z^2_{11',2}}   $ \\[.5cm]
(5)\quad $\vcenter{\hbox{\epsfbox{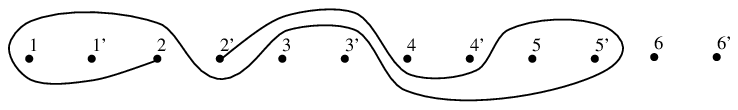}}}$ \>   corresponding to $
 (Z_{22'})^{Z^2_{11',2}\bZ^2_{2',55'}Z^{2}_{44',55'}} $\\[.5cm]
(6)\quad $\vcenter{\hbox{\epsfbox{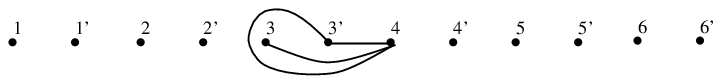}}}$ \>   corresponding to $
Z^{3}_{33',4}$\\[.5cm]
(7)\quad $\vcenter{\hbox{\epsfbox{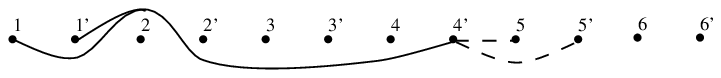}}}$ \>   corresponding to $
 (\uZ^2_{i4'})^{ Z^2_{11',2}}_{i= 11',55'}$\\[.5cm]
(8)\quad $\vcenter{\hbox{\epsfbox{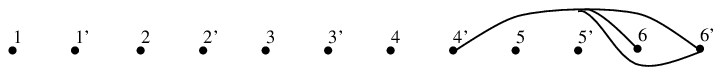}}}$ \>  corresponding to $ \bZ^3_{4',66'} $ \\[.5cm]
(9)\quad $\vcenter{\hbox{\epsfbox{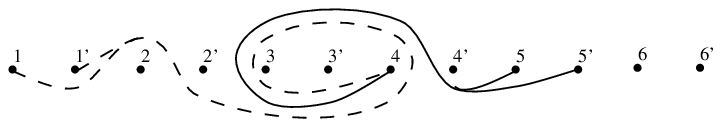}}}$ \>  corresponding to $
(\uZ^2_{i4})^{Z^2_{33',4}Z^2_{11',2}}_{i = 11',55'} $\\[.5cm]
(10)\quad $\vcenter{\hbox{\epsfbox{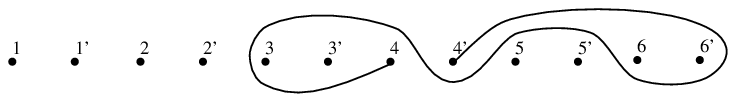}}}$ \>  corresponding to $ (Z_{44'})^{\bZ^2_{4',66'}Z^2_{33',4}}$
\end{tabbing}

}

\noindent
\underline{Proof}: A similar proof as the proof of Theorem 8.6 but changes are applied on (1)$_7$, to obtain a formula for a local braid monodromy of
$S^{(0)}_7$.  Moreover, all changes were applied on indices 5,6 to obtain
$H_{V_7}$.
\hfill
$\Box$

\noindent
{\bf Invariance Property 10.7.}  {\em $H_{V_7}$ is invariant under
$(Z_{22'} Z_{44'})^p (Z_{11'}Z_{33'})^q (Z_{55'} Z_{66'})^r \quad \forall p,q,r\in \Z$.}

\noindent
\underline{Case 1}: $p = q = r$. \\
A similar proof as Case 1 in Invariance Property 8.7. \\
\underline{Case 2}: $p = 0$.\\
Denote: $\epsilon = (Z_{11'} Z_{33'})^q (Z_{55'}Z_{66'})^r$.
We want to prove that  $H_{V_7}$ is invariant under $\epsilon $.\\
\underline{Step 1}:  Product of the factors outside of $\hF_1(\hF_1)^{\rho^{-1}}$.
\\
$  (Z_{22'})^{Z^2_{11',2}\bZ^2_{2',55'}Z^{2}_{44',55'}}    $
 and
$ (Z_{44'})^{\bZ^2_{4',66'}Z^{2}_{33',4}}$
commute with $\epsilon $. \\
The degree 3 factors are invariant under $\epsilon $ by Invariance Rule III.
All conjugations are invariant under $\epsilon$ too.
The
degree 2 factors are of the forms $Z^2_{\alpha,\beta \beta ' }$ where $\alpha
= 2,2',4,4'$ and $Z^2_{\alpha \alpha ', \beta }$ where $\beta = 4,4'$.
 So   they are invariant under $\epsilon $ by Invariance Rule II.\\
\noindent
\underline{Step 2}:   $\hF_1(\hF_1)^{\rho^{-1}}$.\\
A similar proof as step 2 in Invariance Property 8.7, but for:\\
Case 2.1: $q = 0 \ ; \ \epsilon = (Z_{55'}Z_{66'})^r$.\\
Case 2.2: $r = 0, \ q = 1 \ ; \ \epsilon = Z_{11'} Z_{33'}$ . \\
Case 2.3: $q = 2q' \ ; \ \epsilon =   (Z_{11'} Z_{33'})^{2q'} (Z_{55'}Z_{66'})^r$. \\
Case 2.4: $q = 2q' + 1 \ ; \ \epsilon = (Z_{11'} Z_{33'})^{2q'+1} (Z_{55'}Z_{66'})^r$. \\
\noindent
\underline{Case 3}:  $p,q,r$ arbitrary ; \ $\epsilon = (Z_{22'}Z_{44'})^p(Z_{11'}Z_{33'})^q(Z_{55'}Z_{66'})^r$.\\
By Case 1:$H_{V_7}$ is invariant under $(Z_{22'}Z_{44'})^p(Z_{11'}Z_{33'})^p(Z_{55'}Z_{66'})^p$.By Case 2: $H_{V_7}$ is invariant under
$ (Z_{11'} Z_{33'})^{q-p} (Z_{55'}Z_{66'})^{r-p}$.
By Invariance Remark (v),  $H_{V_7}$ is invariant under $\epsilon $.
\hfill
$\Box$

\noindent
{\bf Theorem 10.8.}  {\em Consider $H_{V_7}$ from Theorem 10.6.  The paths corresponding to the factors in $H_{V_7}$ , are shown below considering the Invariance Property 10.7.
Moreover, below to some paths, their complex conjugates appear too.}

\noindent
\underline{Remark}:  By abuse of notation, the simple braids are denoted by
$\underline{z}_{mk}$ and  the more complicated paths are denoted by $\tilde{z}_{mk}$.

\noindent
\underline{Proof:}  By the Invariance Property 10.7,  $H_{V_7}$ is invariant under $\epsilon =(\rho_2 \rho_4)^p (\rho_1 \rho_3)^q (\rho_5 \rho_6)^r \quad \\
 \forall p, q, r \in \Z$.
There are two possible applications for presenting the braids: (a) if $m,k$ are two indices in two different parts of $\epsilon $, then the braid is
$\rho^j_k \rho^i_m z_{mk} \rho^{-i}_m \rho^{-j}_k, \ i, j \in \Z, \ i \neq j$.
(b) If $m,k$ are in the same part of $\epsilon $, then the braid is $\rho^i_k \rho^i_m z_{mk} \rho^{-i}_m \rho^{-i}_k$, i.e.: when conjugating $z_{mk}$ by $\rho^i_m$ for some $i$, we must conjugate $z_{mk}$ also by $\rho^i_k$ for the same $i$.

We apply Complex Conjugation (8),(10).
 Their complex conjugates
appear below them.

\begin{center}
\begin{tabular}[hb]{|l|p{3in}|l|c|}\hline
& The paths corresponding & The braids & The exponent \\
[-.2cm] & to the braids and their  & &  (according to  \\
[-.2cm] & complex conjugates & & singularity type)\\  \hline (1) &
$\vcenter{\hbox{\epsfbox{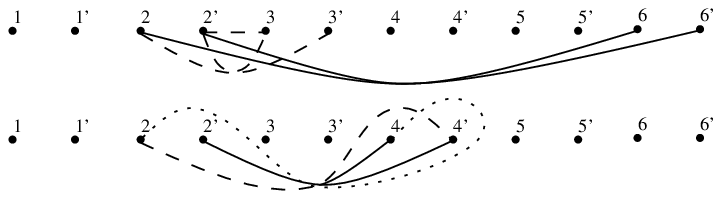}}}$  &
$\begin{array}{lll}\rho^j_m \rho^i_2
\underline{ z}_{2'm} \rho^{-i}_2 \rho^{-j}_m \\[-.2cm]
m= 3,6 \\ [-.2cm]
\rho^i_4 \rho^i_2
\underline{ z}_{2'4} \rho^{-i}_2 \rho^{-i}_4\end{array}$ & 2 \\
    \hline
(2) & $\vcenter{\hbox{\epsfbox{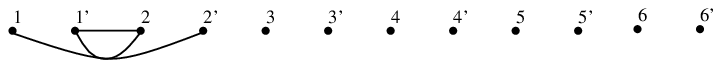}}}$ & $\rho^j_2 \rho^i_1 \underline{ z}_{12} \rho^{-i}_1 \rho^{-j}_2$ & 3  \\
   \hline
(3) & $\vcenter{\hbox{\epsfbox{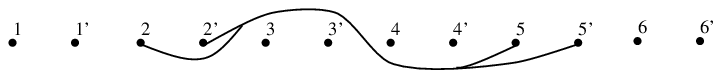}}}$ & $ \rho^j_5 \rho^i_2 \tilde{z}_{2'5} \rho^{-i}_2 \rho^{-j}_5$ & 3
\\
   \hline
(4) & $\vcenter{\hbox{\epsfbox{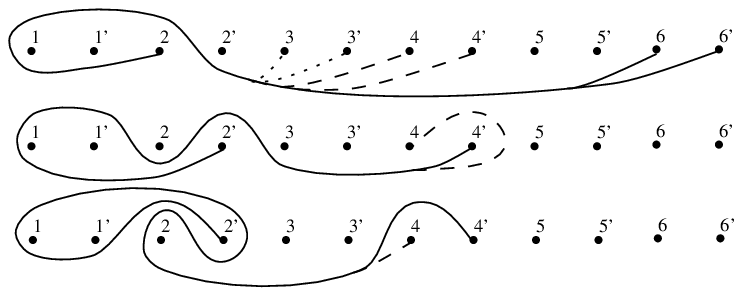}}} $ &
$\begin{array}{lll}\rho^j_m \rho^i_2 \tilde{z}_{2m} \rho^{-i}_2 \rho^{-j}_m \\
  [-.2cm]
   m = 3,6 \\  [-.2cm]
\rho^i_4 \rho^i_2 \tilde{z}_{24} \rho^{-i}_2 \rho^{-i}_4\end{array}$ & 2
\\ \hline
(5) & $\vcenter{\hbox{\epsfbox{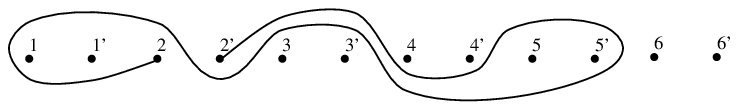}}}$  & $\tilde{z}_{22'} $ & 1\\
    \hline
(6) & $\vcenter{\hbox{\epsfbox{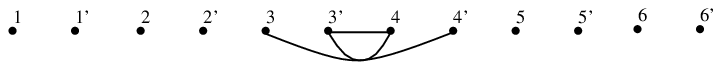}}}$  & $\rho^j_4 \rho^i_3 \underline{ z}_{34} \rho^{-i}_3 \rho^{-j}_4    $ & 3  \\
   \hline
(7) & $\vcenter{\hbox{\epsfbox{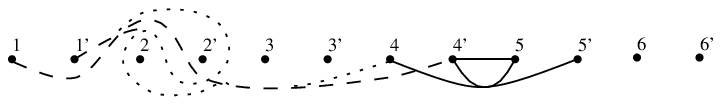}}}$  & $
\begin{array}{ll}\rho^j_m \rho^i_4 z_{4'm}\rho_4^{-i} \rho_m^{-j} \\
  [-.2cm]
m = 1,5\end{array}$ & 2 \\  \hline
(8) & $\vcenter{\hbox{\epsfbox{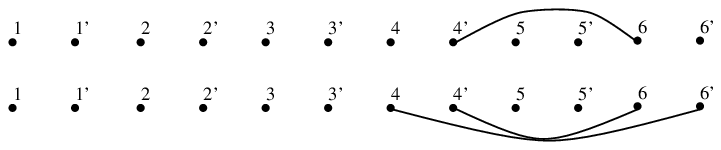}}}$  & $\rho^j_6 \rho^i_4 \underline{ z}_{4'6} \rho^{-i}_4 \rho^{-j}_6$ & 3 \\
    \hline
(9) & $\vcenter{\hbox{\epsfbox{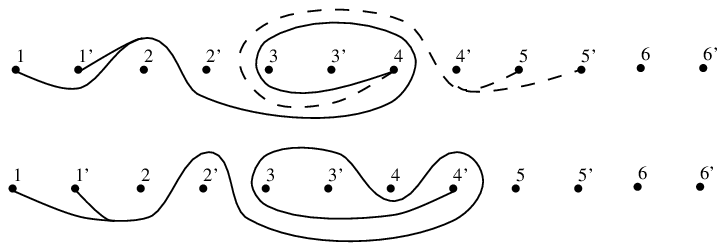}}}$  & $ \begin{array}{ll}\rho^j_4 \rho^i_m \tilde{z}_{m4} \rho^{-i}_m \rho^{-j}_4 \\
  [-.2cm]
  m = 1,5\end{array}$ & 2 \\  \hline
(10) & $\vcenter{\hbox{\epsfbox{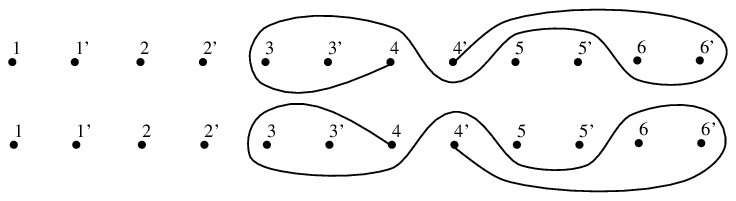}}}$ & $ \tilde{z}_{44'}$ & 1 \\
   \hline
\end{tabular}
\end{center}
\medskip

Consider now $\hF_1 (\hF_1)^{\rho^{-1}} (\rho^{-1} = \rho^{-1}_3 \rho_1^{-1})$ in $H_{V_7}$.
In a similar proof of Theorem 8.8 we obtain the following table:
\begin{center}
\begin{tabular}[H]{|l|c|c|}\hline
$\vcenter{\hbox{\epsfbox{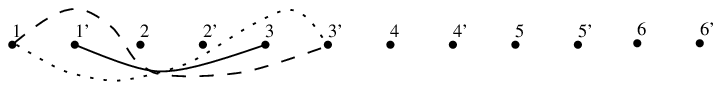}}}$ & $\rho^i_3 \rho^i_1 \underline{ z}_{1'3} \rho^{-i}_1 \rho^{-i}_3$ & 2\\
    \hline
$\vcenter{\hbox{\epsfbox{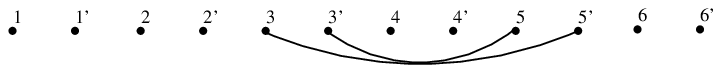}}}$  & $\rho^j_5 \rho^i_3  \underline{z}_{3'5} \rho^{-i}_3 \rho^{-j}_5$ & 3\\
    \hline
$\vcenter{\hbox{\epsfbox{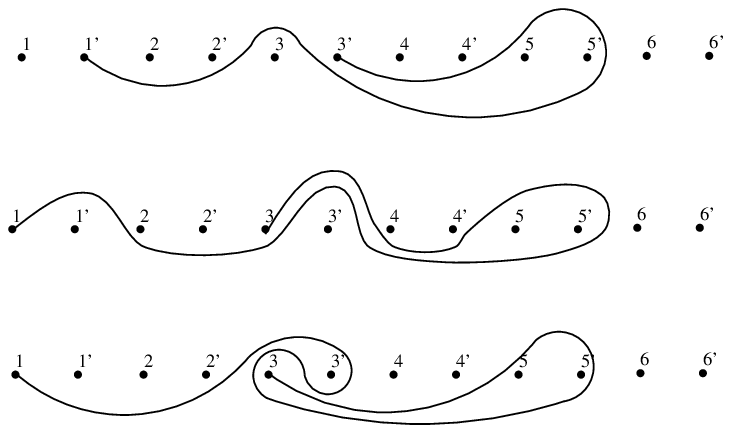}}}$ & $\rho^i_3 \rho^i_1 \tilde{z}_{1'3'} \rho^{-i}_1 \rho^{-i}_3$ & 2\\
  \hline
$\vcenter{\hbox{\epsfbox{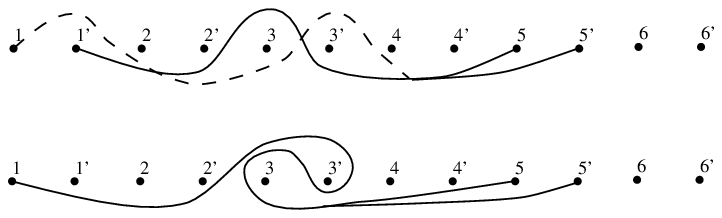}}}$ & $\rho^j_5 \rho^i_1 \tilde{z}_{1'5} \rho^{-i}_1 \rho^{-j}_5$ & 3\\
  \hline
$\vcenter{\hbox{\epsfbox{7v096.ps}}}$ & $\alpha _1(m_1)$ & 1\\
    \hline
$\vcenter{\hbox{\epsfbox{7v097.ps}}}$  & $\alpha _2(m_1)$ & 1\\
     \hline
$\vcenter{\hbox{\epsfbox{7v098.ps}}}$  & $(\alpha _1(m_1))^{\rho^{-1}}$ & 1 \\
  \hline
$\vcenter{\hbox{\epsfbox{7v099.ps}}}$  & $(\alpha  _2(m_1))^{\rho^{-1}}$ & 1 \\   \hline
\end{tabular}
\end{center}
\hfill
$\Box$

\section{The regenerated factors $C_i$}

We computed  $\tC_i$ for $i = 1,\cdots, 9$ in Section \ref{ssec:232}.  
Recall that: $\tC_1 = \prod\limits_{t=1,2,4,6,13,22} D_t \ , \ \tC_2 = \prod\limits_{t=3,7,9,14,17} D_t \ , 
\tC_3 = \prod\limits_{t = 5,10,18,23} D_t \ , \ \tC_4 = \prod\limits_{t=8,11,15,19} D_t \ ,\ \tC_5 = \prod\limits_{t=12,20,24} D_t \ ,  \tC_6 = \prod\limits_{t=16,25} D_t \ ,  \ \tC_7 = \prod\limits_{t=21,26} D_t \ , 
\tC_8 = D_{27} \ , \ \tC_9 = Id$.

$D_t ,  1 \leq t \leq 27$ were computed in  Section \ref{ssec:232}. 
We use the Complex Conjugation to obtain the following results 
(denoted as above):\\
\\
$D_1 = D_2 = D_3 = Id \ , \ D_4 = Z^2_{33',44'} \ , \
D_5 = \stackrel{\scriptstyle (4)(4')}{\uZ^2}_{\hspace{-.3cm}11',55'} \ ,  \
 \ D_6 = \uZ^2_{33',66'} \cdot \ Z^2_{55',66'}$ ,\\
$D_7 = \prod\limits_{i=2,4} \stackrel{\scriptstyle (6)(6')}{\uZ^2}_{\hspace{-.3cm}ii',77'} \cdot \ \uZ^2_{55',77'} \ ,
\ D_8 = \prod\limits^3_{i=1} \stackrel{\scriptstyle (6)-(7')}{\uZ^2}_{\hspace{-.3cm}ii',88'} \ , \   D_9 = \prod\limits_{i=2, 4-6,8} \uZ^2_{ii',99'} \ , \ \\
D_{10} = \prod\limits_{i=1,4,6,7} \stackrel{\scriptstyle (9)(9')}
{\uZ^2}_{\hspace{-.3cm}ii',10 \; 10'} \cdot \ \uZ^2_{88',10 \; 10'} \ ,
\ D_{11} = \prod\limits_{i=1-3, 6,7} \stackrel{\scriptstyle (9)-(10')}{\uZ^2}_{\hspace{-.4cm}ii',11 \; 11'} \ ,\\
D_{12} = \prod\limits^5_{i=1} \stackrel{\scriptstyle (9)-(11')}
{\uZ^2}_{\hspace{-.4cm}ii',12 \; 12'} \ , \
D_{13} = \prod\limits^{12}_{\stackrel{i=3}{i \neq 4,6}}   \uZ^2_{ii', 13 \; 13'} \ , \
 D_{14} = \prod\limits^{11}_{\stackrel{i=2}{i \neq 3,7,9}}
\stackrel{\scriptstyle (13)(13')}{\uZ^2}_{\hspace{-.4cm}ii',14 \; 14'} \cdot \
\uZ^2_{12 \; 12', 14 \; 14'} \ , \\
D_{15} = \prod\limits^{10}_{\stackrel{i=1}{i \neq 4,5,8}} \stackrel{\scriptstyle (13)-(14')}
{\uZ^2}_{\hspace{-.4cm}ii',15 \; 15'} \cdot \
\uZ^2_{12 \; 12' ,\; 15 \; 15'} \ , \
D_{16} = \prod\limits^{8}_{i=1} \stackrel{\scriptstyle (13)-(15')}
{\uZ^2}_{\hspace{-.4cm}ii',16 \; 16'} \ , \
D_{17} = \prod\limits^{16}_{\stackrel{i=2}{i \neq 3,7,9,14}}
\uZ^2_{ii',17 \; 17'} \ , \ \\
D_{18} = \prod\limits^{15}_{\stackrel{i=1}{i \neq 2,3,5,10}}
\stackrel{\scriptstyle (17)(17')}
{\uZ^2}_{\hspace{-.4cm}ii',18 \; 18'} \cdot \ \uZ^2_{16 \; 16', \; 18 \; 18'}
\ , \
D_{19} = \prod\limits^{14}_{\stackrel{i=1}{i \neq 4,5,8,11}}
\stackrel{\scriptstyle (17)-(18')}
{\uZ^2}_{\hspace{-.4cm}ii',19 \; 19'} \cdot \ \uZ^2_{16 \; 16', \; 19 \; 19'}
 \ , \ \\
D_{20} = \prod\limits^{15}_{\stackrel{i=1}{i \neq 6-8,12}}
\stackrel{\scriptstyle (17)-(19')}
{\uZ^2}_{\hspace{-.4cm}ii',20 \; 20'} \cdot \ \uZ^2_{16 \; 16', \; 20 \; 20'}
 \ , \
D_{21} = \prod\limits^{12}_{i=1}
\stackrel{\scriptstyle (17)-(20')}
{\uZ^2}_{\hspace{-.3cm}ii',21 \; 21'} \ , \
D_{22} = \prod\limits^{21}_{\stackrel{i=3}{i \neq 4,6,13}}
\uZ^2_{ii' , 22 \; 22'}\ , \\$
$D_{23} = \prod\limits^{20}_{\stackrel{i=1}{i \neq 2,3,5,10,18}}
\stackrel{\scriptstyle (22)(22')}
{\uZ^2}_{\hspace{-.4cm}ii',23 \; 23'} \cdot  \uZ^2_{21 \; 21', 23 \; 23'} \ , \
D_{24} = \prod\limits^{19}_{\stackrel{i=1}{i \neq 6-8, 12}}
\stackrel{\scriptstyle (22)-(23')}
{\uZ^2}_{\hspace{-.4cm}ii',24 \; 24'} \cdot  \ \uZ^2_{21 \; 21', 24 \; 24'} \ ,\\
D_{25} = \prod\limits^{20}_{\stackrel{i=1}{i \neq 9-12,16}}
\stackrel{\scriptstyle (22)-(24')}
{\uZ^2}_{\hspace{-.4cm}ii',25 \; 25'} \cdot \ \uZ^2_{21 \; 21', 25 \; 25'} \ , \
D_{26} = \prod\limits^{19}_{\stackrel{i=1}{i \neq 13-16}}
\stackrel{\scriptstyle (22)-(25')}
{\uZ^2}_{\hspace{-.4cm}ii',26 \; 26'} \cdot  \
\stackrel{\scriptstyle (21)(21')}
{\uZ^2}_{\hspace{-.4cm}20 \; 20', 26 \; 26'} \ , \  \\
D_{27} = \prod\limits^{16}_{i=1}
\stackrel{\scriptstyle (22)-(26')}
{\uZ^2}_{\hspace{-.4cm}ii',27 \; 27'} $.

During the regeneration, each $\tilde{C}_i$ is regenerated to $C_i$, $1 \leq i \leq 9$.
Each $C_i$ is now a product of the certain regenerated $D_t$ (as shown for $\tilde{C}_i$).

\section{Results}

In this paper we computed the braids $C_i$. We also computed the local 
braid monodromies $\varphi^{(0)}_i = H_{V_i}$ for $i=1,4,7$ in Sections 
\ref{sec:31}, \ref{sec:34}, \ref{sec:37}. The other resulting braid monodromies 
are shown in [AmTe2].

We get the regenerated braid monodromy factorization 
$\Delta^2_{54} = \prodl^9_{i=1} C_i H_{V_i}$.

In order to compute the fundamental group of the Galois cover of $T \times T$ w.r.t. a generic projection to $\C\P^2$, we need to compute first the fundamental group of the complement of $S$ in $\C^2$.

For that we have to apply the van Kampen Theorem 
on $\Delta^2_{54} = \prodl^9_{i=1} C_i H_{V_i}$. 
This is done in [AmTe2].

\section{Notations}\label{not}
\begin{list}{}{}
\item $(A)_B = B^{-1} AB = A^B.$
\item $X$ an algebraic surface, $X \subseteq \C \P^n$.
\item $X_0$ a degenerated object of a surface $X \ , \ X_0 \subset \C \P^N$.
\item $S$ an algebraic curve defined over $\R, \ S \subset \C^2$.
\item$E(\mbox{ resp}.D)$ be a closed disk on the $x$-axis (resp. y-axis) with the center on the real part of the $x$-axis (resp. y-axis), \st
$\{$singularities of $\pi_1\} \subseteq E \ast (D-\partial D)$.
\item $\pi  : S \rightarrow E$.
\item  $K(x) = \pi^{-1}(x)$.
\item $N = \{x \in E \suchthat \#K(x) < n \}$.
\item $u$ real number \st  $x << u \;\;\;  \forall  x \in N$.
\item $\C_u = \pi^{-1}(u)$.
\item $\rho_j = Z_{jj'}$.
\item $\vp$ = the braid monodromy of an algebraic curve $S$ in $M$.
\item $B_p[D,K]$ = the braid group.
\item $\Delta^2_p = (H_1 \cdots H_{p-1})^p$.
\item $\pi_1(\C^2 - S, M)$ = the fundamental group of a complement of a
branch curve $S$.
\item $\tP_1 = \frac{\pi_1(\C^2 - S, M)}{\la \G^2_j, \G^2_{j'}\ra}$ \ .
\item $T$ = complex torus.
\item $\underline{z}_{ij}$ = a path from $q_i \mbox{ to } q_j$ below the real line.
\item $\bar{z}_{ij}$ = a path from $q_i \mbox{ to } q_j$ above the real line.
\end{list}

The corresponding halftwists are:
$
H(\underline{ z}_{ij}) = \underline{Z}_{ij} \; ; \; H(\bar{z}_{ij}) = \bar{Z}_{ij}.$

\end{document}